\documentclass[10pt,leqno]{article}

\usepackage[francais]{babel}
\usepackage{amsmath,amssymb,amscd}
\usepackage{a4wide}
\usepackage[latin1]{inputenc}
\usepackage[all,dvips]{xy}

\newenvironment{paragr}[1][]{\refstepcounter{subsection} \noindent \textbf{\thesubsection . \ #1}}{\medskip}

\newenvironment{theoreme}{ \medskip\refstepcounter{theo}  \noindent\textbf{Th\'eor\`eme \thetheo}. ---\em}{\em \medskip}
\newenvironment{proposition}{\medskip\refstepcounter{theo}   \noindent\textbf{Proposition \thetheo}. ---\em}{\em\medskip}
\newenvironment{corollaire}{\medskip\refstepcounter{theo}  \noindent\textbf{Corollaire \thetheo}. ---\em}{\em\medskip}
\newenvironment{remarque}{\medskip \noindent \textbf{Remarque}. --- }{}
\newenvironment{remarques}{\medskip \noindent \textbf{Remarques}. --- }{}
\newenvironment{lemme}{\medskip\refstepcounter{theo}   \noindent\textbf{Lemme \thetheo}. ---\em}{\em\medskip}

\newenvironment{definition}{\medskip\refstepcounter{theo}  \noindent\textbf{D\'efinition \thetheo}. ---}{\medskip}

\newenvironment{preuve}[1][]{\noindent \textbf{Démonstration.} #1 --- }{\hfill
  \ensuremath{\square} \medskip}

\DeclareMathOperator{\reg}{reg}

\DeclareMathOperator{\vol}{vol}
\DeclareMathOperator{\rang}{rang}
\DeclareMathOperator{\Ad}{Ad}

\DeclareMathOperator{\ad}{ad}

\DeclareMathOperator{\der}{der}
\DeclareMathOperator{\scnx}{sc}
\DeclareMathOperator{\SCNX}{SC}
\DeclareMathOperator{\Norm}{Norm}

\DeclareMathOperator{\el}{ell}

\DeclareMathOperator{\Aut}{Aut}
\DeclareMathOperator{\Gal}{Gal}
\DeclareMathOperator{\Hom}{Hom}

\DeclareMathOperator{\Int}{Int}

\DeclareMathOperator{\stab}{stab}

\DeclareMathOperator{\Ker}{Ker}

\DeclareMathOperator{\val}{val}

\DeclareMathOperator{\Spec}{Spec}
\DeclareMathOperator{\Spf}{Spf}

\DeclareMathOperator{\car}{\mathfrak{car}}

\DeclareMathOperator{\card}{card}

\newcommand{\ZZ}{\mathbb{Z}}
\newcommand{\Gm}{\mathbb{G}_m}
\newcommand{\Gmk}{\mathbb{G}_{m,k}}

\newcommand{\NN}{\mathbb{N}}
\newcommand{\RR}{\mathbb{R}}
\newcommand{\AAA}{\mathbb{A}}

\newcommand{\FF}{\mathbb{F}}
\newcommand{\Fq}{\mathbb{F}_q}

\newcommand{\oc}{\mathcal{O}}

\newcommand{\tc}{\mathcal{T}}

\newcommand{\vc}{\mathcal{V}}

\newcommand{\ec}{\mathcal{E}}

\newcommand{\mc}{\mathcal{M}}

\newcommand{\lc}{\mathcal{L}}

\newcommand{\fc}{\mathcal{F}}
\newcommand{\pc}{\mathcal{P}}

\newcommand{\Ac}{\mathcal{A}}

\newcommand{\Qc}{\mathcal{Q}}
\newcommand{\xc}{\mathcal{X}}

\newcommand{\ggo}{\mathfrak{g}}

\newcommand{\mgo}{\mathfrak{m}}
\newcommand{\ngo}{\mathfrak{n}}
\newcommand{\ago}{\mathfrak{a}}
\newcommand{\cgo}{\mathfrak{c}}
\newcommand{\pgo}{\mathfrak{p}}
\newcommand{\qgo}{\mathfrak{q}}

\newcommand{\tgo}{\mathfrak{t}}
\newcommand{\lgo}{\mathfrak{l}}

\newcommand{\zgo}{\mathfrak{z}}

\newcommand{\bgo}{\mathfrak{b}}

\newcommand{\gb}{\bar{g}}

\newcommand{\hb}{\bar{h}}
\newcommand{\xb}{\bar{x}}

\newcommand{\Yb}{\overline{Y}}

\newcommand{\al}{\alpha}

\newcommand{\la}{\lambda}

\newcommand{\back}{\backslash}

\newcommand{\bg}{\langle}
\newcommand{\bd}{\rangle}

\newcommand{\eps}{\varepsilon}

\renewcommand{\leq}{\leqslant}
\renewcommand{\geq}{\geqslant}

\title{Le lemme fondamental pondéré I \\ Constructions géométriques}
\author{Pierre-Henri Chaudouard et Gérard Laumon}
\date{}

\begin{document}

\makeatletter
\@addtoreset{equation}{section}         
\def\theequation{\thesection.\arabic{equation}}
\makeatother

\maketitle

\begin{abstract}Ce travail est la partie géométrique de notre démonstration du lemme fondamental pondéré  qui prolonge celle du  lemme fondamental de Langlands-Shelstad due à Ngô Bao Châu. 

L'approche de Ngô repose sur l'étude de partie elliptique de la fibration de Hitchin. Cette fibration a pour espace total  le champ des fibrés de Hitchin et pour  base l'espace affine des ``polynômes caractéristiques''. Au-dessus de l'ouvert  elliptique, elle est propre et le nombre de points de ses fibres  sur un corps fini s'exprime en termes d'intégrales orbitales.

  Dans cet article, on étudie la fibration de Hitchin au-dessus d'un ouvert plus gros que l'ouvert elliptique, le lieu ``génériquement semi-simple régulier''. Les fibres ne sont en général ni de type fini ni même séparées. Par analogie avec les troncatures d'Arthur, nous introduisons le champ des fibrés de Hitchin $\xi$-stables. Nous montrons que celui-ci est un champ de Deligne-Mumford,  lisse sur le corps de base et propre au-dessus de la base des polynômes caractéristiques. Nous exprimons le  nombre de points d'une fibre $\xi$-stable  sur un corps fini en termes d'intégrales orbitales pondérées d'Arthur  qui apparaissent dans la formule des traces d'Arthur-Selberg.
\end{abstract}

\renewcommand{\abstractname}{Abstract}

\begin{abstract}
This work is the geometric part of our proof of the weighted fundamental lemma, which is an extension of Ngô Bao Châu's proof of the Langlands-Shelstad fundamental lemma. 

Ngô's approach is based on a study of the elliptic part of the Hichin fibration. The total space of this fibration is the algebraic stack of Hitchin bundles and its base space is the affine space of ``characteristic polynomials''. Over the elliptic set, the Hitchin fibration is proper and the number of points of its fibers over a finite field can be expressed in terms of orbital integrals.

In this paper, we study the Hitchin fibration over an open set bigger than the elliptic set, namely the ``generically regular semi-simple set''. The fibers are in general neither of finite type nor separeted. By analogy with Arthur's truncation, we introduce the substack of $\xi$-stable Hitchin bundles. We show that it is a Deligne-Mumford stack, smooth over the base field and proper over the base space of  ``characteristic polynomials''. Moreover, the number of points of the $\xi$-stable fibers over a finite field can be expressed as a sum of weighted orbital integrals, which appear in the Arthur-Selberg trace formula.

\end{abstract}

\tableofcontents

\section{Introduction}

\begin{paragr}Ce travail est la partie géométrique de notre démonstration du lemme fondamental pondéré dont la stratégie suit celle élaborée par Ngô Bao Châu dans sa démonstration du le lemme fondamental de Langlands-Shelstad.

Cette partie géométrique poursuit deux objectifs. D'une part, on tronque la fibration de Hitchin de sorte que les fibres de Hitchin tronquées soient propres. D'autre part,  on exprime le nombre de points sur un corps fini d'une telle fibre en termes d'intégrales orbitales pondérées globales d'Arthur qui interviennent dans la formule des traces d'Arthur-Selberg ; plus exactement, il s'agit de variantes de ces intégrales pour les algèbres de Lie.

Nos constructions s'appliquent à tout groupe algébrique semi-simple et nos résultats valent dans cette généralité.  Cependant, pour les besoins de cette introduction, nous allons nous limiter au cas du groupe $\mathrm{SL}(n)$.
\end{paragr}

\begin{paragr}Soit $C$ une courbe projective, lisse, connexe, de genre $g$ sur un
corps algébriquement clos $k$, $n$ un entier $>0$ et $D$ un
diviseur effectif sur $C$ de degré $d>2g$. On suppose que la
caractéristique de $k$ est soit nulle soit $>n$.

Un fibré de Hitchin pour le groupe $\mathrm{SL}(n)$ est un couple
$(\mathcal{E},\theta )$ où $\mathcal{E}$ est fibré vectoriel
de rang $n$ sur $C$ muni d'une trivialisation de son déterminant
et $\theta:\mathcal{E}\rightarrow \mathcal{E}(D)$ est un endomorphisme
tordu de $\mathcal{E}$ de trace nulle.

Soit $\mathbb{M}$ le champ algébrique sur $k$ des fibrés de
Hitchin. La fibration de Hitchin est le morphisme
$$
f:\mathbb{M}\rightarrow\mathbb{A}
$$
de base l'espace affine
$$
\mathbb{A}=\bigoplus_{i=2}^{n}H^{0}(C,\mathcal{O}_{C}(iD)).
$$ 
qui envoie $(\mathcal{E},\theta )$ sur le polynôme caractéristique de $\theta$ noté 
$$
P_{a}(u):=u^{n}+a_{2}u^{n-2}+\cdots +a_{n},
$$
avec $a_i \in H^{0}(C,\mathcal{O}_{C}(iD))$ pour $2\leq i\leq n$.

Pour tout $a=(a_2,\ldots,a_n)\in \mathbb{A}$, on dispose de la courbe spectrale $Y_{a}$
d'équation 
$$P_{a}(u)=0$$ 
tracée sur la surface réglée $\Sigma_{D}=\mathbb{V}(\mathcal{O}_{C}(-D))$.  C'est une courbe
projective  qui n'est en générale ni lisse  ni irréductible, ni même réduite. La projection canonique $\pi_{a}:Y_{a}\rightarrow C$
est un revêtement fini de degré $n$.

Soit $\mathbb{A}^{\mathrm{reg}}$ l'ouvert où $Y_{a}$ est réduite, c'est-à-dire où l'équation $P_{a}(u)=0$ n'a que des racines simples au point générique de $C$. La fibre en $a\in \mathbb{A}^{\mathrm{reg}}$ de la fibration de Hitchin peut s'interpréter comme le champ  des $\oc_{Y_a}$-modules  sans torsion $\mathcal{F}$ qui sont  de rang $1$ en tout point générique de $Y_a$ et qui sont munis d'une trivialisation du déterminant de $\pi_{a,\ast}\mathcal{F}$ (cf. \cite{BNR}). L'ouvert $\mathbb{A}^{\mathrm{reg}}$ contient l'ouvert elliptique $\mathbb{A}^{\mathrm{ell}}$ des $a$ tels que $Y_{a}$ est irréductible et donc intègre.  La fibre de Hitchin en tout point de l'ouvert elliptique est un champ de Deligne-Mumford propre. Par contre, la fibre de Hitchin en un point $a$ non elliptique de $\mathbb{A}^{\mathrm{reg}}$ est un champ d'Artin qui n'est ni séparé ni de type fini.
\end{paragr}

\begin{paragr}Dans ce travail nous tronquons les fibres de Hitchin non elliptiques
par une notion convenable de stabilité. Notre construction
dépend de données auxiliaires que nous allons maintenant introduire.

Soit $\infty$ un point fermé de $C$ qui n'est pas dans le support de $D$. On se restreint à l'ouvert $\mathbb{A}^{\infty\textrm{-reg}}$ de $\mathbb{A}^{\mathrm{reg}}$ où $P_{a}(u)$ n'a que des racines simples au point $\infty$. Cela n'est pas très restrictif puisque les ouverts  $\mathbb{A}^{\infty\mathrm{-reg}}$ recouvrent $\mathbb{A}^{\textrm{reg}}$ quand le point $\infty$ varie. Soit
$$\Ac\to  \mathbb{A}^{\mathrm{reg}}$$
le revêtement fini étale galoisien de groupe de Galois le groupe
symétrique $\mathfrak{S}_{n}$ qui consiste, pour chaque point $a$ de $\mathbb{A}^{\mathrm{reg}}$, à se donner un ordre total sur les racines de $P_{a}(u)$ au point $\infty$. Soit
$$
f:\mathcal{M}\rightarrow \mathcal{A}
$$
la fibration déduite de la fibration de Hitchin par le changement de base $\mathcal{A}\rightarrow \mathbb{A}$.  Un point de $\mathcal{M}$ est donc un triplet $(\mathcal{E},\theta ,t)$ formé
\begin{itemize}
\item d'un fibré de Hitchin $(\mathcal{E},\theta)$ ;
\item d'un point  $t=(t_{1},\ldots ,t_{n})$ de $k^{n}$ dont les coordonnées sont deux à deux distinctes et pour lequel  la somme $t_{1}+\cdots +t_{n}$ est nulle ;
\end{itemize}
qui vérifient la condition suivante : le $n$-uplet $(t_{1},\ldots ,t_{n})$  est la collection ordonnée des valeurs propres de $\theta_{\infty}\in \mathrm{End}(\mathcal{E}_{\infty})$.

Soit $\xi\in \mathbb{R}^{n}$ tel que $\xi_{1}+\cdots +\xi_{n}=0$. La fibre en $\infty$ d'un sous-fibré vectoriel  $\mathcal{F}$ de $\mathcal{E}$ qui est stable par $\theta$, est une somme d'espaces propres pour $\theta_\infty$ c'est-à-dire on a 
$$
\mathcal{F}_{\infty}=\bigoplus_{i\in
I}\mathrm{Ker}(\theta_{\infty}-t_{i})
$$
pour un certain ensemble $I\subset\{1,\ldots ,n\}$. On pose alors
$$
\mathrm{deg}_{\xi}(\mathcal{F})=\mathrm{deg}(\mathcal{F})+\sum_{i\in
I}\xi_{i}.
$$

\begin{definition}
On dit que $(\mathcal{E},\theta ,t_{\infty})\in \mathcal{M}$ est
$\xi$-stable si pour tout sous-fibré vectoriel non nul
$\mathcal{F}\subsetneq \mathcal{E}$ stable par $\theta$ on a
$$
\frac{\mathrm{deg}_{\xi}(\mathcal{F})}{\mathrm{rang}(\mathcal{F})}
<0
$$
\end{definition}

Pour $\xi=0$, on retrouve la notion usuelle de stabilité. Nous démontrons que la condition de $\xi$-stabilité est
ouverte et définit donc un sous-champ ouvert
$\mathcal{M}^{\xi}$ de $\mathcal{M}$.

\begin{theoreme}\label{thm:sln}
Supposons que $\sum_{i\in I}\xi_{i}\notin \mathbb{Z}$ pour tout $I\subsetneq\{1,\ldots ,n\}$ non vide. Alors  le champ $\mathcal{M}^{\xi}$
est un champ de Deligne-Mumford propre sur $\mathcal{A}$ qui est lisse sur $k$.
\end{theoreme}

Notre définition de $\xi$-stabilité est réminiscente d'une définition de stabilité avec poids due à Esteves \cite{Esteves} dans le cadre des jacobiennes compactifiées.  Pour un groupe réductif quelconque,  notre définition  s'inspire des troncatures d'Arthur mais aussi du travail de Behrend (cf. \cite{Behrend}) sur la stabilité et la réduction canonique des schémas en groupes.
Notre définition de stabilité est proche de celle des fibrés ordinaires ou de Hitchin avec structure parabolique étudiée dans \cite{Heinloth-Schmitt} et \cite{Boden-Yokogawa}.

Le théorème ci-dessus se démontre à l'aide de la méthode de Langton \cite{Langton}, reprise par Nitsure \cite{Nitsure} dans le cadre des fibrés de Hitchin. En caractéristique nulle et pour un groupe général, on peut utiliser la méthode de  Faltings (cf. \cite{Faltings}). Notre approche est en fait une adaptation d'un argument  de nature adélique dû à Heinloth (cf. \cite{Heinloth}), qui a l'avantage de fonctionner encore en caractéristique non nulle.

\end{paragr}

\begin{paragr}Le corps de base $k$ est désormais une clôture algébrique d'un corps fini $\mathbb{F}_{q}$. On suppose que toutes nos données sont définies sur $\mathbb{F}_{q}$. Soit $F$ le corps de fonctions de $C$, $\mathbb{A}$ son anneau
des adèles et $\mathcal{O}$ le sous-anneau compact maximal de
$\mathbb{A}$. On a
$$
H^{0}(C,\mathcal{O}_{C}(iD))=F\cap \varpi^{-iD}\mathcal{O}
$$
où $\varpi^{D}$ est l'adèle correspondant à $D$.

Soit $(a,t)$ un point de $\mathcal{A}$ rationnel sur
$\mathbb{F}_{q}$. Soit $P_{1},\ldots ,P_{s}$ les facteurs
irréductibles de $P_{a}(u)$ vu comme un élément de $F[u]$.
Ces facteurs s'écrivent au point $\infty$ 
$$
P_{j}(\infty )(u)=\prod_{i\in I_{j}}(u-t_{i})
$$
pour une unique partition $\{1,\ldots ,n\}=I_{1}\amalg\cdots\amalg I_{s}$. Soit $(e_i)_{1\leq i\leq n}$ la base canonique de $k^n$ et, pour tout $I\subset\{1,\ldots,n\}$, $V_I$ le sous-espace de $k^n$ engendré par $e_i$ pour $i\in I$. Soit $M\subset \mathrm{SL}(n)$ le stabilisateur des sous-espaces $V_{I_j}$ pour $1\leq j\leq s$. C'est un sous-groupe de Levi. Soit $\mathfrak{m}$ et $\mathfrak{sl}(n)$ les algèbres de Lie de $M$ et $\mathrm{SL}(n)$. Soit $X\in \mathfrak{m}(F)$ un élément semi-simple. On suppose  que $X$ est régulier au sens où son centralisateur $T_X$ dans $\mathrm{SL}(n)$ est un tore maximal. Notons qu'on a $T_X\subset M$. Soit $\varphi$ une fonction localement constante à support compact sur
$\mathfrak{sl}(n,\mathbb{A})$.

L'intégrale orbitale pondérée $J_{M}(X,\varphi)$ est définie par la formule
$$
J_{M}(X,\varphi)=\int_{T_{X}(\mathbb{A})\backslash \mathrm{SL}(n,\mathbb{A})}
\varphi (g^{-1}X g)
\mathrm{v}_{M}(g)\frac{\mathrm{d}g}{\mathrm{d}t}
$$
où
\begin{itemize}
\item[-] $\mathrm{d}g$ est la mesure de Haar sur
$\mathrm{SL}(n,\mathbb{A})$ qui donne le volume $1$ au sous-groupe compact maximal
$ \mathrm{SL}(n,\mathcal{O})$ ;
\item[-] $\mathrm{d}t$ est une mesure de Haar sur le tore
$T_{X}(\mathbb{A})$ ;
\item[-] $\mathrm{v}_{M}(g)$ est la fonction poids d'Arthur définie en terme de volume qui est  invariante à gauche par $M(\mathbb{A})$.
\end{itemize}
Le volume
$$\mathrm{vol}(T_{X}(F)\backslash
T_{X}(\mathbb{A})^{1},\mathrm{d}t)
$$
est fini, où $T_{X}(\mathbb{A})^{1}\subset
T_{X}(\mathbb{A})$ est l'intersection des noyaux de tous les
homomorphismes $T_{X}(\mathbb{A})\rightarrow \mathbb{Z}$ qui
sont obtenus en composant un caractère $F$-rationnel
$T_{X}\rightarrow \mathbb{G}_{m,F}$ avec l'application degré $\mathbb{A}^{\times}\rightarrow \mathbb{Z}$.

\begin{theoreme}\label{thm:slncomptage}
Supposons  $\sum_{i\in I}\xi_{i}\notin \mathbb{Z}$ pour tout $I\subsetneq\{1,\ldots ,n\}$ non vide. Alors le nombre de points rationnels sur
$\mathbb{F}_{q}$ de la fibre en $(a,t)$  de la
fibration de Hitchin tronquée $\mathcal{M}^{\xi}\rightarrow \mathcal{A}$ ne
dépend pas de $\xi$. À un coefficient près qui ne dépend que de la normalisation de la fonction poids d'Arthur, il est égal à la somme suivante d'intégrales orbitales
pondérées d'Arthur
$$
\sum_{X}\mathrm{vol}(T_{X}(F)\backslash
T_{X}(\mathbb{A})^{1},\mathrm{d}t)
J_{M}(X,\mathbf{1}_D)
$$
où  $\mathbf{1}_D$ est la fonction caractéristique de $\varpi^{-D}\mathfrak{sl}(n,\mathcal{O})$ et où  la somme est prise sur les classes de $M(F)$-conjugaison des éléments $X\in \mathfrak{m}(F)$ tels que le polynôme caractéristique de la restriction de $X$ à $V_{I_j}$ est égal à $P_j$ pour $1\leq j\leq s$. 
\end{theoreme}
\end{paragr}

\begin{paragr}
  Décrivons rapidement l'organisation de cet article. Quelques rappels et notations sont donnés dans la section \ref{sec:morcar}. La section \ref{sec:schema-car} introduit la base $\Ac$ de notre fibration de Hitchin et donne les propriétés essentielles. On définit dans la section \ref{sec:fib-Hitchin} le champ algébrique $\mc$ des triplets de Hitchin et le morphisme de Hitchin $f\ : \ \mc\to \Ac$. On y explique également la notion de réduction à un sous-groupe parabolique d'un triplet de Hitchin et on en donne un critère d'existence et d'unicité. Dans la section \ref{sec:cvx}, on associe à chaque triplet de Hitchin un convexe qui est construit à l'aide des degrés des réductions de ce triplet aux sous-groupes paraboliques semi-standard. À l'aide de ce convexe, on définit dans la section \ref{sec:xi-stab} la $\xi$-stabilité d'un triplet de Hitchin. On y énonce également le théorème principal de cet article qui généralise le théorème \ref{thm:sln} ci-dessus (cf. théorème \ref{thm:principal}). La section \ref{sec:description} fournit une  description adélique des triplets de Hitchin. La section \ref{sec:existence} est entièrement consacrée à la démonstration  de la partie existence du critère valuatif de propreté pour le morphisme $f$ restreint au champ des triplets de Hitchin $\xi$-semi-stables.  À la section \ref{sec:separation}, on démontre que la restriction du morphisme $f$ au champ des triplets de Hitchin $\xi$-stables est séparé.  À la section \ref{sec:DM},  on prouve que le champ des triplets de Hitchin $\xi$-stables, pour des $\xi$ en ``position générale'', est un champ de Deligne-Mumford. L'article s'achève à la section \ref{sec:comptage} par le comptage des fibres de Hitchin $\xi$-semi-stables en termes d'intégrales orbitales pondérées d'Arthur (cf. le théorème \ref{thm:comptage} qui généralise le théorème \ref{thm:slncomptage}).
\end{paragr}

\begin{paragr}[Remerciements.] --- Une partie de cet article a été écrit lors d'un séjour du premier auteur nommé à l'\emph{Institute for Advanced Study} de Princeton à l'automne 2008. Il souhaite remercier cet institut pour son hospitalité ainsi que la \emph{National Science Foundation}  pour le soutien (agreement No. DMS-0635607) qui a rendu ce séjour possible.
  
\end{paragr}

\section{Morphisme caractéristique}\label{sec:morcar}

\begin{paragr}
Soit $k$ une clôture algébrique d'un corps fini $\FF_q$. Soit $G$ un groupe algébrique réductif et connexe sur $k$ et $\ggo$ son algèbre de Lie. De manière générale, si une lettre majuscule désigne un groupe, son algèbre de Lie est notée par la lettre gothique minuscule correspondante. Soit $\Int$ l'action de $G$ sur lui-même par automorphisme intérieur et $\Ad$ la représentation adjointe. Soit $T$ un sous-tore maximal de $G$  et $\tgo$ son algèbre de Lie. Soit 
$$W=W^G=W^G_T$$
le groupe de Weyl de $T$ dans $G$ et  
$$\Phi=\Phi^G=\Phi^G_T$$ 
l'ensemble des racines de $T$ dans $G$. On note $|X|$ le cardinal d'un ensemble fini $X$.  On suppose que l'ordre $|W|$ du groupe de Weyl est inversible dans $k$.

\end{paragr}

\begin{paragr}[Morphisme caractéristique.] \label{S:morcar}
  Soit $k[\ggo]$ la $k$-algèbre des fonctions régulières sur $\ggo$. Le groupe $G$ agit sur son algèbre de Lie par l'action adjointe et par dualité sur $k[\ggo]$. Soit    $k[\ggo]^{G}$ la sous-algèbre des fonctions $G$-invariantes et $\car$ le quotient catégorique de $\ggo$ par $G$
$$\car=\Spec(k[\ggo]^{G}).$$
On a donc un morphisme canonique $G$-invariant
\begin{equation}
  \label{eq:chi}
  \chi \ :\ \ggo \to \car
\end{equation}
qu'on appelle morphisme caractéristique.
 
Soit $k[\tgo]$ la $k$-algèbre des fonctions régulières sur $\tgo$ et $k[\tgo]^W$ la sous-algèbre des fonctions invariantes sous le groupe de Weyl $W$. Le morphisme de restriction $k[\ggo] \to k[\tgo]$ induit un morphisme
$$k[\ggo]^G \to k[\tgo]^W$$
qui est en fait un isomorphisme d'après le théorème de Chevalley. 
Il s'ensuit que d'une part la restriction du morphisme $\chi$ à $\tgo$, encore notée $\chi$, induit un isomorphisme 
$$\tgo//W=\Spec(k[\tgo]^W)\simeq \car$$
et d'autre part il existe $n$ éléments de $k[\ggo]^G$ (avec $n=\rang(G)$), homogènes et algébriquement indépendants, qui engendrent la $k$-algèbre  $k[\ggo]^G$ (\cite{Bki} \S\S 5 et 6).  Les degrés de ces éléments ne dépendent que du groupe  $G$. En particulier, le schéma $\car$ est isomorphe à l'espace affine standard de dimension $n=\rang(G)$.
  \end{paragr}

\begin{paragr}[Discriminant.]  Pour toute racine $\al\in \Phi$, soit $d\al\in k[\tgo]$ sa dérivée. Le discriminant $D^G$ défini par  
$$D^G=\prod_{\al\in \Phi} d\al$$
appartient à $k[\tgo]^{W}$. Soit $\car^{\reg}$ l'ouvert régulier de $\car$ c'est-à-dire l'ouvert où $D^G$ ne s'annule pas. 
\end{paragr}

\begin{paragr}[Quelques propriétés du morphisme caractéristique.]
Par définition, un élément de $\ggo$ est \emph{régulier} si  son orbite sous $G$ est de dimension maximale c'est-à-dire égale à la différence entre la dimension et le rang de $G$. Soit $\ggo^{\reg}$ l'ouvert de $\ggo$ formé des éléments réguliers et 
$$\tgo^{\reg}=\tgo\cap \ggo^{\reg}.$$

 Le morphisme 
$$\chi \ : \ \tgo \to \car$$
est un revêtement fini et  galoisien de groupe $W$. Il induit sur les ouverts réguliers un revêtement étale $\tgo^{\reg}\to \car^{\reg}$.  
\end{paragr} 

\begin{paragr}[Section de Kostant.] ---  \label{S:Kostant} Pour tout $\al\in \Phi$, le sous-espace radiciel 
$$\ggo_\al=\{X\in \ggo \ | \ \forall t\in T \ \Ad(t)X=\al(t)X\}$$
est de dimension $1$ sur $k$.

Soit $B$ un sous-groupe de Borel de $G$ qui contient $T$ et $\Delta\subset \Phi$ l'ensemble des racines simples de $T$ dans $B$. Pour tout $\al\in \Delta$, soit $X_\al$ un élément non nul de $\ggo_\al$. Soit $\al^\vee$ la coracine de $\al$. Soit $X_{-\al}$ l'unique élément de $\ggo_{-\al}$ tel que $[X_{\al},X_{-\al}]$ est égale à la coracine $\al^\vee$ vue comme élément de $\tgo$. Soit
$$X_-=\sum_{\al\in \Delta}X_{-\al}.$$
et 
$$\ggo_{X_+}=\{Y\in \ggo \ | \ [Y,X_+]=0\}.$$

D'après un théorème de Kostant, l'espace affine $X_- +\ggo_{X_+}$ est inclus dans l'ouvert régulier $\ggo^{\reg}$ et  le morphisme caractéristique $\chi \ :\ \ggo \to \car$ induit un isomorphisme 
$$X_- +\ggo_{X_+} \to \car.$$
Le morphisme inverse 
$$\eps \ : \ \car \to  X_- +\ggo_{X_+}$$
est noté simplement $\eps$ ou plus précisément $\eps_{(B,\{X_\al\}_{\al\in \Delta})}$ si l'on veut rappeler la dépendance en l'épinglage $(B,\{X_\al\}_{\al\in \Delta})$.
\end{paragr}

\begin{paragr}[Sous-groupes paraboliques et sous-groupes  de Levi.] --- \label{S:parab}Pour tous sous-groupes $M$ et $Q$ de $G$, on note $\fc^Q(M)$ l'ensemble des  sous-groupes paraboliques $P$ de $G$ qui vérifient $M\subset P\subset Q$. Soit $\fc=\fc^G(T)$ l'ensemble des sous-groupes paraboliques \emph{semi-standard} de $G$ (au sens où ils contiennent $T$).

On appelle \emph{sous-groupe de Lévi} de $G$ un facteur de Levi d'un sous-groupe parabolique de $G$. Soit $\lc=\lc^G(T)$ l'ensemble des sous-groupes de Levi \emph{semi-standard} de $G$ (au sens où ils contiennent $T$).

  Pour tout $P\in\fc$, soit $N_P$ le radical unipotent de $P$ et $M_P\in \lc$ l'unique sous-groupe de Levi de $P$ qui contient $T$.  Pour tout sous-groupe de Lévi $M\in \lc$, soit $\pc(M)$, resp. $\lc(M)$ le sous-ensemble des $P\in  \fc$ tels que $M_P=M$,  resp. des $L\in \lc$  tels que $M\subset L$.
\end{paragr}

\begin{paragr} Soit $M\in \lc$. Les notations des paragraphes précédents affublées d'un indice ou d'un exposant $M$ valent pour le groupe réductif $M$. On prendra garde à ne pas confondre  les ouverts $\car_M^{G\textrm{-}\reg}$ et $\car_M^{M\textrm{-}\reg}$ de $\car_M$ : le premier est le lieu où $D^G$ ne s'annule pas alors que le second est défini par $D^M\not=0$. Comme ce dernier ouvert n'intervient pas dans la suite, on pose
$$\car_M^{\reg}=\car_M^{G\textrm{-}\reg}.$$

De même, on distinguera les ouverts $G$-régulier $\mgo^{G\textrm{-}\reg}$ et $M$-régulier $\mgo^{M\textrm{-}\reg}$ de $\mgo$. On ne s'intéressera qu'au premier. Aussi on pose
$$\mgo^{\reg}=\mgo^{G\textrm{-}\reg}.$$

Soit $P\in \fc$ un sous-groupe parabolique semi-standard et $M=M_P$. Soit $k[\pgo]$ l'algèbre des fonctions régulières sur $\pgo$ et $k[\pgo]^P$ la sous-algèbre des fonctions $P$-invariantes. Soit le $k$-schéma affine $\car_P=\Spec(k[\pgo]^P)$ et 
\begin{equation}
  \label{eq:chiP}
  \chi_P \ :\ \pgo \to \car_P
\end{equation}
le morphisme $P$-invariant donné par l'inclusion $k[\pgo]^P\subset k[\pgo]$.

\begin{lemme} Le morphisme de restriction $k[\pgo]\to k[\mgo]$ induit un isomorphisme d'algèbre  $k[\pgo]^P\simeq k[\mgo]^M$ et un isomorphisme de schéma $\car_P\simeq \car_M$. 
\end{lemme}

\begin{preuve}
  Soit $N=N_P$ et $\mgo'\in \mgo$ l'ouvert formé des éléments semi-simples $G$-réguliers. L'ouvert $\mgo'\oplus \ngo$ de $\pgo$  est la réunion des conjugués de $\mgo'$ sous $N$. Il s'ensuit que le morphisme $k[\pgo]^P\to k[\mgo]^M$ est injectif. Il possède une section donnée par $\phi\in  k[\mgo]^M \mapsto \phi\circ p\in k[\pgo]^P$ où $p$ est la projection $\mgo\oplus\ngo\to \mgo$.
\end{preuve}

\end{paragr}

\section{Schémas  caractéristiques}\label{sec:schema-car}

\begin{paragr}[Conventions.] --- \label{S:Conventions} Soit $S$ un $k$-schéma. Pour tout $k$-schéma $U$, soit 
$$U_S= U\times_k S$$
le produit fibré de $U$ et $S$ au-dessus de $k$. Pour tous $k$-schémas $U$ et $V$ et tout $k$-morphisme 
$\phi \ : \ U \to V$, soit $\phi_S \ : \ U_S \to V_S$ le $S$-morphisme obtenu par changement de base. Lorsque $S=\Spec(K)$ est le spectre d'un corps $K$ extension de $k$, on note simplement $U_K$ le $K$-schéma et $\phi_K$ le $K$-morphisme correspondant. Si $s\in S$, on note $k(s)$ le corps local en $s$ et on pose $U_s=U_{k(s)}$ etc.

Soit $S$ un $k$-schéma et $H$ un schéma en groupes  sur $k$. Pour tout $H$-torseur $\ec$ au-dessus de $S$ et tout $S$-schéma $U$ muni d'une action à gauche de $H$, le groupe $H$ agit à droite sur le produit fibré $\ec\times_S U$ par $(e,u).h=(eh,h^{-1}u)$ pour tous $h\in H$ et $(e,u)\in \ec\times_S U$. Soit $\ec \times_S^H U$ le produit contracté par $H$ c'est-à-dire le quotient du produit $\ec \times_S U$ par $H$. Lorsque $S=\Spec(k)$, on omet l'indice $S$.
 \end{paragr}

\begin{paragr} \label{S:lacourbe}Soit $C$ une courbe projective, lisse et connexe sur $k$ de genre $g$. Soit $D$ un diviseur  effectif de degré $>2g$ et $\infty$ un point de  $C(k)$ qui ne rencontre pas le support de $D$.
 Soit $\lc_D$ "le" $\Gmk$-torseur sur $C$ associé à $D$.
\end{paragr}

\begin{paragr}[Les fibrés $\tgo_D$ et $\car_{M,D}$.] --- \label{S:tgoD} Pour tout $k$-schéma $V$ muni d'une action de $\Gmk$, soit
$$V_D=\lc_D\times_k^{\Gmk} V$$
le produit contracté. Si $V$ est un $k$-espace vectoriel et si l'action de $\Gmk$ est linéaire, $V_D$ est un fibré vectoriel sur $C$.

Soit $M\in \lc$. L'action par homothétie du groupe multiplicatif $\Gmk$ sur $\tgo$ induit une action de $\Gmk$ sur $\car_M$ pour laquelle le morphisme caractéristique $\chi_M$ est  $\Gmk$-équivariant. En outre, les actions respectives de $\Gmk$ sur $\tgo$ et $\car_M$ respectent les ouverts $\tgo^{\reg}$ et  $\car^{\reg}_M$. Par la construction précédente, on obtient un fibré vectoriel $\tgo_D$ et un fibré $\car_{M,D}$ ainsi que des ouverts
$$\tgo^{\reg}_D\subset \tgo_D$$
et
$$\car^{\reg}_{M,D}\subset\car_{M,D}.$$
On obtient également un morphisme, par abus encore noté $\chi_M$, 
$$\chi_{M} \ : \ \tgo_D \to \car_{M,D}$$
qui est  un revêtement fini, galoisien de groupe $W^M$ et étale au-dessus de l'ouvert $\car^{\reg}_D$.

Lorsqu'on prend $M=G$, on omet l'indice $M$ dans les notations.

\end{paragr}

\begin{paragr}[Schéma caractéristique $\Ac_M$.] --- \label{S:A}Soit $M\in \lc$. Soit $\Ac_M^{G\textrm{-}\reg}$ le $k$-schéma tel que pour tout $k$-schéma $S$, l'ensemble $\Ac_M^{G\textrm{-}\reg}(S)$ des $S$-points est l'ensemble des couples 
$$(a,t)$$ 
formés d'une section $a$ du fibré $\car_{M,D,S}$ au-dessus de $C_S$ et d'un point $t \in \tgo^{G\textrm{-}\reg}_D(S)$ qui vérifient 
$$a(\infty_S)=\chi_{M}(t).$$
Dans la suite, on appelle $\Ac_M^{G\textrm{-}\reg}$ le schéma caractéristique de $M$.

\begin{remarque}  L'exposant $G\textrm{-}\reg$ rappelle que le point $t \in \tgo^{G\textrm{-}\reg}_D(S)$ est $G$-régulier au sens où il appartient à $\ggo^{\reg}(S)$. Dans la suite, on omet cet exposant et on pose, abusivement,
$$\Ac_M=\Ac_M^{G\textrm{-}\reg}.$$
Lorsque $M=G$, on omet l'indice $M$ et on pose
$$\Ac=\Ac_G$$
\end{remarque}

\begin{proposition} \label{prop:irredA_M}Le schéma  $\Ac_M$ est lisse et irréductible de dimension
$$\dim(\Ac_M)=\deg(D)(|\Phi^M|/2+\rang(M)) + \rang(M)(1-g).$$
\end{proposition}

\begin{preuve}
  Le morphisme d'oubli de la donnée $t$ fait de $\Ac_M$ un revêtement fini, étale et galoisien de groupe de Galois  $W^M$  au-dessus du schéma des sections globales de $\car_{M,D}$ qui sont $G$-régulières au point $\infty$. Mais ce dernier schéma est un ouvert du schéma des sections  globales de $\car_{M,D}$, qui est non-canoniquement isomorphe à un espace vectoriel dont on voit qu'il est de la dimension annoncée par la formule de Riemann-Roch (dès que $\deg(D)>2g-2$, cf. \cite{Ngo2} lemme 4.4.1).

Le morphisme d'oubli de la donnée $a$ est un morphisme surjectif et ouvert de $\Ac_M$ sur $\tgo^{\reg}$ dès que $ \deg(D)>2g-1$ (cf.  \cite{Ngo2} preuve du lemme 5.5.2). Ses fibres sont isomorphes à des espaces affines. Il s'ensuit que $\Ac_M$ est irréductible.
\end{preuve}

\end{paragr}

\begin{paragr}[Morphismes entre schémas caractéristiques.] --- Soit $M\subset L$ deux sous-groupes de Levi semi-standard. Soit 
$$\chi^M_L  \ :\ \car_M\to\car_L$$
le morphisme défini par l'inclusion $k[\tgo]^{W^L}\subset k[\tgo]^{W^M}$. C'est un morphisme fini qui est étale au-dessus de $\car_L^{\reg}$. Comme ce morphisme est $\Gmk$-équivariant, il induit un morphisme fini, étale au-dessus de $\car_{L,D}^{\reg}$
$$\chi^M_{L,D}  \ :\ \car_{M,D}\to\car_{L,D}$$ 
Par abus de notations, on note encore $\chi^M_L$ le morphisme
$$\chi^M_L \ : \ \Ac_M \to \Ac_L$$
qui, pour tout $k$-schéma affine $S$, associe à $(a,t)\in \Ac_M(S)$ le couple $(\chi^M_{L,D}(a),t)\in \Ac_L$.

\begin{proposition}\label{prop:immersionfermee}
  Le morphisme 
$$\chi^M_L \ : \ \Ac_M \to \Ac_L$$
est une immersion fermée.
\end{proposition}

On reporte au \S \ref{S:immersionfermee} la preuve de cette proposition. 

\end{paragr}

\begin{paragr}[Caractéristiques elliptiques.] ---  \label{S:car-elliptique}Soit $M\in \lc(T)$. On pose 
$$\Ac_{M,\el}=\Ac_M- \bigcup_{L\in \lc, L\subsetneq M} \Ac_L.$$
On l'appelle le schéma  caractéristique elliptique de $M$.

\begin{proposition} Le schéma $\Ac_{M,\el}$ est un sous-schéma ouvert non vide de $\Ac_M$ et un sous-schéma localement fermé de $\Ac$. Son adhérence dans $\Ac$ est le sous-schéma fermé $\Ac_M$.
\end{proposition}

\begin{preuve} Tout d'abord, $\Ac_{M,\el}$ est non vide pour des raisons de dimension (cf. la formule de dimension de la proposition \ref{prop:irredA_M}). Le reste de la première assertion résulte de la proposition \ref{prop:immersionfermee}. La seconde assertion résulte de l'irréductibilité de $\Ac_M$ (cf. proposition \ref{prop:irredA_M}).
\end{preuve}

\begin{proposition} \label{prop:reunionAM}
 Le schéma $\Ac$ est la réunion disjointe des sous-schémas localement fermés $\Ac_{M,\el}$
pour $M\in \lc(T)$  
$$\Ac=\bigcup_{M\in \lc(T)} \Ac_{M,\el}.$$
\end{proposition}

La démonstration de cette proposition est reportée au \S \ref{S:reunionAM}.
\end{paragr}

\begin{paragr}[Preuve de la proposition \ref{prop:immersionfermee}.] --- \label{S:immersionfermee}C'est une variante des preuves du lemme 7.3 de \cite{Ngo1} et de la proposition 6.3.5 de \cite{Ngo2}. Pour la commodité du lecteur, on rappelle ici les principaux arguments fondés en partie sur le lemme suivant (cf. lemme   7.3 de \cite{Ngo1}).
  
\begin{lemme}\label{lem:relevtNgo}
  Soit $S$ un $k$-schéma normal et intègre et $U\subset S$ un sous-schéma ouvert non vide. Soit $V'$ et $V$ deux $S$-schémas. Pour tout diagramme commutatif
 $$\xymatrix{U \ar[d]_{i} \ar[r]^{h'} & V' \ar[d]^{\pi} \\ S  \ar[r]^{h}& V}$$
où $i$ est l'inclusion canonique, $h'$ et $h$ sont des sections et $\pi$ est un $S$-morphisme \emph{fini}, la section $h'$ se prolonge d'une manière unique en une section $S \to V'$ qui relève $h$.
\end{lemme}

 Montrons tout d'abord que le morphisme $\chi^M_L$ est radiciel. Il s'agit de voir que, pour tout corps $K$ extension de $k$, le morphisme  $\chi^M_L$ induit une injection sur les ensembles correspondants de $K$-points. Pour $i=1,2$ soit $(a_i,t_i)\in \Ac_M(K)$ deux éléments qui ont même image $(a,t)$ dans $\Ac_L$.  On a donc $t=t_1=t_2$ et un diagramme commutatif 
$$\xymatrix{C_K \ar@<1ex>[r]^{a_1} \ar[r]_{a_2} \ar[rd]_{a} & \car_{M,D,K} \ar[d]^{\chi_{L}^M} \\ & \car_{L,D,K}}.$$

Rappelons que le morphisme $\chi^M_{L,D}$ est fini. D'après le lemme \ref{lem:relevtNgo} ci-dessus pour que $a_1=a_2$ il suffit que $a_1$ et $a_2$ coïncident sur un ouvert de $C_K$. Or $a_1$ et $a_2$ prennent la même valeur au point $\infty_K$ à savoir $\chi_M(t)$. Comme $t$ est $G$-régulier et que le morphisme $\chi_L^M$ est étale au-dessus de $\car_{L,D,K}^{\reg}$, les sections $a_1$ et $a_2$ coïncident sur le spectre du complété de l'anneau local de $C_K$ au point $\infty_K$ donc elles coïncident sur un ouvert de $C_K$.
\medskip

Montrons ensuite que le morphisme $\chi^M_{L,D}$ est propre. Pour cela, on applique le critère valuatif de propreté. Soit $A$ un anneau de valuation discrète de corps des fractions $K$· Soit $S=\Spec(A)$. Soit $(a_M,t_M)\in \Ac_M(K)$ et $(a,t)\in \Ac_L(S)$ tel que $\chi_L^M((a_M,t_M))=(a,t)$. On a donc $t_M=t$. D'après le lemme \ref{lem:relevtNgo} ci-dessus, le morphisme $a_1$ se prolonge de manière unique en un morphisme de $C_S$ dans $\car_{M,D}$ qui s'inscrit (en pointillés) dans le diagramme commutatif suivant :

$$\xymatrix{C_K   \ar[r]^{a_1}  \ar[d] &    \car_{M,D,S} \ar[d]^{\chi_{L}^M} \\ C_S \ar[r]^{a} \ar@{.>}[ur] &      \car_{L,D,S}    }.$$ 
Une fois le prolongement de $a_1$ acquis, il reste à vérifier que $a_1\circ\infty_S=\chi^M_{L}(t)$. Cette égalité résulte de la propreté du  morphisme $\chi_{L}^M$ puisque les deux sections $a_1\circ\infty_S$ et $\chi^M_{L,D}(t)$ s'inscrivent (en pointillés) dans le diagramme commutatif suivant :

$$\xymatrix{\Spec(K)   \ar[r]^{a_1\circ \infty_K}  \ar[d] &    \car_{M,D,S} \ar[d]^{\chi_{L}^M} \\ S \ar[r]^{a\circ\infty_S} \ar@{.>}[ur] &      \car_{L,D,S}    }.$$

Montrons ensuite que le morphisme $\chi^M_{L,D}$ est non ramifié. Puisque $\Ac_M$ et $\Ac_L$ sont lisses, il revient au même de montrer que le morphisme tangent est injectif. Il s'agit donc de montrer  que $\chi_L^M$ induit une application injective $\Ac_M(k[\eps])\to \Ac_L(k[\eps])$ (où $\eps$ est une indéterminée de carré $\eps^2=0$). Soit donc $(a_1,t)$ et $(a_2,t)$ deux éléments de   $\Ac_M(k[\eps])$ qui ont même image $(a,t)$ dans  $\Ac_L(k[\eps])$. Dans la suite, on note par un indice $\eps$ le changement de base de $k$ à $k[\eps]$. Soit $a_0$ la section de $\car_{L,D}$ défini par restriction de $a$ à $C$ et $U\subset C$ l'ouvert défini comme l'image inverse par $a_0$ du lieu régulier $\cgo_{L,D}^{\reg}$. Notons que $\infty\in U$. La section $a$ induit alors un morphisme $U_\eps\to  \cgo_{L,D,\eps}^{\reg}$. Or le morphisme $\chi_{L,\eps}^M$ est étale donc non ramifié au-dessus de $\cgo_{L,D,\eps}^{\reg}$. Il s'ensuit que les sections $a_1$ et $a_2$ coïncident sur $U_\eps$ donc sur $C_\eps$ par platitude de $k[\eps]$ sur $k$.

\medskip
En conclusion, le morphisme $\chi^M_{L,D}$ qui est radiciel, propre et  non ramifié est une immersion fermée.

\end{paragr}

\begin{paragr} Dans ce paragraphe et le suivant, on rassemble quelques résultats utiles pour la suite et pour la  preuve de la  proposition \ref{prop:reunionAM}. 

 \begin{lemme} \label{lem:cjsep}Soit $P$ un sous-groupe parabolique ou un sous-groupe de Levi de $G$. Soit $F$ une extension de $k$ et $F_s$ une clôture séparable de $F$. Pour tous   $X$ et $Y$ dans $\pgo(F_s)$ et tout $a\in \car_P^{\reg}(F_s)$ tels  que
$$\chi_P(X)=\chi_P(Y)=a$$
il existe $p\in P(F_s)$ tel que $X=\Ad(p)Y$.

Si, on suppose, de plus, que $F$ est une extension de degré de transcendance $1$ sur un corps algébriquement clos et que les éléments $X$ et $Y$ appartiennent à  $\pgo(F)$, alors on peut  même prendre $p\in P(F)$. 

\end{lemme}

\begin{preuve}
  Le $F_s$-schéma défini par 
$$\{p\in P \ | \ \Ad(p)X=Y\}$$
n'est pas vide puisqu'il possède des points sur une clôture algébrique de $F$. C'est clairement un torseur sous le centralisateur $T_X$ de $X$ dans $G\times_k F_s$. Or $X$ est semi-simple et régulier donc ce dernier est un tore sur $F_s$. On en déduit que le schéma ci-dessus est lisse sur $F_s$. Il possède donc des points dans $F_s$. 

Supposons de plus que $X$ et $Y$ appartiennent à  $\pgo(F)$. Soit $p\in P(F)$  tel que $X=\Ad(p)Y$. Alors pour tout $\tau\in\Gal(F_s/F)$, l'élément $p\tau(p)^{-1}$ appartient à $T_X(F_s)$. On obtient ainsi un $1$-cocycle de $\Gal(F_s/F)$ dans $T_X$.  Si le groupe de cohomologie 
$H^1(F_s/F,T_X)$ est trivial, on peut supposer qu'on a $p\in P(F)$, quitte à modifier $p$ par un élément de $T_X(F_s)$. Comme $T_X$ est connexe, dès que $F$ est de dimension $\leq 1$, on sait que le groupe $H^1(F_s/F,T_X)$ est trivial (cf. \cite{Serre} chap. III \S2 théorème 1' et les remarques qui suivent). Cela permet de conclure puisqu'une extension de degré de transcendance $1$ sur un corps algébriquement clos est de dimension $\leq 1$.

\end{preuve}

On utilisera la définition suivante d'un élément elliptique.

\begin{definition}  Soit $F$ une extension de  $k$ et $X\in \ggo^{\reg}(F)$ un élément $G$-régulier et semi-simple. Soit $T_X$ son centralisateur dans $G_F$ (c'est un sous-$F$-tore maximal de $G_F$) et $A_X\subset T_X$ le sous-$F$-tore déployé maximal de $T_X$. On dit que $X$ est un élément \emph{elliptique} de $\ggo(F)$ si l'on a l'égalité 
$$A_X=A_G$$ 
où $A_G$ est la composante neutre du centre de $G_F$. 
\end{definition}
\end{paragr}

\begin{paragr} \label{S:lemmesAM} Soit $K$ une extension de $k$. Soit $F$ le corps des fonctions de la courbe $C_K$ sur $K$. Soit $\oc_\infty$ le complété de l'anneau local de $C_K$ en $\infty_K$ et $F_\infty$ le corps des fractions de $\oc_\infty$. Avec ces  notations, on peut énoncer la proposition suivante.

\begin{proposition}\label{prop:caractdesAM}
Pour tout $(a,t)\in \Ac_G(K)$, soit $t_a\in \tgo^{\reg}(\oc_\infty)$ l'unique relèvement de $t$ tel que la restriction de $a$ à $\Spec(\oc_\infty)$ soit égale à $\chi(t_a)$. Soit $a_\eta \in \car^{\reg}(F)$ la restriction de $a$ au point générique de $C_K$.

Pour tout $M \in \lc(T)$ et $(a,t)\in \Ac_G(K)$, les deux assertions suivantes sont équivalentes :
  \begin{enumerate}
  \item $(a,t) \in \chi^M_G(\Ac_{M}(K))$ ;
  \item il existe $X$ un élément semi-simple $G$-régulier de $\mgo(F)$ et $m\in M(F_\infty)$ tels que 

    \begin{enumerate}
    \item $\chi_G(X)=a_\eta$ ;
    \item  $\Ad(m)X=t_a$.
    \end{enumerate}
   \end{enumerate}
\end{proposition}

\begin{preuve}Soit $(a,t)\in \Ac_G(K)$. Un point $t_a$ avec les propriétés ci-dessus existe et est unique puisque $\chi$ est étale au-dessus de  $\car^{\reg}$.

Prouvons que la première assertion implique la seconde. On suppose donc qu'on a  $(a,t)\in \chi^M_G(\Ac_M(K))$, c'est-à-dire qu'il existe $(a_M,t)\in \Ac_M(K)$ tel que 

\begin{equation}
  \label{eq:chi(aM)=a}
  \chi_G^M(a_M)=a
\end{equation}
et 
\begin{equation}
  \label{eq:am=chi(t)}
  a_M(\infty_K)=\chi_M(t).
\end{equation}
Soit $a_{\infty}$ et $a_{M,\infty}$ les restrictions de $a$ et $a_M$ à $\Spec(\oc_\infty)$. Les morphismes donnés par $a_{M,\infty}$ et $\chi_M(t_a)$ coïncident en fibre spéciale par (\ref{eq:am=chi(t)})  et s'inscrivent tous deux dans le diagramme commutatif
$$\xymatrix{\Spec(\oc_\infty) \ar@<1ex>[rr]^{\chi_M(t_a)} \ar[rr]_{a_{M,\infty}} \ar[rrd]_{a_\infty} & & \car_{M}^{\reg} \ar[d]^{\chi_{G}^M} \\ & & \car^{\reg}}.$$ 
Comme $\chi_{G}^M$ est étale au-dessus de $\car^{\reg}$, on a 
$$a_{M,\infty}=\chi_M(t_a).$$
Soit $a_{M,\eta}\in \car_{M}^{\reg}(F)$ la restriction de $a_M$ au point générique.  Il existe $X\in \mgo(F)$ tel que 
\begin{equation}
  \label{eq:chiM(gamma)=aM}
  \chi_M(X)=a_{M,\eta}
\end{equation}
(pour le voir, on peut utiliser une section de Kostant relative à $M$).  Comme $a_{M,\eta}$ est $G$-régulier, $X$ est nécessairement semi-simple et $G$-régulier. On a donc l'égalité suivante dans $\car_M^{\reg}(F_\infty)$
$$ \chi_M(t_a)=\chi_M(X).$$

D'après le lemme \ref{lem:cjsep}, il existe $L$ une extension galoisienne finie de $F_\infty$ et $m\in M(L)$ tel que 
$$\Ad(m)X=t_a.$$
On déduit de cette dernière égalité que pour tout $\sigma\in \Gal(L/F_\infty)$, l'élément $m\sigma(m)^{-1}$ normalise $t_a$ : il appartient donc à $T(L)$. On obtient ainsi une $1$-cochaîne du groupe de Galois  $\Gal(L/F_\infty)$ à valeurs dans $T(L)$ qui est un cobord puisque, $T$ étant un tore déployé, le groupe $H^1(L/F_\infty,T(L))$ est trivial. Il existe donc $x\in T(L)$ tel que  $m\sigma(m)^{-1}=x\sigma(x)^{-1}$ pour tout $ \sigma\in \Gal(L/F_\infty)$. Quitte à remplacer $m$ par $x^{-1}m$, on peut supposer que $m\in M(F_\infty)$.  D'où la condition 2. (b). La condition 2.(a) résulte de (\ref{eq:chi(aM)=a}) et (\ref{eq:chiM(gamma)=aM}).

\medskip 

  Prouvons l'implication réciproque.  On suppose donc qu'il existe $X\in \mgo(F)$ semi-simple $G$-régulier et $m\in M(F_\infty)$ qui vérifient les assertions 2.(a) et 2.(b). Soit
  \begin{equation}
    \label{eq:defa_Meta}
    a_{M,\eta}=\chi_M(X)\in \cgo_M^{\reg}(F).
  \end{equation}
  Ce $F$-point s'inscrit dans le diagramme commutatif
$$\xymatrix{\Spec(F) \ar[r]^{a_{M,\eta}} \ar[d] & \car_{M,D} \ar[d]^{\chi_{G,D}^M} \\ C_K \ar[r]^{a}& \car_{G,D}}.$$
La propreté du morphisme $\chi_{G}^M$ implique que $a_{M,\eta}$ se prolonge en un morphisme 
$$a_M \ : \ C_K \to \car_{M,D}$$ 
qui relève $a$.
On a l'égalité suivante dans $\cgo_M^{\reg}(K)$ 
\begin{equation}
  \label{eq:chi(t)=aM}
  \chi_M(t)=a_M(\infty_K)
\end{equation}
En effet, vu (\ref{eq:defa_Meta}), la restriction de $a_M$ à $\Spec(F_\infty)$ est égale à $\chi_M(X)$ donc égale à $\chi_M(t_a)$ par 2.(b). Il s'ensuit que la restriction de $a_M$ à $\Spec(\oc_\infty)$ est aussi  égale à $\chi_M(t_a)$. L'égalité (\ref{eq:chi(t)=aM}) en résulte par évaluation au point $\infty_K$. 

Le couple $(a_M,t)$ définit donc un élément de $\Ac_M(K)$. Par construction de $a_M$, on a
$$\chi_G^M(a_M,t)=(a,t)$$
d'où l'assertion 1.
\end{preuve}

\begin{corollaire}\label{cor:caractdesAM}
  Avec les notations de la proposition \ref{prop:caractdesAM}, les deux assertions suivantes sont équivalentes :
  \begin{enumerate}
  \item $(a,t) \in \chi_G^M(\Ac_{M,\el}(K))$ ;
  \item il existe $X$ un élément semi-simple, $G$-régulier et elliptique de $\mgo(F)$ et $m\in M(F_\infty)$ tels que 

\begin{enumerate}
    \item $\chi_G(X)=a_\eta$ ;
    \item  $\Ad(m)X=t_a$.
    \end{enumerate}
  \end{enumerate}
\end{corollaire}

\begin{preuve}Soit $(a,t)\in \Ac_G(K)$. 

 Prouvons que la première assertion implique la seconde. On suppose donc $(a,t)\in \chi_G^M(\Ac_{M,\el}(K))$. D'après la proposition \ref{prop:caractdesAM}, il existe $X$ un élément semi-simple $G$-régulier de $\mgo(F)$ et $m\in M(F_\infty)$ qui vérifient les assertions 2.(a) et 2.(b) de la proposition \ref{prop:caractdesAM}. On va montrer que $X$ est nécessairement elliptique dans $\mgo(F)$. Pour cela on raisonne par contradiction. Soit $A_X$ le sous-$F$-tore déployé maximal de $T_X$ le centralisateur de $X$ dans $G_F$. Quitte à conjuguer $X$ par un élément de $M(F)$ et translater $m$ par ce même élément, on peut et on va supposer que $A_X$ est inclus dans le tore ``standard'' $T_F$. Il existe alors $L\in\lc^G(T)$ tel que le centralisateur de $A_X$ dans $G_F$ soit $L_F$. Notons qu'on a $L\subsetneq M$ et cette inclusion est stricte puisque $X$ n'est pas elliptique dans $M$. Mais l'assertion 2.(b) de la proposition \ref{prop:caractdesAM} implique 
$$m T_{X,F_\infty}m^{-1}=T_{F_\infty}.$$
Or les tores $T_{X,F_\infty}$ et $T_{F_\infty}$ sont deux sous-$F_\infty$-tores maximaux et $F_\infty$-déployés de $L_{F_\infty}$. Ils sont donc conjugués par un élément de $L(F_\infty)$.  On a donc 
$$m\in \Norm_{M(F_\infty)}(T) L(F_\infty).$$

Comme $T$ est un $k$-tore déployé, on a  $\Norm_{M(F_\infty)}(T)= \Norm_{M(k)}(T) T(F_\infty)$. Il existe donc $w\in  \Norm_{M(k)}(T)$ et $l\in L(F_\infty)$ tels que 
$$m=wl.$$
Soit  $L_1=wLw^{-1}$, $l_1=wlw^{-1}$ et $X_1=\Ad(w)X$. Par construction, $L_1\in \lc(T)$ est un sous-groupe de Levi propre de $M$ ; en outre, $l_1\in L_1(F_\infty)$ et $X_1\in \lgo_1(F)$ vérifient les relations $\chi_G( X_1)=a_\eta$ et $\Ad(l_1)X_1=t_a$. Il résulte alors de la proposition \ref{prop:caractdesAM} qu'on a  $(a,t)\in \chi^{L_1}_G(\AAA_{L_1}(K))$. Cela contredit l'assertion 1 et montre que $X$ est elliptique dans $\mgo(F)$.

\medskip

Montrons que l'assertion 2 implique la 1. Soit $X$ un élément semi-simple, $G$-régulier et elliptique de $\mgo(F)$ et $m\in M(F_\infty)$ qui satisfont 2.(a) et 2.(b). D'après la proposition \ref{prop:caractdesAM}, on a $(a,t)\in \chi_G^M(\Ac_M(K))$. Soit $L\in \lc(T)$ un sous-groupe de Levi semi-standard inclus dans $M$ tel que $(a,t)\in \chi^L_G(\Ac_{L,\el}(K))$. On vient de voir qu'il existe alors   $Y$ un élément semi-simple, $G$-régulier et elliptique de $\lgo(F)$ et $l\in L(F_\infty)$ tels que $\chi_G(Y)=a_\eta$ et $\Ad(l)Y=t_a$.  D'après le lemme \ref{lem:cjsep}, il existe $g\in G(F_s)$ un point de $G$ dans une clôture séparable $F_s$ de $F$ tel que $\Ad(g)X=Y$. Le centralisateur de $Y$ dans $G_F$ est un sous-$F$-tore maximal. Pour tout $\sigma\in \Gal(F_s/F)$, on a donc $\sigma(g)g^{-1}\in T_{Y}$ et l'automorphisme intérieur $\Int(g)$ de $G_{F_s}$ induit un $F$-isomorphisme de $T_X$ sur $T_{Y}$ et donc un $F$-isomorphisme entre les sous-tores $F$-déployés maximaux $A_{X}$ et $A_{Y}$. Or ces derniers ne sont autres que les centres connexes $A_{M,F}$ et $A_{L,F}$ des sous-groupes de Levi $M_F$ et $L_F$. On a donc $\dim(A_M)=\dim(A_L)$. Comme on aussi $A_M\subset A_L$, il vient $A_M=A_L$ et $L=M$. Donc  $(a,t)\in \chi_G^M(\Ac_{M,\el}(K))$, ce qu'il fallait démontrer.
\end{preuve}
  
\end{paragr}

\begin{paragr}[Preuve de la proposition \ref{prop:reunionAM}.] --- \label{S:reunionAM} Seul le fait que la réunion soit disjointe n'est pas évident. On reprend les notations du \S \ref{S:lemmesAM}. Soit $(a,t)\in \Ac_G(K)$ et $M_1$ et $M_2$ deux sous-groupes de Levi dans $\lc(T)$ tels que $a$ appartiennent à l'intersection des images de $\Ac_{M_1,\el}(K)$ et $\Ac_{M_2,\el}(K)$. D'après le corollaire \ref{cor:caractdesAM}, pour $i=1,2$, il existe $X_i\in \mgo_i(F)$ un élément elliptique et $G$-régulier et $m_i\in M_i(F_\infty)$ tels que $\Ad(m_i)X_i=t_a$ et $\chi_G(X_i)=a_\eta$. D'après le lemme \ref{lem:cjsep},   l'égalité $\chi_G(X_1)=\chi_G(X_2)$ entraîne l'existence d'un point $g\in G(F_S)$ à valeurs dans une clôture séparable $F_s$ de $F$ tel que $\Ad(g)X_1=X_2$. En raisonnant comme dans la preuve du corollaire  \ref{cor:caractdesAM}, on voit qu'on  a
$$ g A_{M_1} g^{-1}=A_{M_2}.$$
Mais $g$ conjugue $X_1$ et $X_2$ donc $m_2gm_1^{-1}$ centralise $t_a$ qui est $G$-régulier. On a donc $m_2gm_1^{-1}\in T(F_s\otimes_F F_\infty)$ puis $A_{M_1}=A_{M_2}$ et $M_1=M_2$ ce qu'il fallait voir. 
\end{paragr}

\section{La fibration de Hitchin}\label{sec:fib-Hitchin}

On poursuit avec les notations des sections précédentes.
\medskip

\begin{paragr}[Notations.] --- \label{S:VD} On complète les notations des \S\S \ref{S:Conventions} à \ref{S:tgoD}.  Soit $S$ un $k$-schéma et $\ec$ un  $H$-torseur  sur  $S$. Soit 
  $$\Aut_G(\ec)= \ec\times_k^{G,\Int} G$$
le schéma en groupes sur $S$ des automorphismes $G$-équivariants de $\ec$.

Soit $V$ un $k$-espace vectoriel  de dimension finie. Pour toute représentation algébrique linéaire 
$$\rho \ :\ G \to GL(V),$$
soit
$$\rho(\ec)=G\times_k^{G,\rho} V$$
le fibré vectoriel sur $S$ que l'on obtient en poussant $\rho$ par la représentation $\rho$.

Lorsque $\ec$ est un $G$-torseur sur $C_S=C\times_k S$ (où $C$ est la courbe  du \S\ref{S:lacourbe}) on pose
$$\rho_D(\ec)=\ec \times_{C_S}^{G,\rho} V_{D,S}$$
où $V_{D,S}=V_D\times_k S$ est le fibré vectoriel sur $C_S$ déduit du fibré vectoriel $V_D$ sur $C$ (cf.\S\ref{S:tgoD}) par changement de base. 
\end{paragr}

\begin{paragr}[Trivialisation générique des $G$-torseurs.] --- Soit $K$ un extension de $k$ et $\ec$ un $G$-torseur sur $C_K$. Dans la suite, on fera appel plusieurs fois au lemme suivant.

  \begin{lemme}\label{lem:Sen}
    Il existe $K'$ une extension finie de $K$ telle que le $G$-torseur $\ec\times_{C_K} C_{K'}$ soit trivial au point générique de $C_{K'}$.
  \end{lemme}

  \begin{preuve} Clairement, il suffit de prouver le lemme lorsque $K$ est algébriquement clos. Dans ce cas, $F$ le corps des fonctions de $C_K$, est de  dimension $\leq 1$. Soit $F_s$ une clôture séparable de $F$. Le torseur $\ec$, qui est lisse sur $C_K$, admet des sections au-dessus de $F_s$. Comme $G$ est connexe et que $F$ est  de dimension $\leq 1$, on sait que l'ensemble  $H^1(F_s/F,G)$ est réduit à la classe triviale  (cf. \cite{Serre} chap. III \S2 théorème 1' et les remarques qui suivent).
 \end{preuve}

\end{paragr}

\begin{paragr}[Réductions des torseurs à un sous-groupe.] \label{S:reduc} --- On reprend les notations du paragraphe précédent. Soit $P$ un sous-groupe algébrique de $G$. Une réduction du $G$-torseur $\ec$ sur $S$ à $P$ est la donnée d'un $P$-torseur $\ec_P$ sur $S$ et d'un isomorphisme du $G$-torseur $$\ec_P\times^{P} G,$$
où $G$ est muni de l'action de $P$ par translation à gauche, sur $\ec$. Le morphisme $P$-équivariant évident $\ec_P\to  \ec_P\times^{P} G\simeq \ec$ donne, par passage au quotient par $P$, une section 
$$\sigma_P \ :\ S \to \ec/P.$$
Réciproquement, $\ec$ est un $P$-torseur au-dessus de $\ec/P$ qui, lorsqu'on le tire en arrière par une section  de $\ec/P$, donne une réduction de $\ec$ à $P$. De cette manière, on obtient une bijection naturelle entre les sections de $\ec/P$ et les classes d'isomorphisme de réductions de $\ec$ à $P$. 
\end{paragr}

\begin{paragr}[Morphisme caractéristique.] --- Rappelons qu'on a fixé en \S\ref{S:lacourbe} une courbe $C$ sur $k$ et un point $\infty\in C(k)$. Soit $S$ un $k$-schéma et $C_S$ et $\infty_S$ les objets obtenus par changement de base. Soit $\ec$ un $G$-torseur sur $C_S$. La construction du paragraphe \ref{S:VD} donne un fibré vectoriel $\Ad_D(\ec)$ sur $C_S$. Le morphisme   caractéristique défini en \S\ref{S:morcar} (\ref{eq:chi})
$$\chi \ :\ \ggo \to \car$$ 
est $G$-invariant et $\Gmk$-équivariant. Il induit donc un morphisme par abus encore noté $\chi$
$$\chi  \ :\ \Ad_D(\ec) \to \car_{D,S}$$
et donc un morphisme toujours noté de la même façon
$$\chi \ :\ H^0(C_S,\Ad_D(\ec)) \to H^0(C_S,\car_{D,S}).$$

Le schéma en groupes $\Aut_G(\ec)$ agit sur le fibré $\Ad_D(\ec)$ par l'action adjointe. Pour toute section $\theta\in H^0(C_S,\Ad_D(\ec))$, soit
$$\Aut_G(\ec,\theta)$$
le sous-schéma en groupes de  $\Aut_G(\ec)$ qui centralise $\theta$.
\end{paragr}

\begin{paragr}[Le champ algébrique de Hitchin.] --- \label{S:Cham-Hitchin}
Introduisons la définition suivante.

  \begin{definition}
Le champ de Hitchin $\mc_G$ est le $k$-champ algébrique dont la catégorie fibre en un $k$-schéma $S$ est le groupoïde dont les objets sont les triplets $(\ec,\theta,t)$ qui vérifient les conditions suivantes :
  \begin{enumerate}
  \item $\ec$ est un $G$-torseur sur $C\times_k S$ ;
  \item $\theta$ est une section du fibré vectoriel $\Ad_D(\ec)$ ;
  \item $t$ appartient à $\tgo^{G\textrm{-}\reg}(S)$ et satisfait 
$$\chi(\theta)(\infty_S)=\chi(t)$$
  \end{enumerate}
et l'ensemble des morphismes d'un objet  $(\ec,\theta,t)$ sur un autre $(\ec',\theta',t')$ est
\begin{itemize}
\item vide si $t\not= t'$ ;
\item l'ensemble des isomorphismes $G$-équivariants de $\ec$ sur $\ec'$ tels que l'isomorphisme induit de $\Ad_D(\ec)$ sur $\Ad_D(\ec')$ envoie $\theta$ sur $\theta'$.  

\end{itemize}
\end{definition}

\begin{remarque} Ce n'est pas la définition usuelle du champ de Hitchin qui ne considère que des couples $(\ec,\theta)$. Ici la donnée supplémentaire d'un élément $t$ qui est $G$-régulier force la section $\theta$ à être semi-simple régulière en $\infty_S$ et donc à l'être aussi génériquement. Le morphisme d'oubli de $t$ fait de notre champ de Hitchin un revêtement étale et galoisien de groupe $W$ d'un ouvert du champ de Hitchin usuel.
\end{remarque}

\end{paragr}

\begin{paragr}[La fibration de Hitchin.] --- Il s'agit là du morphisme 
$$f_G \ : \ \mc_G \to \Ac_G$$
défini pour tout $k$-schéma $S$ et tout triplet  $m=(\ec,\theta,t)$ dans $\mc(S)$ par
$$f(m)=(\chi(\theta),t).$$
Dans la suite, on omet le plus souvent l'indice $G$ et on note le morphisme ci-dessus ainsi
$$f\ :\ \mc\to \Ac.$$

\end{paragr}

\begin{paragr} Soit $K$ une extension de $k$ et $P\in \fc$ un sous-groupe parabolique semi-standard. Soit $\ec$ un $P$-torseur sur $C_K$. Soit $\Ad$ l'action adjointe de $P$ sur son algèbre de Lie $\pgo$. La construction du paragraphe \ref{S:VD} appliquée au groupe $P$ fournit un fibré vectoriel $\Ad_D(\ec)$ sur $C_K$. Le morphisme   caractéristique défini en \S\ref{S:parab} (\ref{eq:chiP})
$$\chi_P  \ :\ \pgo \to \car_P$$ 
est $P$-invariant et $\Gmk$-équivariant. Il induit donc des morphismes par abus encore notés $\chi_P$
$$\chi_{P}  \ :\ \Ad_D(\ec) \to \car_{P,D,S}$$
et 
$$\chi_P \ :\ H^0(C_S,\Ad_D(\ec)) \to H^0(C_S,\car_{P,D,S}).$$

\end{paragr}

\begin{paragr}[Réduction à un sous-groupe parabolique.] --- Soit $S$ un $k$-schéma affine et $m=(\ec,\theta,t)\in \mc(S)$. Introduisons la définition suivante.

  \begin{definition} \label{def:reduc}
    On appelle réduction de $m$ à $P$ tout triplet $m_P=(\ec_P,\theta_P,t_P)$ formé de

\begin{itemize}
    \item une réduction $\ec_P$ de $\ec$ à $P$ autrement dit un $P$-torseur $\ec_P$ sur $C\times_k S$ muni d'un isomorphisme 
\begin{equation}
        \label{eq:onpousse}
        \ec_P\times_k^P G \simeq \ec \ ;
      \end{equation}
    \item une section  $\theta_P$ de $\Ad_D(\ec_P)$ sur $C\times_k S$;
    \item un point $t_P\in \tgo^{G\textrm{-}\reg}(S)$ qui vérifie
$$\chi_P(\theta_P)(\infty_S)=\chi_P(t_P)$$
\end{itemize}
qui vérifie les deux conditions suivantes :
  \begin{enumerate}
  \item l'isomorphisme, qui se déduit de (\ref{eq:onpousse}),
$$\Ad_D(\ec_P\times_k^P G)\simeq \Ad_D(\ec)$$
composé avec le morphisme 
$$\Ad_D(\ec_P)\hookrightarrow \ec_P\times_k^{P,\Ad}\ggo \simeq  \Ad_D(\ec_P\times_k^P G)$$
envoie la section $\theta_P$ sur $\theta$ ;
\item les points $t_P$ et $t$ sont égaux.
\end{enumerate}
  \end{definition}

Dans la suite, on dit que deux réduction $(\ec_P,\theta_P,t)$ et  $(\ec_P',\theta_P',t)$ de $m$ sont isomorphes s'il existe un isomorphisme de $P$-torseurs de $\ec_P$ sur $\ec_P'$ tel que l'isomorphisme qui s'en déduit $\Ad_D( \ec_P) \to \Ad_D( \ec_P')$ envoie $\theta$ sur $\theta'$. 
\end{paragr}

\begin{paragr}[Existence de réductions sur un ouvert.] --- \label{S:reductouvert} Soit $S$ un $k$-schéma test et $(\ec,\theta,t)\in \mc(S)$. Soit $M\in \lc$ et $(a_M,t)\in \Ac_M(S)$ tels que 
$$f(m)=\chi_G^M((a_M,t)).$$
  Soit $U$ l'ouvert de $C_S$, image réciproque par $a_M$ de l'ouvert $\car_M^{\reg}$.
Soit $X$ une section de $\mgo\times_k U$ tel que $\chi_M(X)$ coïncide avec $a_M$ sur l'ouvert $U$. L'existence d'une section de Kostant assure qu'une telle section existe. Notons que $X$ est partout semi-simple $G$-régulier. Soit $T_X$ le sous-schéma en tores de $M_U$ centralisateur de $X$ dans $G_U$. Introduisons alors $\ec'$ le $U$-schéma défini ainsi : pour $U$-schéma $\Omega$, l'ensemble des sections de $\ec'$ au-dessus de $\Omega$ est l'ensemble des sections de $\ec$ au-dessus de $\Omega$ telles que l'isomorphisme $\Ad_D(\ec_{\Omega}) \to \ggo \times_k \Omega$ qui s'en déduit envoie $\theta_\Omega$ sur $X_\Omega$. Comme localement pour la topologie étale les sous-schémas en tores de $G\times_k U$ sont conjugués au tore constant $T\times_k U$ et qu'on a  l'égalité 
$$\chi_G(\theta)=\chi_G^M(a_M)=\chi_G(X)$$ 
sur $U$, on voit que $\ec'$ possède des sections localement pour la topologie étale. Il est clair alors que $\ec'$ est  un $T_X$-torseur sur $U$ qui devient, une fois poussé par le morphisme $T_X \to G_U$, isomorphe à $\ec_U$. En particulier, pour tout sous-groupe parabolique $P\in \fc(M)$ soit $\ec_P$ le $P$-torseur sur $U$ qu'on déduit de $\ec'$ en poussant par le morphisme $T_X\to P_U$. On vérifie alors que $(\ec_P,\theta_U,t_U)$ est une réduction à $P$ (en un sens évident) du triplet déduit de $(\ec,\theta,t)$ par changement de base à $U$.
  \end{paragr}

\begin{paragr} Voici le principal résultat sur l'existence et l'unicité de réductions à un sous-groupe parabolique d'un triplet de Hitchin sur le spectre d'un corps $K$ extension de $k$.

  \begin{proposition}\label{prop:Preduc}
    Soit $m\in \mc(K)$ et $M\in \lc$ tel que $f(m)\in \chi^M_G(\Ac_{M,\el}(K))$ (cf. proposition \ref{prop:reunionAM}). Soit $P\in \fc$. On a l'alternative suivante :
\begin{itemize}
\item soit $P$ ne contient pas $M$ : il n'y a alors aucune réduction de $m$ à $P$ ;
\item soit $P$ contient $M$ : dans ce cas, il existe une et une seule réduction de $m$ à $P$ (à isomorphisme près).
\end{itemize}
\end{proposition}

\begin{preuve} Soit $P\in \fc$. Montrons que la condition $M\subset P$ est nécessaire à l'existence d'une réduction à $P$.  Soit $m_P=(\ec_P,\theta_P,t)$ une réduction de $m$ à $P$. Le couple $(\chi_P(\theta_P),t)$ définit un élément de $\Ac_{M_P}(K)$ dont l'image dans $\Ac_G(K)$ est égale à $f(m)$. L'hypothèse $f(m)\in \chi^M_G(\Ac_{M,\el}(K))$ et la proposition \ref{prop:reunionAM} impliquent qu'on a $M\subset M_P$ d'où $M\subset P$. 
\medskip

On suppose désormais $M\subset P$. Montrons tout d'abord l'existence d'une réduction à $P$. La construction du paragraphe \S\ref{S:reductouvert} fournit une réduction de $m$ à $P$ mais seulement sur un ouvert $U$ de $C_K$. Il s'agit alors de voir que cette réduction se prolonge autrement dit que le torseur $P$-torseur $\ec_P$ se prolonge à $C_K$. Par la bijection décrite au \S\ref{S:reduc}, il revient au même de montrer que la section $U\to \ec/P$ déduite de $\ec_P$ se prolonge à tout $C_K$. Or une telle section se prolonge par propreté  du quotient $\ec/P$  sur $C$ d'où l'existence.

\medskip

Montrons enfin l'unicité d'une réduction de $m$ à $P$. Soit $m_P=(\ec_P,\theta_P,t)$   $m'_P=(\ec_P',\theta_P',t)$ deux telles réductions.  Soit $\sigma_P$ et $\sigma_P'$ les sections de $\ec/P$ déduites de $\ec_P$ et $\ec_P'$ par la construction décrite au \S \ref{S:reduc}. Il s'agit de prouver que ces sections sont égales et il suffit, par propreté de $\ec/P$ sur $C_K$, de le prouver au point générique $\Spec(F)$ de $C_K$ et même au point $\Spec(F_s)$ où $F_s$ est une clôture séparable du corps des fonctions $F$ de la courbe $C_K$.

Or section $\sigma_P$ donne une section $\Spec(F_s)\to \ec/P$ et celle-ci est la $P(F_s)$-orbite de l'image d'un point  $e\in \ec_P(F_s)$ par le morphisme $\ec_P\to \ec$. Un tel point $e\in \ec_P(F_s)$ donne un isomorphisme $P$-équivariant
$$\Phi_P \ : \ \ec_{P,F_s}\to P\times_k F_s$$
un isomorphisme de schémas en groupes
$$\iota_P \ : \ \Aut_P(\ec_{F_s}) \to P\times_k F_s$$
et un isomorphisme dérivé
 $$d\iota_P \ : \ \Ad_P(\ec_{F_s}) \to \pgo \times_k F_s.$$
Soit $Y=d\iota_P(\theta_{F_s})\in \pgo(F_s)$. Par définition de la réduction, $\chi_P(Y)$ est égal à $a_{M,F_s}$ la restriction de la section $a_M  \ : \ C \to \cgo_M$ à $\Spec(F_s)$. En utilisant la proposition \ref{prop:caractdesAM} pour le groupe non pas $G$ mais $M$, on obtient l'existence de $X\in \mgo(F)$ semi-simple tel que
\begin{enumerate}
\item $\chi_M(X)$ est égal à la restriction $a_{M,F}$ de $a_M$ au point générique ;
\item il existe $m\in M(F_\infty)$ tel que $\Ad(m)X=t_a$ où $t_a\in \tgo^{\reg}(\oc_\infty)$ relève $t$ et vérifie $\chi_M(t_a)=t$.
\end{enumerate}
 Notons que comme $a_{M,F}$ appartient à $\car_M^{G\textrm{-}\reg}$, les éléments $X$ et $t_a$ sont automatiquement $G$-réguliers. Comme on a l'égalité $\chi_P(Y)=\chi_P(X)$ dans $\car_M^{G\textrm{-}\reg}(F_s)$, par le lemme \ref{lem:cjsep} les éléments $X$ et $Y$ sont conjugués sous $P(F_s)$. Quitte à translater $e$ par un élément de $P(F_s)$ on peut et on va supposer que $Y=X$.  L'image de $e$ dans $\ec(F_s)$ donne un  isomorphisme $G$-équivariant
 \begin{equation}
   \label{eq:Phi}
   \Phi \ : \ \ec_{F_s}\to G\times_k F_s
 \end{equation}
 un isomorphisme de schémas en groupes
 \begin{equation}
   \label{eq:iota}\iota \ : \ \Aut_G(\ec_{F_s}) \to G\times_k F_s
 \end{equation}
 et un isomorphisme dérivé
 \begin{equation}
   \label{eq:diota}d\iota \ : \ \Ad_G(\ec_{F_s}) \to \ggo \times_k F_s.
 \end{equation}
 Ce dernier est tel que $X=d\iota(\theta)$ est un élément de $\mgo(F)$ semi-simple $G$-régulier qui vérifie les assertions 1 et 2 ci-dessus. Notons par un $'$ les mêmes objets attachés à la réduction $m_P'$. Soit $x\in G(F_s)$ tel que $e'=e\cdot x$. L'automorphisme de $\ggo\times_k F_s$ donné par $d\iota'\circ d\iota^{-1}$ est alors égal à $\Ad(x)$. Par conséquent, $\Ad(x)X=X'$. On a donc les égalités
$$t_a=\Ad(m')X'=\Ad(m'x)X=\Ad(m'xm^{-1})t_a$$
ce qui implique $m'xm^{-1}\in T(F_\infty\otimes_F F_s)$ puisque $t_a$ est $G$-régulier. Il s'ensuit qu'on a $x\in M(F_s)$. Donc les sections $\sigma_P$ et $\sigma'_P$ coïncident au point générique. Ainsi il existe un isomorphisme de $\ec_P$ sur $\ec_P'$. Avec les notations précédentes, on a vu que $X$ et $X'$ sont conjugués sous $M(F_s)$. Donc on peut choisir cet isomorphisme de sorte que $\theta_P$ s'envoie sur $\theta'_P$ (au moins au point $\Spec(F_s)$ donc partout).  Cela termine la démonstration de l'unicité et de la proposition.
\end{preuve}
\end{paragr}

\section{Convexe associé à un triplet de Hitchin. }\label{sec:cvx}
\begin{paragr}\label{S:espacesa}
  Pour tout groupe $M$ défini sur $k$ soit $X^*(M)$ le groupe des caractères rationnels de $M$ et
$$X_*(M)=\Hom_{\ZZ}(X^*(M),\ZZ).$$
Lorsque $M$ est un tore, on note $X_*(M)$ s'identifie au groupe des cocaractères de $M$ par la dualité naturelle entre les groupes de caractères et de cocaractères.

Soit $M\in \lc$ un sous-groupe de Levi semi-standard et soit $A_M$ la composante neutre du centre de $M$. Le morphisme de restriction
$$X^*(M)\to X^*(A_M)$$
est injectif et son conoyau est fini. Soit $\ago_M^*$ le $\RR$-espace vectoriel défini par
$$\ago_M^*=X^*(M)\otimes_\ZZ \RR=X^*(A_M)\otimes_\ZZ \RR.$$
Son dual noté $\ago_M$ s'identifie à 
$$\ago_M=\Hom_\ZZ(X^*(M),\RR).$$
On ne confondra pas $\ago_M$ avec l'algèbre de Lie de $A_M$ qui est aussi celle du centre $Z_M$ de $M$ et qui sera notée $\zgo_M$.
Pour tout $L\in \lc(M)$, les morphismes de restriction
$$X^*(L)\to X^*(M)\to X^*(A_L)$$
induisent des morphismes
$$\ago_L^*\to \ago_M^*\to \ago_L^*$$
dont le composé est l'identité. De la sorte, $\ago_L^*$ s'identifie à un sous-espace de $\ago_M^*$. Soit
$$(\ago_M^L)^*=\Ker(\ago_M^*\to \ago_L^*).$$
On a donc une décomposition
$$\ago_M^*=\ago_L^*\oplus (\ago_M^L)^*$$
d'où dualement une décomposition
$$\ago_M=\ago_L\oplus \ago_M^L$$
où $\ago_M^L$ est l'orthogonal de $\ago_L^*$ dans $\ago_M$. Tout $\xi\in \ago_M$ s'écrit 
\begin{equation}
  \label{eq:decompoLevi}
  \xi=\xi_L+\xi^L
\end{equation}
suivant cette décomposition.

\medskip
On généralise les définitions précédentes à un sous-groupe parabolique $P\in \fc$ de la façon suivante. Soit $A_P=A_{M_P}$ et $\ago_P=\ago_{M_P}$. Soit $Q\in \fc(P)$ et  $\ago_P^Q=\ago_{M_P}^{M_Q}$; On note par un exposant $*$ le dual de $\ago_P^Q$. On a donc une décomposition
$$\ago_P=\ago_P^Q\oplus \ago_Q$$
et tout $\xi\in \ago_P$ s'écrit 
\begin{equation}
  \label{eq:decompoParab}
  \xi=\xi^Q+\xi_Q
\end{equation}
suivant cette décomposition.
\end{paragr}

\begin{paragr}\label{S:deltachapeau}
 Soit $B\in \fc(T)$ un sous-groupe de Borel de $G$. Soit $\Delta_B$ l'ensemble des racines simples de $T$ dans $B$. Soit $\Delta_B^\vee$ l'ensemble des coracines simples. Les parties $\Delta_B$ et $\Delta_B^\vee$ forment des bases respectives des espaces $\ago_T^*$ et $\ago_T^*$. 

Soit $P\in \fc(B)$ et $\Delta_{P,B}\subset \Delta_B$ le sous-ensemble formé de racines dans $N_P$. Soit $\Delta_{P,B}^\vee\subset \Delta_B^\vee$ le sous-ensemble correspondant de coracines. L'application de restriction $X^*(T)\to X^*(A_P)$ induit une bijection de $\Delta_{P,B}$ sur un ensemble noté $\Delta_P$. De même, la projection $\ago_T\to \ago_P$ induit une bijection de $\Delta_{P,B}^\vee$ sur une partie notée $\Delta_P^\vee$. Comme tout sous-groupe de Borel inclus dans $P$ qui contient $T$ appartient à l'orbite de $B$ sous le groupe de Weyl $W^{M_P}(T)$, les ensembles $\Delta_P$ et $\Delta_P^\vee$ ne dépendent pas du choix de $B$. En outre, $\Delta_P$ et $\Delta_P^\vee$ forment des bases respectives de $\ago_P^*$ et $\ago_P$. Soit $\hat{\Delta}_P$ la base de $\ago_P^*$ duale de $\Delta_P^\vee$.

\begin{lemme} \label{lem:deltachapeau}
Soit $Q\in \fc(P)$. On a alors l'inclusion naturelle $\hat{\Delta}_Q\subset\hat{\Delta}_P$.
\end{lemme}

\begin{preuve}Soit $B\subset P$ un sous-groupe de Borel qui contient $T$. On a l'inclusion $\Delta_{Q,B}^\vee\subset \Delta_{P,B}^\vee$ et le complémentaire $\Delta_{P,B}^\vee-\Delta_{Q,B}^\vee$ est inclus dans $\ago^Q_T$. Soit $\Delta_{Q,P}^\vee$ l'image de  $\Delta_{Q,B}^\vee$ par la projection $\ago_T\to \ago_P$. Le complémentaire $\Delta_{P}^\vee-\Delta_{Q,P}^\vee$ est alors inclus dans $\ago^Q_P$. La projection $\ago_P\to \ago_Q$ induit une bijection de $\Delta_{Q,P}^\vee$ sur $\Delta_Q^\vee$.

Soit $\varpi\in \hat{\Delta}_Q\subset \ago_Q^*$. Soit $\al^\vee\in \Delta_Q^\vee$ tel que $\varpi(\al^\vee)=1$ et $\varpi$ est nul sur $\Delta_Q^\vee-\{\al^\vee\}$. Soit $\beta^\vee\in \Delta_P^\vee$. On a alors l'alternative suivante. Soit $\beta^\vee\notin \Delta_{Q,P}^\vee$ et alors $\varpi(\beta^\vee)=0$. Soit   $\beta^\vee\in \Delta_{Q,P}^\vee$ et alors $\varpi(\beta^\vee)$ vaut $1$ si $\beta^\vee$ est l'unique élément de $\Delta_{Q,P}^\vee$ et $0$ sinon. Cela prouve que $\varpi\in \hat{\Delta}_P$.
  
\end{preuve}

\end{paragr}

\begin{paragr}[Cônes dans $\ago_P$.] --- Soit $P\in \fc$. Soit $^+\ago_P$ le cône \emph{obtus}, ouvert dans $\ago_P$, défini par
$$^+\ago_P=\{H\in \ago_P \ |\ \varpi(H)>0 \ \forall \varpi\in \hat{\Delta}_P \}.$$
Soit $\overline{^+\ago_P}$ l'adhérence de $^+\ago_P$ dans $\ago_T$ c'est-à-dire
$$\overline{^+\ago_P}=\{H\in \ago_P \ |\ \varpi(H)\geq 0 \ \forall \varpi\in \hat{\Delta}_P \}.$$

\begin{lemme} \label{lem:cones} Soit $P\subsetneq G$ un sous-groupe parabolique semi-standard. On a les assertions suivantes :
  \begin{enumerate}
  \item pour tout $Q\in \fc(P)$ on a $^+\ago_P \subset \, ^+\ago_Q + \ago_P^Q$ ;
  \item le cône $^+\ago_P$ est l'intersection des cônes $^+\ago_Q + \ago_P^Q$ lorsque $Q$ parcourt les sous-groupes paraboliques propres maximaux de $G$ qui contiennent $P$.
  \end{enumerate}
Ces assertions valent aussi pour les cônes fermés.
\end{lemme}

\begin{preuve}
    Soit $Q\in \fc(P)$. Il est clair qu'on a 
$$^+\ago_Q + \ago_P^Q=\{H\in \ago_P \ |\ \varpi(H)>0 \ \forall \varpi\in \hat{\Delta}_Q \}.$$
La première assertion résulte alors de l'inclusion $\hat{\Delta}_Q\subset\hat{\Delta}_P$ du lemme  \ref{lem:deltachapeau}.
Pour montrer la seconde assertion, il suffit d'après la première assertion de prouver que l'intersection des cônes $^+\ago_Q + \ago_P^Q$ lorsque $Q$ parcourt les éléments maximaux de $\fc(P)$ est inclus dans $^+\ago_P$. Soit $H$ dans cet intersection. Il s'agit de voir que $\varpi(H)>0$ pour tout $\varpi\in \hat{\Delta}_P$. Mais c'est évident car pour tout $\varpi\in \hat{\Delta}_P$ il existe un unique sous-groupe parabolique $Q$ propre et maximal qui contient $P$ tel que $\hat{\Delta}_Q=\{\varpi\}$.
  \end{preuve}

\end{paragr}

\begin{paragr}[Degré d'une réduction.] --- Soit $K$ une extension de $k$, $m\in \mc(K)$ un triplet de Hitchin, $P\in \fc$ un sous-groupe parabolique semi-standard et $m_P=(\ec,\theta,t)$ une réduction de $m$ à $P$ (cf. définition \ref{def:reduc}. Le \emph{degré} de la réduction $m_P$ est l'unique élément $\deg(m_P)\in \ago_P$ qui satisfait pour tout caractère $\mu\in X^*(P)$ l'égalité suivante
$$ \mu(\deg(m_P))=\deg(\mu(\ec))$$
 où $\mu(\ec)$ est le fibré en droites sur $C_K$ obtenu lorsqu'on pousse le torseur $\ec$ par la représentation 
$$\mu \ : \ P \to \mathrm{GL}(1)$$
et $\deg(\mu(\ec))$ est son degré.

\begin{lemme}\label{lem:projdeg}
  Soit $P\in \fc$ et $m_P$ une réduction de $m$ à $P$. Soit $Q\in \fc(Q)$ et $m_Q$ une réduction de $m$ à $Q$. Alors le degré $\deg(m_Q)$ est égal à la projection de $\deg(m_P)$ sur $\ago_Q$ suivant la décomposition $\ago_P=\ago_Q\oplus \ago_P^Q$. 
\end{lemme}

\begin{preuve}
Rappelons que la réduction $m_P$ est  un certain triplet $(\ec_P,\theta_P,t)$ (cf. définition \ref{def:reduc}). Soit $\ec_Q$ défini par $\ec_P\times_k^P Q$ et $\theta_Q$ la section de $\Ad_D(\ec_Q)$ obtenue lorsqu'on compose $\theta_P$ avec le morphisme $\Ad_D(\ec_P)\hookrightarrow \Ad_D(\ec_Q)$. On vérifie que le triplet $(\ec_Q,\theta_Q,t)$ est une réduction de $m$ à $Q$ donc il est isomorphe à $m_Q$ (cf.   proposition \ref{prop:Preduc}). Soit $\mu\in X^*(Q)$. Il est clair que $\mu(\ec_Q)$ est isomorphe à $\mu(\ec_P)$. On a donc
$$\mu(\deg(m_P))=\deg(\mu(\ec_P))=\deg(\mu(\ec_Q))=\mu(\deg(m_Q)).$$
Il s'ensuit que $\deg(m_P)-\deg(m_Q)$ appartient à $\ago_P^Q$ d'où le lemme.

\end{preuve}

\end{paragr}

\begin{paragr}[Convexes associés à un triplet de Hitchin.] --- \label{S:lescvxes}  On poursuit avec les notation sur paragraphe précédent. Soit $M$ l'unique sous-groupe de Levi semi-standard tel que $f(m)$ appartienne à $\Ac_{M,\el}(K)$ (cf. proposition \ref{prop:reunionAM}). Pour tout $P$ dans l'ensemble $\fc(M)$ des sous-groupes paraboliques de $G$ contenant $M$ soit $m_P$ l'unique réduction de $m$ à $P$ (cf. \ref{prop:Preduc}) et soit
$$ \,^+\mathcal{C}_{P,m}=-\deg(m_{P}) - \, ^+\ago_P +\ago^P_T$$
le cylindre ouvert de $\ago_T$ de base la chambre obtuse $-\deg(m_{P}) - \, ^+\ago_P$
et
$$ \overline{\,^+\mathcal{C}_{P,m}}=-\deg(m_{P}) - \, \overline{^+\ago_P} +\ago^P_T$$
 son adhérence dans $\ago_T$.
Notons que, pour $P=G$, on a  $\,^+\mathcal{C}_{G,m}=\ago_T$.

Soit  $\mathcal{C}_m$ le polyèdre  ouvert inclus dans $\ago_T$ défini par 
$$\mathcal{C}_m = \bigcap_{P\in \fc(M)}  \,^+\mathcal{C}_{P,m} $$
et soit
$$ \overline{\mathcal{C}}_m=\bigcap_{P\in \fc(M)} \overline{\,^+\mathcal{C}_{P,m}}$$
son adhérence.
\medskip
 
Si $M=G$ alors $\fc(G)=\{G\}$ et $\mathcal{C}_m= \overline{\mathcal{C}}_m=\ago_T$.

\begin{remarques}
Nous verrons une définition équivalente des polyèdres précédents à la proposition \ref{prop:envpecvx}.
  Ces définitions sont directement suggérées par la définition d'Arthur des poids qui interviennent dans la définition des intégrales orbitales pondérées (cf. \cite{dis_series} \S\S 2,3). Le lien avec les constructions d'Arthur sera plus évident lorsque nous disposerons d'une description adélique des triplets de Hitchin (cf. section \ref{S:description-adelique}). Lorsque $M=T$, on retrouve le ``polyèdre complémentaire'' que Behrend associe au schéma en groupes $\Aut_G(\ec)$ et au tore $T$ (cf. \cite{Behrend}\S6).
\end{remarques} 
  
\begin{proposition}\label{prop:descripCm} Avec les notations ci-dessus, dans la définition de $\mathcal{C}_m$ ou $\overline{\mathcal{C}}_m$, on peut remplacer $\fc(M)$ par $\pc(M)$. Si $M\subsetneq G$ on peut remplacer $\fc(M)$ par les éléments maximaux de $\fc(M)-\{G\}$.
\end{proposition}

\begin{preuve} C'est une conséquence évidente du lemme suivant.
\end{preuve}

\begin{lemme} Soit $P\in \fc(M)$ un sous-groupe parabolique semi-standard. On a les assertions suivantes :
  \begin{enumerate}
  \item pour tout $Q\in \fc(P)$ on a $ \,^+\mathcal{C}_{P,m} \subset \,^+\mathcal{C}_{Q,m}$ ;
  \item le cône $ \,^+\mathcal{C}_{P,m}$  est l'intersection des cônes $\,^+\mathcal{C}_{Q,m}$ lorsque $Q$ parcourt les sous-groupes paraboliques propres maximaux de $G$ qui contiennent $P$.
  \end{enumerate}
Ces assertions valent aussi pour les cônes fermés.
  
\end{lemme}

\begin{preuve} Soit $Q\in \fc(P)$.  Le lemme \ref{lem:projdeg} implique qu'on a l'égalité
 $$\,^+\mathcal{C}_{Q,m}=-\deg(m_{P})- \, ^+\ago_Q +\ago_T^Q.$$
Le lemme est alors une conséquence immédiate du lemme  \ref{lem:cones}.
\end{preuve}
\end{paragr}

\begin{paragr} Soit $S$  un $k$-schéma affine et $m\in \mc(S)$. Pour tout $s\in S$, soit $m_s\in \mc(k(s))$ le triplet de Hitchin sur $k(s)$ le corps local en $s$ déduit de $m$ par changement de base. Soit $\mathcal{C}_m$, resp.  $\overline{\mathcal{C}}_m$, l'application sur $S$ 
qui à tout $s\in S$ associe le convexe ouvert $\mathcal{C}_{m_s}$, resp. le convexe fermé $\overline{\mathcal{C}}_{m_s}$.

\begin{proposition} \label{prop:scontinue} Supposons $S$ noethérien. L'application $\mathcal{C}_m$  de $\mc(k(s))$ dans l'ensemble des parties de $\ago_T$ ordonné par l'inclusion est semi-continue infé\-rieurement c'est-à-dire pour toute partie $\Xi$ de $\ago_T$ l'ensemble
$$\{s\in S \ | \Xi \subset \mathcal{C}_{m_s} \}$$ 
  est ouvert.
Il en est de même pour l'application $\overline{\mathcal{C}}_m$.
\end{proposition}

\begin{preuve}Soit $m=(\ec,\theta,t)\in \mc(S)$. On ne traite que l'application  $\mathcal{C}_m$ car la même démonstration s'applique à $\overline{\mathcal{C}}_m$.

Montrons tout d'abord que $\mathcal{C}_m$ est constructible et ne prend qu'un nombre fini de valeurs sur $S$ . Pour cela, quitte à remplacer $S$ par une composante irréductible, on suppose que $S$ est irréductible. Soit $\eta$ le point générique de $S$ et $C_\eta=C\times_k k(\eta)$. Soit $m_\eta=(\ec_{\eta},\theta_{\eta},t_\eta)$ le triplet de Hitchin sur $C_{\eta}$ déduit de $m$ par le changement de base $\Spec(k(\eta))\to S$. Soit $M\in \lc$ tel que $f(m_\eta)\in \chi_G^M(\Ac_{M,\el}(k(\eta)))$ (cf. proposition \ref{prop:reunionAM}). Soit $P\in \fc(M)$  et $ m_{P,\eta}=(\ec_{P,\eta},\theta_{\eta},t_\eta)$ ``la'' réduction de $m_{\eta}$ à $P$ (cf. proposition \ref{prop:Preduc}). On a associé au \S\ref{S:reduc} à la $P$-réduction $\ec_{P,\eta}$ une section  $C_\eta \to \ec_{\eta}/P$. Celle-ci se prolonge se prolonge à tout $C_S$ privé d'un certain fermé. L'image de ce fermé par le morphisme propre $C_S\to S$ est un fermé qui ne contient pas le point générique $\eta$. Soit $U$ le complémentaire de ce fermé de $S$. Par construction, la section   $C_\eta \to \ec_{\eta}/P$ se prolonge en une section  $C_U\to \ec_{\eta}/P$. En utilisant la bijection du \S\ref{S:reduc}, on déduit de cette section un $P$-torseur $\ec_{P,U}$ qui est un prolongement de $\ec_{P,\eta}$ et qui est une réduction de $\ec_U$ à $P$. On en déduit que le triplet $(\ec_{P,U},\theta_{U},t_{|U})$ est une réduction de $m$ à $P$ sur $U$. Il s'ensuit que l'application qui à $u\in U$ associe le degré de $m_{P,u}$ est constante. Quitte à réduire $U$, l'application précédente est constante pour tout $P\in \fc(M)$. Quitte à réduire encore $U$, on peut supposer que $f(m_u)\in \chi_G^M(\Ac_{M,\el}(k(u)))$ pour tout $u\in U$. Mais alors l'application $\mathcal{C}_m$ est constante sur $U$. Par récurrence noethérienne, elle ne prend qu'un nombre fini de valeurs et il existe une partition de $S$ en des ensembles constructibles sur lesquels $\mathcal{C}$ est constante.
\medskip

Montrons ensuite que pour tout $\Xi\subset \ago_T$ l'ensemble 
  $$\{s\in S \ | \Xi \subset \mathcal{C}_{m_s} \}$$ 
 est ouvert. D'après ce qui précède, il  est constructible. Il suffit donc de montrer qu'il est stable par générisation. Pour cette question, on peut supposer que $S$ est le spectre d'un anneau de valuation discrète. Soit $s$ le point spécial et $\eta$ le point générique de $S$. Il suffit de montrer l'inclusion
 \begin{equation}
   \label{eq:incluCm}
   \mathcal{C}_{m_s}\subset \mathcal{C}_{m_\eta}.
 \end{equation}
 Soit $(a,t)=f(m)\in \Ac_G(S)$. Soit $M\in \lc$ tel que  $(a_\eta,t_\eta)$ appartienne à $\chi_G^M(\Ac_{M,\el}(k(\eta)))$ (cf. proposition \ref{prop:reunionAM}). Comme le morphisme $\chi_G^M$ est une immersion fermée (cf. proposition \ref{prop:immersionfermee}), il existe un couple  $(a_M,t)\in \Ac_M(S)$ d'image $(a,t)$ dans $\Ac_G(S)$.   
Si $M=G$ l'inclusion (\ref{eq:incluCm}) est triviale puisqu'alors  $\mathcal{C}_{m_\eta}=\ago_T$.  Supposons donc $M\subsetneq G$. En utilisant la proposition \ref{prop:descripCm}, on voit qu'il suffit montrer que pour tout $P\in \fc(M)$ maximal, on a
 $$ \,^+\mathcal{C}_{P,m_s}\subset  \,^+\mathcal{C}_{P,m_\eta}$$
ce qui revient à montrer que 
$$\varpi(\deg(m_{P,s})-\deg(m_{P,\eta}))\geq 0$$
où $\varpi$ est l'unique élément de $\hat{\Delta}_P$. Soit 
$$\rho \ : \ G \to V$$
la représentation irréductible de dimension finie de plus haut poids $\varpi$. Soit $V_\varpi\subset V$ la droite engendrée par un vecteur de plus haut poids $\varpi$. En particulier, cette droite est stable sous l'action de $P$. 

 Soit $U$ l'ouvert de $C_S$ défini comme l'image réciproque par $a_M$ de l'ouvert $G$-régulier $\car^{\reg}_M$. On a construit au \S\ref{S:reductouvert} un triplet $(\ec_{P,U},\theta_U,t_U)$ formé d'une réduction à $P$ du triplet déduit de $m$ par changement de base à $U$. On déduit de $\ec_{P,U}$ une section $U \to \ec/P$ qui, par projectivité de $\ec/P$, se prolonge à un ouvert $U_1\supset U$ qui contient tous les points de codimension $\leq 1$ de $C_S$. Notons que $U_1$ contient $C_\eta$ ainsi que le point générique de $C_s$. On obtient alors un $P$-torseur $\ec_{P,U_1}$ qui prolonge $\ec_{P,U}$ (cf. \S\ref{S:reduc}). On vérifie que le triplet $(\ec_{P,U_1},\theta_{U_1},t_{U_1})$ est une réduction de $m$ à $P$ sur $U_1$. 

 Le fibré en droite $\ec_{P,U_1}\times_k^{P,\rho} V_\varpi$ est un sous-faisceau localement facteur direct de la restriction à $U_1$ du fibré vectoriel $\rho(\ec)=\ec\times_k^{P,\rho} V$. Il se prolonge de manière unique en un fibré en droites $\vc_\varpi$ qui est un sous-faisceau de $\rho(\ec)$. Toutefois, il se peut que hors de $U_1$ ce fibré ne soit pas localement  facteur direct. Comme $U_1$ contient $C_\eta$ on a l'égalité
\begin{equation}
  \label{eq:egVeta}
\deg( \vc_{\varpi,\eta})=\varpi(\deg(m_{P,\eta})).
\end{equation}

Sur l'ouvert $U_1\cap C_s$ de $C_s$, on obtient par restriction du triplet $(\ec_{P,U_1},\theta_{U_1},t_{U_1})$ une réduction de $m_s$ à $P$. Toujours par les mêmes arguments, cette réduction se prolonge à $C_s$. On obtient donc un triplet $(\ec_{P,s},\theta_{s},t_s)$. Le fibré en droites
$$  \vc_{\varpi,s}'=\ec_{P,s}\times_k^{P,\rho} V_\varpi$$
est un sous-fibré localement facteur direct de $\rho(\ec_s)$ et comme ci-dessus on a 
\begin{equation}
  \label{eq:egVs}
\deg( \vc_{\varpi,s}')=\varpi(\deg(m_{P,s})).
\end{equation}
Mais ce fibré $\vc_{\varpi,s}'$  coïncide avec $\vc_{\varpi,s}$ sur l'ouvert $U_1\cap C_s$. Il s'ensuit que $\vc_{\varpi,s}'$  est le saturé de $\vc_{\varpi,s}$ dans $\rho(\ec_s)$. On a donc
$$\deg(\vc_{\varpi,s}')\geq \deg(\vc_{\varpi,s})$$
avec égalité si et seulement si $\vc_{\varpi,s}'=\vc_{\varpi,s}$. Avec (\ref{eq:egVeta}) et (\ref{eq:egVs}) cela donne 
$$\varpi(\deg(m_{P,s}))=\deg(\vc_{\varpi,s}')\geq \deg(\vc_{\varpi,s}) =\varpi(\deg(m_{P,\eta}))=\varpi(\deg(m_{P,\eta}))$$
comme voulu.

\end{preuve}

\end{paragr}

\begin{paragr}[Réductions à des sous-groupes paraboliques adjacents.] --- Soit $K$ une extension de $k$. Avant d'énoncer le principal résultat de ce paragraphe, nous allons l'illustrer dans le cas de $GL(2)$ muni de son sous-tore maximal standard $T$. Soit $m=(\ec,\theta,t)$ un triplet de Hitchin dans $\mc_{GL(2)}(K)$. La donnée de $(\ec,\theta)$ est équivalente à celle d'un fibré vectoriel $\vc$ de rang $2$ sur $C_K$, d'un endomorphisme tordu 
$$\theta \ : \ \vc\to\vc(D)$$ 
encore noté $\theta$. On suppose que $f(m)$ n'est pas elliptique. Il s'ensuit qu'au point générique $\eta$ de $C_K$, l'endomorphisme $\theta$ est diagonal dans une base notée $(e_1,e_2)$ de l'espace vectoriel $\vc_\eta$.  Soit $\vc_i\subset \vc$ le sous-fibré en droites déterminé par $e_i$. Le point $t\in \tgo(K)\simeq K^2$  s'écrit comme un couple   $t=(t_1,t_2)$.  Il y a une indétermination sur la numérotation des fibrés $\vc_i$ qui est levée lorsqu'on exige qu'au point $\infty$ l'endomorphisme $\theta$ ait comme valeur propre $t_i$ sur l'espace $\vc_{i,\infty}$.

Soit $B\subset GL(2)$ le sous-groupe de Borel standard et $\overline{B}$ son opposé. Il résulte de la proposition \ref{prop:Preduc} qu'il existe des réductions $m_B$ et $m_{\bar{B}}$ de $m$ à $B$ et $\overline{B}$. On identifie $\ago_T$ à $\RR^2$ de manière évidente. On vérifie les formules suivantes
$$\deg(m_B)=(\deg(\vc_1), \deg(\vc/\vc_1))$$
et
$$\deg(m_{\bar{B}})=(\deg(\vc/\vc_2), \deg(\vc_2)).$$
Il s'ensuit qu'on a
\begin{equation}
  \label{eq:GL2a}
  \deg(m_{\bar{B}})-\deg(m_B)= \mathrm{long}(\vc/(\vc_1+\vc_2))(1,-1)
\end{equation}
où la longueur notée $ \mathrm{long}$ du module de torsion $\vc/(\vc_1+\vc_2)$ est donnée par
$$\deg(\vc)-\deg(\vc_1)-\deg(\vc_2)$$
vu que le morphisme entre fibrés de rang $1$
$$\vc_1 \to \vc/\vc_2$$
est injectif.

Au point générique de $C_K$, on a une somme directe
$$\vc_\eta=\vc_{1,\eta}\oplus \vc_{2,\eta}.$$
Soit $p_i$ la projection sur $\vc_{i,\eta}$. Celle-ci commute à $\theta$ qui agit sur le facteur $ \vc_{i,\eta}$ par la section globale $\la_i$ de $\oc(D)$. Soit $p_1(\vc)$ le fibré inversible image de $\vc$ par $p_1$. La projection $p_1$ induit un isomorphisme de  $\vc/\vc_2$ sur $p_1(\vc)$. La restriction de $p_1$ à $\vc_1$ est injective et induit un isomorphisme de $\vc_1$ sur un sous-faisceau de $p_1(\vc)$. De même, on introduit le fibré $p_2(\vc)$ et on obtient des isomorphismes
$$p_1(\vc)/\vc_1\simeq \vc /(\vc_1+\vc_2) \simeq p_2(\vc)/\vc_2$$
compatibles aux endomorphismes tordus induit par $\theta$. Or $\theta$ agit par les scalaires $\la_1$ sur $p_1(\vc)$ et $\la_2$ sur $p_2(\vc)$. Il s'ensuit que le scalaire $\la_1-\la_2$, qui est nul puisque $\theta$ est régulier, agit par $0$ sur $p_1(\vc)/\vc_1 $ d'où
$$(\la_1-\la_2) p_1(\vc) \subset \vc_1(D)$$
dont on déduit l'inégalité
$$\mathrm{long}(p_1(\vc)/\vc_1) \leq 2\deg(D)$$
soit encore
\begin{equation}
  \label{eq:GL2b}
\mathrm{long}(\vc/(\vc_1+\vc_2))\leq 2\deg(D).
\end{equation}

\medskip

Ce sont les lignes (\ref{eq:GL2a}) et (\ref{eq:GL2b}) que nous allons généraliser à un groupe quelconque.
Pour tout sous-groupe de Levi $M\in \lc$, deux sous-groupes paraboliques $P$ et $Q$ dans $\pc(M)$ sont dits \emph{adjacents} si l'intersection 
$$\Delta_P^\vee \cap (-\Delta_Q^\vee)$$
 est  réduite à un élément. 

\begin{lemme}\label{lem:pos} Soit $m\in \mc(K)$ et $M\in \lc$ tel que $f(m)\in \Ac_{M,\el}(K)$.

Soit $P$ et $Q$ deux sous-groupes paraboliques adjacents dans $\pc(M)$. Soit $\al^\vee$ l'unique coracine qui vérifie
$$\Delta_P^\vee \cap (-\Delta_Q^\vee)=\{\al^\vee\}.$$ 

Soit $m_P$ et $m_Q$ les réductions respectives de $m$ à $P$ et $Q$ (cf. proposition \ref{prop:Preduc}).

Il existe un unique rationel $x_\al$ tel que
    \begin{equation}
      \label{eq:positiv}
      -\deg(m_P)+\deg(m_Q) = x_\al \al^\vee.
    \end{equation}
En outre, $x_\al$ satisfait l'inégalité
\begin{equation}
  \label{eq:positiv2}
  0 \leq x_\al \leq  2\mu(\al^\vee)^{-1}(\dim(\ngo_P/(\ngo_P\cap \ngo_Q)))^2 \dim(\qgo/(\ngo_P\cap \ngo_Q)) \deg(D).
\end{equation}
où le caractère $\mu\in X^*(M)$ est donné par $\det(\Ad(\cdot)_{|\ngo_P/(\ngo_P\cap \ngo_Q)})$.

En particulier, $x_\al$ est borné sur $\mc(K)$.
      \end{lemme}

      \begin{remarques}
        L'inégalité de gauche dans \ref{eq:positiv2} est réminiscente du lemme 3.6 de \cite{dis_series}. Dans le contexte de la réduction de schémas en groupes réductifs à un sous-groupe de Borel, elle apparaît dans \cite{Behrend} proposition 6.6.
      \end{remarques}

\medskip

 \begin{preuve} Soit $m=(\ec,\theta,t)$ le point considéré de $\mc(K)$. Soit $m_P=(\ec_P,\theta,t)$ et $m_Q=(\ec_Q,\theta,t)$ ses réductions à $P$ et $Q$. L'énoncé est trivial si $M=G$. On suppose donc $M\subsetneq G$ dans la suite.

        Soit $R \in \fc(M)$ défini par $\hat{\Delta}_R=\hat{\Delta}_P-\{\al^\vee\}$. On a donc $\hat{\Delta}_R \cup \{-\al^\vee\}=\hat{\Delta}_Q$. Ainsi $R$ contient $P$ et $Q$. Pour tout caractère $\mu\in X^*(R)$, on a donc 
$$\deg(\mu(m_R))=\deg(\mu(m_P))=\deg(\mu(m_Q)).$$
Par conséquent, $-\deg(m_P)+\deg(m_Q)$ appartient à l'espace $\ago_M^R$ qui n'est autre que la droite engendrée par $\al^\vee$. On obtient ainsi l'existence de $x_\al$ qui vérifie (\ref{eq:positiv}).

Il reste à vérifier l'inégalité (\ref{eq:positiv2}). On commence par le cas où l'on a $R=G$ c'est-à-dire  $M$ est maximal dans $G$. Dans ce cas, $P\cap Q=M$. Pour alléger les notations, posons
$$\vc=\Ad_D(\ec),$$
$$\vc_1=\ec_P\times_k^{P,\Ad} \ngo,$$
où $\ngo$ est l'algèbre de Lie du radical unipotent de $P$, et
 $$\vc_2=\Ad_D(\ec_Q)=\ec_Q\times_k^{Q,\Ad}\qgo.$$
Alors $\vc_1$ et $\vc_2$ sont deux sous-fibrés vectoriels du fibré vectoriel $\vc$ tels que $\vc_1\cap \vc_2=0$. Au point générique de $C_K$, on a une somme directe
$$\vc_\eta=\vc_{1,\eta}\oplus \vc_{2,\eta}.$$
Soit $p_i$ la projection sur $\vc_{i,\eta}$. 
La section $\theta$ de $\vc(D)$ se factorise par $\vc_2(D)$. L'action adjointe $\ad(\theta)$ de $\theta$ induit des endomorphismes tordus $\vc_1\to \vc_1(D)$ et $\vc_2\to \vc_2(D)$. Les fibrés vectoriels $\vc_1$ et $\vc / \vc_2$ sont de même rang égal à la dimension de $\ngo$. 
Le morphisme évident
$$\vc_1 \to \vc /\vc_2$$
est injectif. Il s'ensuit qu'on a 
$$\deg( \vc_1)\leq  \deg(\vc /\vc_2)$$
avec égalité si et seulement si le morphisme est bijectif. Soit $\mu$ le caractère de $X^*(P)=X^*(M)$ donné par $\det(\Ad_{|\ngo})$. On a donc $\mu(\al^\vee)>0$.  Notons que $\mu$ est encore égal au caractère de $M$ donné par $\det(\Ad_{|\ggo/\qgo})$. Il s'ensuit qu'on a 
$$\mu(\deg(m_P))=\deg(\vc_1)\leq \deg(\vc/\vc_2)=\mu(\deg(m_Q))$$
et
$$\mu(\al^\vee) x_\al= \deg(\vc/\vc_2)-\deg(\vc_1)=\mathrm{long}(\tc)$$
où $\tc$ est le faisceau de torsion défini par
$$ \tc= \vc/ (\vc_1+\vc_2).$$
On a donc la positivité dans l'inégalité (\ref{eq:positiv2}). Il nous reste à majorer la longueur de $\tc$. Pour cela, on utilise les isomorphismes
$$p_1(\vc)/\vc_1\simeq \tc \simeq p_2(\vc)/\vc_2$$
qui sont compatibles aux endomorphismes tordus induit par $\ad(\theta)$. Soit $\chi_i$ le polynôme caractéristique de la restriction de $\ad(\theta)$ à $\vc_{i,\eta}$ et soit $\la$ le résultant de ces deux polynômes. En fait, $\la$ est même une section non nulle de $\oc(ND)$ où
$$N=\dim(\ngo)\dim(\qgo).$$
L'endomorphisme tordu $p_1(\vc)\to p_1(\vc)(ND)$ donné par la multiplication induit un endomorphisme tordu trivial de $p_1(\vc)/\vc_1$. On a donc 
$$\la p_1(\vc)\subset \vc_1(ND)$$
et la majoration
$$\mathrm{long}(\tc) \leq 2\dim(\ngo)N \deg(D)$$
ce qui montre la majoration de droite de l'inégalité (\ref{eq:positiv2}) dans le cas où $M$ est maximal.

Si $M$ n'est pas maximal, on utilise les résultats précédents au groupe réductif $M_R$ qui est l'unique élément de $\lc$ qui est un facteur de Levi de $R$.  L'algèbre de Lie de $R$ est $\pgo+\qgo$ et  celle de $M_R$ est isomorphe à  $\pgo+\qgo/(\ngo_P\cap\ngo_Q)$.

\end{preuve}
\end{paragr}

\begin{paragr}[Enveloppes convexes associés à un triplet de Hitchin.] --- On continue avec les notations du paragraphe précédent. Soit $m\in \mc(K)$ et $M\in \lc$ tel que $f(m)\in \Ac_{M,\el}(K)$. Pour tout $P\in \pc(M)$, soit $m_P$ une réduction de $m$ à $P$.
$$\mathcal{CV}_m$$
l'enveloppe convexe fermée des points
$$-\deg(m_P)$$
pour $P\in \pc(M)$.   

La positivité dans le lemme \ref{lem:pos} va nous permettre de donner une autre description du convexe $\overline{\mathcal{C}}_m$ (pour des résultats analogues dans des contextes similaires, cf. \cite{dis_series} \S\S 2,3 et \cite{Behrend} \S2).

 \begin{proposition} \label{prop:envpecvx} Pour tout $m\in \mc(K)$, on a l'égalité suivante
$$\overline{\mathcal{C}}_m=\mathcal{CV}_m+\ago_T^M+\ago_G.$$
  \end{proposition}

  \begin{preuve}    Il s'agit de prouver l'égalité entre deux ensembles qui sont clairement des parties convexes et fermées de $\ago_T$, invariantes par translation par $\ago_T^M+\ago_G$. Soit $X_P$ la projection de $-\deg(m_P)$ sur $\ago_M^G$. Il suffit donc de prouver l'égalité suivantes entre parties de $\ago_M^G$
$$\bigcap_{P\in \pc(M)} \big(X_P -\overline{ ^+\ago_P}\cap \ago_M^G \big)=\mathrm{cvx}(X_P)$$
où $\mathrm{cvx}(X_P)$ est l'enveloppe convexe des $X_P$ pour $P\in \pc(M)$.

Prouvons l'inclusion $\subset$. Soit $H$ un point de $\ago_M^G$ qui n'appartienne pas à $\mathrm{cvx}(X_P)$. Alors, il existe $\la\in (\ago_M^G)^*$ tel que
$$\la(H) > \sup_{P\in \pc(M)} \la(X_P).$$

Soit $P\in \pc(M)$ tel que $\la$ appartienne au cône positivement engendré par $\hat{\Delta}_P$. L'inégalité $\la(H) > \la(X_P))$ entraîne qu'il existe $\varpi\in \hat{\Delta}_P$ tel que 
 $$\varpi(H) > \varpi(X_P).$$
Par conséquent, $H$ n'appartient pas à $X_P -\overline{ ^+\ago_P^G}$.\\

Réciproquement, prouvons l'inclusion $\supset$. Il s'agit de voir que pour tous sous-groupes paraboliques $P$ et $Q$ dans $\pc(M)$ le vecteur $X_P-X_Q$ appartient au cône $\overline{^+\ago_P}\cap \ago_M^G$. Si $P$ et $Q$ sont adjacents, cela résulte du lemme \ref{lem:pos}. Plus généralement, en prenant une suite de longueur minimale $P_i\in \pc(M)$ pour $i=0,\ldots,n$ telle que $P_{i+1}$ et $P_i$ sont adjacents, $P_0=P$ et $P_n=Q$, on voit que l'on peut écrire  $X_P-X_Q$
comme une combinaison linéaire à coefficients positifs de racines  de $A_M$ dans $N_P$ ce qui donne le résultat voulu.

  \end{preuve}

\end{paragr}

\section{La $\xi$-stabilité}\label{sec:xi-stab}

\begin{paragr}[Champs de Hitchin $\xi$-stable.] ---  Soit $\xi\in \ago_T$.

\begin{definition} Soit $S$ un $k$-schéma affine et $m\in \mc(S)$. On dit que $m$ est $\xi$-stable, resp. $\xi$-semi-stable, si pour tout $s\in S$ on a
$$\xi\in \mathcal{C}_{m_s},$$ 
resp.
$$\xi\in \overline{\mathcal{C}}_{m_s}.$$
\end{definition}

\begin{remarques}
 Comme  $ \mathcal{C}_{m_s}$ est invariant par translation par $\ago_G$, la $\xi$-(semi-)stabilité ne dépend que de la projection sur $\ago_T^G$ de $\xi$.

Dans le cas $G=GL(n)$, cette définition de $\xi$-stabilité est réminiscente d'une définition de stabilité avec poids due à Esteves \cite{Esteves} dans le cadre des jacobiennes compactifiées. Lorsque $\xi=0$, on retrouve la notion usuelle de stabilité pour les fibrés de Hitchin. On notera également que cette définition est très proche de celle utilisée pour les fibrés (ordinaire ou de Hitchin) avec structure parabolique  étudiée par Boden et Yokogawa dans  \cite{Boden-Yokogawa} et Heinloth et Schmitt dans \cite{Heinloth-Schmitt}.
Elle s'inspire également des troncatures et du poids qu'Arthur introduit dans le contexte de la formule des traces. Notre présentation (en particulier l'utilisation du convexe $\mathcal{C}_{m}$) s'inspire aussi du travail de Behrend sur la stabilité des schémas en groupes réductifs dans \cite{Behrend}.
\end{remarques}

\begin{definition}\label{def:mcxi}
   Soit $\mc^{\xi}$, resp. $\overline{\mc^{\xi}}$ le sous-champ de $\mc$ tel que pour tout $k$-schéma affine $S$ la catégorie $\mc^{\xi}(S)$, resp. $\overline{\mc^{\xi}}$, soit la sous-catégorie pleine de $\mc(S)$ dont les objets sont les points $\xi$-stables, resp. $\xi$-semi-stables.
 \end{definition}

 \begin{definition}
   \label{def:posgenerale}
On dit que $\xi\in \ago_T$ est  \emph{en position générale} si  pour tout $P\in \fc$ tel que $P\subsetneq G$ la projection de $\xi$ sur $\ago_P$ suivant $\ago_T=\ago_T^P\oplus \ago_P$ n'appartient pas au réseau $X_*(P)$.
 \end{definition}

\begin{remarque}
   Il résulte de la définition \ref{def:mcxi} qu'on a $\mc^{\xi}\subset \overline{\mc^{\xi}}$. Si $\xi$ est en position générale,  quels que soient l'extension $K$ de $k$ et le triplet de Hitchin $m\in \mc(K)$, le  point $\xi$ ne peut pas appartenir au bord de $\overline{\mathcal{C}}_{m_s}$. Il s'ensuit que pour tout $\xi$ en position générale, on a
$$\mc^{\xi}= \overline{\mc^{\xi}}.$$
 \end{remarque}

 \begin{proposition}\label{prop:schamp-ouvert}
   Les champs $\mc^{\xi}$ et $\overline{\mc^{\xi}}$ sont des sous-champs ouverts de $\mc$.
 \end{proposition}

 \begin{preuve}
   Il suffit de vérifier que pour tout $k$-schéma affine $S$ et tout $m\in \mc(S)$ l'ensemble des $s\in S$ tel que $m_s$ appartienne à    $\mc^{\xi}$, resp. $\overline{\mc^{\xi}}$, est ouvert. Comme $\mc$ est localement de type fini, on peut se limiter à des $S$ noethériens auquel cas le résultat est donné par la semi-continuité de la fonction $\mathcal{C}_{m}$, resp. $\overline{\mathcal{C}}_{m}$ (cf. proposition \ref{prop:scontinue}). 
 \end{preuve}

 \begin{proposition}\label{prop:schptf}
Si $G$ est semi-simple, les champs   $\mc^{\xi}$ et $\overline{\mc^{\xi}}$ sont des champs de type fini sur $k$. 
 \end{proposition}

 \begin{preuve}
   Soit $B\subset G$ un sous-groupe de Borel contenant $T$. Soit $\tc_G$ le champ des $G$-torseurs sur $C$ et pour tout cocaractère $\delta\in X_*(T)$ soit $\tc_B^\delta$ le champ des $B$-torseurs de degré $\delta$ sur $C$. On sait bien que ce dernier est un champ de type fini sur $k$. On a un morphisme $\tc_B^\delta \to \tc$ qui à un $B$-torseur $\ec$ associe le $G$-torseur $\ec\times^G B$. Par ailleurs, le morphisme $\mc\to\tc$ qui, à un triplet $(\ec,\theta,t)$ associe le $G$-torseur $\ec$ est de type fini. Il suffit donc de prouver que l'image de $\mc^\xi$ ou  $\overline{\mc^{\xi}}$ par ce morphisme est de type fini. Or cette image  est recouverte, d'après le lemme suivant,  par les images de  $\tc_B^\delta$ pour $\delta$ parcourant un ensemble fini. Elle est donc bien de type fini.
 \end{preuve}

 \begin{lemme} Soit $B\subset G$ un sous-groupe de Borel contenant $T$. Il existe un ensemble fini $\mathcal{X}\subset X_*(T)$ tel que pour tout triplet $(\ec,\theta,t)\in \overline{\mc^{\xi}}(k)$ il existe une réduction $\ec_B$ de $\ec$ à $B$ telle que 
$$\deg(\ec_B)\in \mathcal{X}.$$
 \end{lemme}

 \begin{preuve}
D'après Harder (\cite{Harder1} Satz 2.1.1 et \cite{Harder2} p.253), il existe une constante $c>0$ telle que tout $G$-torseur $\ec$ sur $C$ admette une réduction $\ec_B$ à $B$ dont le  degré $\deg(\ec_B)$ appartient au cône aigu
$$\cgo=\{H\in \ago_T \ | \ \al(H)\geq -c \ \forall \al\in \Delta_B \}$$
où $\Delta_B$ est l'ensemble des racines simples de $T$ dans $B$.  Soit $N$ un entier qui majore, pour toute racine de $T$ dans $B$, la somme de ses coefficients dans la base $\Delta_B$. Soit $d$ un entier qui vérifie
\begin{equation}
  \label{eq:inegsurd}
  d> \deg(D) +2N c.
\end{equation}

Soit $P$ un sous-groupe parabolique de $G$ qui contient $B$ et $\Delta_P\subset \Delta_B$ le sous-ensemble des racines simples dans $N_P$. Soit $\cgo_P$ le cône dans $\ago_T$ défini par  
$$\cgo_P=\{H\in \ago_T \ | \ \al(H)\geq d \ \forall \al\in \Delta_P \text{ et } \al(H)\leq d \ \forall \al\in \Delta_B- \Delta_P \}.$$
Lorsque $P$ parcourt l'ensembles des sous-groupes paraboliques standard, les cônes $\cgo_P$ recouvrent $\ago_T$. Il suffit donc de prouver le résultat suivant : il existe un ensemble fini $\mathcal{X}_P\subset X_*(T)$ tel que pour tout  triplet $(\ec,\theta,t)\in \overline{\mc^{\xi}}(k)$ et toute réduction $\ec_B$ de $\ec$ à $B$ dont le degré vérifie 
$$\deg(\ec_B)\in \cgo\cap \cgo_P$$
on a $\deg(\ec_B)\in \mathcal{X}_P$. 

\begin{lemme}
  Soit $(\ec,\theta,t)\in \overline{\mc^{\xi}}(k)$ et $\ec_B$ une réduction de $\ec$ à $B$ dont le degré vérifie 
$$\deg(\ec_B)\in \cgo\cap \cgo_P.$$
Alors $\theta$ se factorise par $\ec_B\times^{B,\Ad} \pgo_D\hookrightarrow \Ad_D(\ec)$. 
\end{lemme}

\begin{preuve}
  L'action adjointe de $\theta$ induit un morphisme de schémas en algèbres de Lie
$$\ec\times^{G,\Ad} \ggo\to  \Ad_D(\ec)=\ec \times^{G,\Ad}\ggo_D.$$
Il suffit de prouver que ce morphisme envoie $\ec_B\times^B \bgo$ dans $\ec_B\times^B \pgo_D$. Considérons alors une filtration de $\ggo$
$$(0)=\bgo_0 \subset \bgo_1 \subset \ldots \subset \bgo_r=\bgo \subset \pgo=\pgo_s \subset \pgo_{s-1} \subset \ldots \subset \pgo_0=\ggo$$
par des sous-espaces $B$-stables de sorte que les quotients $\bgo_{i+1}/ \bgo_{i}$ et  $\pgo_{i}/ \pgo_{i+1}$ soient de dimension $1$. 

Montrons que pour tout $0\leq i\leq r$ et $1\leq j\leq s$ l'action adjointe de  $\theta$ envoie $\ec_B\times^{B,\Ad}\bgo_{i}$ dans $\ec_B\times^{B,\Ad}\pgo_{j,D}$. Le cas $i=r$ et $j=s$ donne le résultat cherché. On raisonne pour cela par récurrence. Le cas $i=0$, trivial, amorce la récurrence. Pour un couple $(i,j)$, l'hypothèse de récurrence est que le résultat vaut pour tout couple $(i',j')$ avec soit $i'<i$ et $j'\leq s$ soit $i'=i$ et $j'\leq j$. Supposons tout d'abord $j<s$. Prouvons l'assertion pour le couple $(i,j+1)$.  
Par hypothèse de récurrence, on sait que  $\theta$ envoie $\ec_B\times^{B,\Ad}\bgo_{i-1}$ dans $\ec_B\times^{B,\Ad}\pgo_{D}$ et $\ec_B\times^{B,\Ad}\bgo_{i}$ dans $\ec_B\times^{B,\Ad}\pgo_{j,D}$. En particulier, $\theta$ induit un morphisme de fibrés en droites
\begin{equation}
  \label{eq:morphfibdroites}
   \ec_B\times^{B,\Ad}\bgo_i/ \bgo_{i-1} \to \ec_B\times^{B,\Ad}\pgo_{j,D}/ \pgo_{j+1,D}.
 \end{equation}
 Si ce morphisme n'est pas nul, le degré du but doit être supérieur au degré de la source. Il existe des caractères  disons $\al\in \phi^B_T\cup \{0\}$ et $\beta\in \Phi^{N_P}_T$ de sorte que $T$ agisse sur les quotients  $\bgo_i/ \bgo_{i-1}$ et $\pgo_{j}/ \pgo_{j+1}$ respectivement par $\al$ et $-\beta$. Il s'ensuit que le degré du but est $-\beta(\deg(\ec_B))+\deg(D)$ et celui de la source est $\al(\deg(\ec_B))$. On doit donc avoir l'inégalité
$$  (\al+\beta)(\deg(\ec_B)) \leq \deg(D)$$
Or $\beta$ est une racine de $T$ dans $N_P$ donc elle a une composante non nulle au moins sur un élément de $\Delta_P$. Comme $\deg(\ec_B)\in \cgo\cap\cgo_P$, on a la minoration 
$$\beta(\deg(\ec_B))\geq d -Nc.$$
En minorant trivialement $\al(\deg(\ec_B))$ par $-Nc$, on voit qu'on contredit l'inégalité (\ref{eq:inegsurd}) ci-dessus. En conclusion, le morphisme (\ref{eq:morphfibdroites}) est nul d'où l'assertion pour le couple $(i,j+1)$. Lorsque $j=s$, on doit prouver que l'assertion vaut pour le couple $(i+1,1)$ ce qui se démontre de manière analogue.

\end{preuve}

Soit le $P$-torseur $\ec_P=\ec_B\times^{B} P$. D'après le lemme précédent, $\theta$ se factorise par $\Ad_D(\ec_P)$. \emph{A priori}, on sait que $\chi_G(t)$ est égal à la valeur de $\chi_{G}(\theta)$ au point $\infty$ mais il n'y a pas de raison pour que $\chi_P(t)$ soit égal à la valeur de $\chi_{P}(\theta)$ au point $\infty$. Cependant, il existe un élément $w$ du groupe de Weyl $W$ tel que $\chi_P(w\cdot t)$  soit égal à la valeur de $\chi_{P}(\theta)$ au point $\infty$. Il s'ensuit que le triplet $(\ec_P,\theta,w\cdot t)$ est une réduction de $(\ec,\theta,w\cdot t)$ à $P$. Il est clair sur les définitions que  $(\ec,\theta,w\cdot t)$ appartient à  $\overline{\mc^{w\cdot\xi}}$. Quitte à changer $\xi$ en $w\cdot\xi$, on suppose dans la suite qu'on a $w=1$.

Si $P=G$, le convexe $\cgo\cap \cgo_G$ est compact et $\deg(\ec_B)$ appartient à l'ensemble fini $X_*(T)\cap \cgo\cap\cgo_G$. Si $P\not=G$, on a par $\xi$-semi-stabilité de $(\ec,\theta,t)$ la condition
$$\xi \in -\deg(\ec_P)- \, ^+\ago_P +\ago^P_T.$$
Ainsi $\deg(\ec_B)$ appartient au convexe compact
$$\cgo\cap \cgo_P \cap (-\xi - \, ^+\ago_P +\ago^P_T)$$
donc de nouveau à un ensemble fini.

 \end{preuve}

\end{paragr}

 \begin{paragr}[Le théorème principal.] --- Introduisons la définition suivante.   

 \begin{definition}
   Soit $f^\xi$, resp. $\overline{f^\xi}$, la restriction du morphisme $f$ de Hitchin à  $\mc^{\xi}$, resp.  $\overline{\mc^{\xi}}$.
 \end{definition}

Nous pouvons énoncer alors notre principal résultat. 

 \begin{theoreme}
   \label{thm:principal}
Soit $\xi\in \ago_T$ en position générale. On suppose que $G$ est semi-simple. Le champ $\mc^\xi$ est un champ de Deligne-Mumford, lisse et de type fini sur $k$. De plus, le morphisme de Hitchin $f^\xi$ est propre.
 \end{theoreme}

 \begin{preuve}
   D'après la proposition \ref{prop:schamp-ouvert}, $\mc^\xi$ est un sous-champ ouvert de $\mc$. Or ce dernier est lisse sur $k$ (Biswas et Ramanan l'ont démontré par un calcul de déformation dans \cite{Biswas-Ramanan} qui a été repris par Ngô cf. \cite{Ngo1} proposition 5.3 et \cite{Ngo2}) d'où la lissité de  $\mc^\xi$. On a vu dans la proposition \ref{prop:schptf} que le champ $\mc^\xi$ est de type fini.
On montrera ultérieurement que $\mc^\xi$ est un champ de Deligne-Mumford (cf. la démonstration du théorème \ref{thm:DM}). La propreté de $f^\xi$ est satisfaite par le critère valuatif qui combine les résultats des théorèmes \ref{thm:existence} et \ref{thm:separation} qui seront énoncés dans la suite.
 \end{preuve}

\end{paragr}

\begin{paragr}[Les $\xi$-points de Harder-Narasimhan.] ---  \label{S:HN} On munit l'espace $\ago_T$ d'un produit scalaire $\bg\cdot,\cdot\bd$ invariant par le groupe de Weyl $W^G(T)$ qui fait de $\ago_T$ un espace euclidien. Alors $\ago_T$ est canoniquement isomorphe à son dual $\ago_T^*$. Dans cet isomorphisme, une racine et sa coracine sont égales à un coefficient strictement positif près.  Soit $\|\cdot\|$ la norme euclidienne sur $\ago_T$ associée à ce produit scalaire. Notons que les espaces $\ago_M$ et $\ago_T^M$ sont alors orthogonaux (pour $M\in \lc$). 

Soit $K$ une extension de $k$ et  $\xi\in \ago_T$.

\begin{definition} Soit $m\in \mc(K)$. Le $\xi$-point de Harder-Narasimhan de $m$ est l'unique point $\varrho_m\in \overline{\mathcal{C}}_m$ qui vérifie
 $$\|\varrho_m-\xi\|=\min_{X\in \overline{\mathcal{C}}_m} \|X-\xi\|.$$
  
\end{definition}

\begin{remarque} Comme $\overline{\mathcal{C}}_m$ est un convexe fermé non vide (ce qui résulte de la proposition \ref{prop:envpecvx}), il existe un et un seul $\xi$-point de Harder-Narasimhan.
\end{remarque}
\medskip 
Pour tout sous-groupe parabolique $P\in \fc$, soit $\ago_P^{+}$ le cône \emph{aigu}, ouvert dans $\ago_P$, défini par
\begin{equation}
  \label{eq:aigu}
  \ago_P^{+}=\{H\in \ago_P \ |\ \al(H)>0 \ \forall \al\in \Delta_P \}.
\end{equation}
On obtient ainsi une partition de $\ago_T$
\begin{equation}
  \label{eq:partition}
  \ago_T=\bigcup_{P\in \fc}\ago_P^{+}.
\end{equation}

La proposition suivante et sa démonstration sont essentiellement une reformulation de \cite{Behrend} proposition 3.13.

\begin{proposition} \label{prop:HN} Soit $m\in \mc(K)$ et  $M\in \lc$ tel que $f(m)\in \Ac_{M,\el}$.  Pour tout $P\in \fc(M)$ soit $m_P$ une réduction de $m$ à $P$. Soit $\varrho\in \overline{\mathcal{C}}_m$.
Les deux assertions suivantes sont équivalentes
  \begin{enumerate}
  \item le point $\varrho$ est le $\xi$-point de Harder-Narasimhan de $m$ ;
 \item il existe un sous-groupe parabolique $Q\in \fc(M)$ tel qu'on ait 
    \begin{enumerate}
    \item $\xi-\varrho\in \ago_Q^+\cap \ago_T^G$ ;
    \item la projection de $\varrho$ sur $\ago_M^G$ appartient à la projection sur $\ago_M^G$ de l'enveloppe convexe des points
$$-\deg(m_P)$$
pour tout $P\in \pc(M)$ inclus dans $Q$ ;
\end{enumerate}
\item il existe un sous-groupe parabolique $Q\in \fc(M)$ tel qu'on ait 
    \begin{enumerate}
    \item $\xi-\varrho\in \ago_Q^+\cap \ago_T^G$ ;
    \item $\varrho$ appartient au sous-espace affine 
$$-\deg(m_Q)+\ago_T^Q+\ago_G.$$
    \end{enumerate}
  \end{enumerate}
  
\end{proposition}
  \begin{preuve} Montrons que la première assertion implique la deuxième. Supposons que $\varrho$ soit le $\xi$-point de Harder-Narasimhan de $m$. Soit $Q$ l'unique élément de $\fc$ tel que  $\xi-\varrho\in \ago_Q^+$, cf. la partition (\ref{eq:partition}). Comme $\overline{\mathcal{C}}_m$ est stable par translation par $\ago_T^M+\ago_G$, on a 
\begin{equation}
      \label{eq:1ereinclusion}
      \xi-\varrho \in \ago_M^G.
    \end{equation}
Par conséquent $Q$ doit contenir $M$. On a donc vérifié 2.(a). Vérifions 2.(b). Soit $P_0\in \pc(M)$ tel que $P_0\subset Q$. D'après la proposition \ref{prop:envpecvx}, le point $-\deg(m_{P_0})$ appartient à $\overline{\mathcal{C}}_m$. Par convexité de $\overline{\mathcal{C}}_m$, pour tout $0\leq \la \leq 1$, le point $(1-\la)\varrho -\la\deg(m_{P_0})$ appartient aussi à $\overline{\mathcal{C}}_m$ et le réel
$$\| (1-\la)\varrho -\la\deg(m_{P_0})-\xi\|^2$$
 atteint son minimum en $\la=0$. En dérivant en $\la=0$, on tombe sur la condition
 \begin{equation}
   \label{eq:derivation}
   \bg \varrho+\deg(m_{P_0}), \xi-\varrho\bd \geq 0.
 \end{equation}
Comme $\xi-\varrho\in \ago_Q^+\cap \ago_T^G$, on peut écrire 
\begin{equation}
  \label{eq:Yvarpi}
\xi-\varrho =\sum_{\varpi \in \hat{\Delta}_Q} y_\varpi \,\varpi
\end{equation}
avec $y_\varpi>0$. La condition (\ref{eq:derivation}) ci-dessus implique
\begin{equation}
  \label{eq:ineg}\sum_{\varpi \in \hat{\Delta}_Q} y_\varpi \,\varpi( \varrho+\deg(m_{P_0}))\geq 0.
\end{equation}
Pour tout  $\varpi\in \hat{\Delta}_Q$, on a 
$$\varpi( \varrho+\deg(m_{P_0}))= \varpi( \varrho+\deg(m_Q))$$
et cette quantité est négative puisque $\varrho\in \overline{\mathcal{C}}_m$. L'inégalité (\ref{eq:ineg}) n'est donc possible que si pour tout $\varpi\in \hat{\Delta}_Q$ on a 
 
\begin{equation}
  \label{eq:0}
\varpi( \varrho+\deg(m_{P_0}))=0.
\end{equation}
Pour tout $P\in \pc(M)$, soit $\la_P\in \RR$ qui satifait les conditions suivantes
\begin{itemize}
\item $0\leq \la_P \leq 1$ ;
\item $\sum_{P\in \pc(M)}\la_P=1$ ;
\item le projeté de $\varrho$ sur $\ago_M^G$ est égal au projeté sur $\ago_M^G$ du point
$$\sum_{P\in \pc(M)} -\la_P \deg(m_P).$$
\end{itemize}
De tels réels $\la_P$ existent par la proposition \ref{prop:envpecvx}. Il s'agit de  voir qu'on peut choisir $\la_P$ de sorte que $\la_P=0$ si $P\not\subset Q$. Par (\ref{eq:0}), on obtient
\begin{equation}
  \label{eq:=}
  \varpi(\deg(m_{P_0}))= \sum_{P\in \pc(M)} \la_P \,  \varpi(\deg(m_P))
\end{equation}
pour tout $\varpi\in \hat{\Delta}_Q$.
Le lemme \ref{lem:pos} implique que pour $P\in \pc(M)$ la différence $\deg(m_P)-\deg(m_Q)$ est une combinaison linéaire à coefficients positifs d'éléments de $\Delta^\vee_P$. Par ailleurs, $\varpi$ ne prend que des valeurs positives sur  $\Delta^\vee_P$ et elle est nulle sur $ $$\Delta^\vee_P- \Delta^\vee_Q$. Soit $P\in \pc(M)$ adjacent à $P_0$. De choses l'une soit $\deg(m_P)=\deg(m_{P_0})$ et dans ce cas on peut bien supposer que $\la_P=0$, soit 
$$ \deg(m_P)-\deg(m_{P_0})=x_\al \al^\vee$$
avec $x_\la>0$ et $\{\al^\vee\}=\Delta_{P_0}^\vee \cap (-\Delta_{P}^\vee)$ (par le lemme \ref{lem:pos}). L'égalité (\ref{eq:=}) n'est donc possible que si pour tout $\varpi\in   \hat{\Delta}_Q$ on a $\al^\vee\notin \Delta_{Q}^\vee$. Mais dans ce cas on a $P\subset Q$. En raisonnant par récurrence sur le cardinal de $\Delta_{P_0}^\vee \cap (-\Delta_{P}^\vee)$, on voit que pour tout $P\in \pc(M)$, on a soit $P\subset Q$ soit on peut supposer que $\la_P=0$. Cela prouve ainsi la deuxième assertion.

\medskip 
La deuxième assertion implique la troisième puisque pour $P\subset Q$ la projection de $\deg(m_P)$ sur $\ago_Q$ est égal à $\deg(m_Q)$. 

\medskip

Finalement montrons que la troisième assertion implique la première. Soit $Q\in \fc(M)$ tel que les conditions 3.(a) et 2.(b) soient satisfaites. 
Montrons que l'assertion 1 est alors vérifiée. Soit $X\in \overline{\mathcal{C}}_m$. Il s'agit de voir que $\|X-\xi\|\geq \|\varrho-\xi \|$. On a 
\begin{eqnarray}\label{eq:carres}
  \|X-\xi\|^2&=& \|\varrho-\xi \|^2 +2 \bg \xi-\varrho  , \varrho-X  \bd +\|X-\varrho \|^2\\ \nonumber
  &\geq &  \|\varrho-\xi \|^2 +2 \bg \xi-\varrho  , \varrho-X  \bd\\ \nonumber
&\geq&  \|\varrho-\xi \|^2 + \sum_{\varpi \in \hat{\Delta}_Q} 2y_\varpi \,\varpi(\varrho-X).
\end{eqnarray}
où les $y_\varpi$, définis en (\ref{eq:Yvarpi}), sont positifs. Par hypothèse, on a $\varrho\in -\deg(m_Q)+\ago_T^Q+\ago_G$ et 
Comme $X$ appartient à $\overline{\mathcal{C}}_m$, on a,  en particulier, 
$$X\in -\deg(m_Q)-\overline{^+\ago_Q}+\ago_T^Q +\ago_G.$$ 
Combiné à 3.(b), cela donne 
On a donc aussi
$$\varrho-X\in  \overline{^+\!\ago_Q}+\ago_T^Q+\ago_G.$$
Donc pour tout $\varpi \in \hat{\Delta}_Q$, on a $\varpi(\varrho-X)\geq 0$ et (\ref{eq:carres}) donne $ \|X-\xi\| \geq \|\varrho-\xi \|$ comme voulu.
\end{preuve}

La proposition précèdente implique immédiatement le corollaire suivant.

\begin{corollaire} \label{cor:HN} Plaçons-nous sous les hypothèses de la proposition \ref{prop:HN}. Soit $\varrho\in \bar{C}_m$ le $\xi$-point de Harder-Narasimhan de $m$. Il existe un unique $Q\in \fc(M)$ tel que $\varrho$ réalise la distance de $\xi$ au sous-espace affine $-\deg(m_Q)+\ago_T^Q+\ago_G$.  
\end{corollaire}

\end{paragr}

\section{Description adélique des fibres  de Hitchin}\label{sec:description}

\begin{paragr}\label{par:adel1} Soit $V$ un ensemble fini de points fermés de $C$ qui contient le point $\infty$ et le support du diviseur $D$. Le diviseur effectif $D$ s'écrit alors comme la somme formelle
$$D=\sum_{v\in V } d_v v$$
où, pour tout $v\in V$, l'entier $d_v$ est positif, nul hors du support de $D$ donc nul en $\infty$. Soit $C^V$ l'ouvert de $C$ complémentaire de $V$. Soit $k[C^V]$ l'algèbre des fonctions régulières sur $C^V$.
\end{paragr}

\begin{paragr}\label{par:adel2} Pour toute $k$-algèbre $A$, soit $C_A=C\times_k \Spec(A)$ et $C_A^V=C^V\times_k \Spec(A)$.  Soit $A[C^V]=k[C^V]\otimes_k A$. Pour tout $v\in V$, le complété de l'algèbre quasi-cohérente $\oc_{C_A}$ le long du diviseur $v\times_\kappa A$ s'identifie à l'anneau $A[[z_v]]$ via le choix d'une uniformisante $z_v$. Soit $A((z_v))$ l'anneau des séries formelles de Laurent en la variable $z_v$ à coefficients dans $A$ : c'est le localisé de $A[[z_v]]$ que l'on obtient en inversant $z_v$. On  a le diagramme commutatif suivant où les morphismes sont les morphismes évidents
  \begin{equation}
    \label{eq:diagdescript}
    \xymatrix{  \Spec(A((z_v))) \ar[d]_{i_{A,v}'}  \ar[r]^{j_{A,v}} & \Spec(A[[z_v]])   \ar[d]^{i_{A,v}} \\ 
 C_A^V  \ar[r]^{j_A}& C_A}
\end{equation}

Pour tout $k$-schéma $S$, on note $S((z_v))$, resp. $S[[z_v]]$ le foncteur qui, à toute $k$-algèbre $A$, associe l'ensemble des points $S(A((z_v)))$, resp. $S(A[[z_v]])$.
  
\end{paragr}

\begin{paragr}[Descente formelle à la Beauville-Laszlo et uniformisation.] ---\label{S:BL} Avec les notations des paragraphes précédents, on peut énoncer la proposition suivante.

 \begin{proposition} \label{prop:BL} Soit un  triplet $(X,(g_v)_{v\in V},t)$ qui vérifie les conditions suivantes :
   \begin{enumerate}
   \item $X$ est un  élément de $\ggo(A[C^V])$ ;
   \item pour tout $v\in V$, l'élément $g_v$ appartient à  $G((z_v))(A)$ et vérifie
$$\Ad(g_v^{-1})X\in z_v^{-d_v}\ggo[[z_v]](A) \ ;$$
en particulier $\chi(X)\in \car[[z_\infty]](A)$ ;
\item $t$ est un élément de $\tgo^{\reg}(A)$ dont la caractéristique $\chi(t)$ est égale à la réduction modulo $z_\infty$ de $\chi(X)$. 
   \end{enumerate}
Il existe un triplet $m=(\ec,\theta,t)\in \mc(\Spec(A))$ ainsi que des isomorphismes $G$-équivariants
$$\al_{A,v}\ :\ i_{A,v}^*\ec \to G\times_k A[[z_v]]$$ 
et 
$$\beta_A \ :\ j_A^*\ec\to    G\times_k A[C^V]$$ 
tels que
\begin{itemize}
\item l'automorphisme  $G$-équivariant de $G\times_k A((z_v))$ défini par $(i_{A,v}^{'*}) (\beta_{A})\circ (j^{*}_{A,v})(\al_{A,v}^{-1})$ est donné par la translation à gauche par $g_v\in G((z_v))(A)$ ;
\item $\beta_A$ induit un isomorphisme entre $\Ad(j_A^*\ec)$ et $\ggo\times_k A[C^V]$ qui envoie $j_A^*\theta$ sur $X$.
\end{itemize}
En outre, le triplet $(m,(\al_{A,v})_{v\in V},\beta_A)$ est uniquement déterminé à un unique isomorphisme près. Tout tel triplet s'obtient à partir d'un unique triplet $(X,(g_v)_{v\in V},t)$ comme ci-dessus.
  \end{proposition}

  \begin{preuve}
    C'est une conséquence du résultat de descente formelle de Beauville-Laszlo (cf. \cite{BL}).
  \end{preuve}

Introduisons alors la catégorie $Cat_V(A)$ dont les objets sont les triplets $(X,(g_v G[[z_v]](A))_{v\in V},t)$  qui vérifient les conditions 1 à 3 de la proposition \ref{prop:BL} (la condition 2 sur $g_v$ est clairement stable par translation par $G[[z_v]](A)$) et l'ensemble des morphismes entre deux objets $(X,(g_vG[[z_v]](A))_{v\in V},t)$ et $(X',(g'_vG[[z_v]](A))_{v\in V},t')$ est
\begin{itemize}
\item vide si $t\not= t'$ ;
\item l'ensemble des $\delta\in G(A[C^V])$ tels que $\Ad(g)X=X'$ et pour tout $v\in V$
$$\delta g_v G[[z_v]](A)= g_v'G[[z_v]](A).$$
\end{itemize}

  \begin{corollaire}\label{cor:BL}   La construction de la proposition \ref{prop:BL}  induit une équivalence de catégories  entre 
    \begin{itemize}
    \item la catégorie $Cat_V(A)$ ;
\item la sous-catégorie pleine de $\mc(\Spec(A))$ formée des triplets $(\ec,\theta,t)\in \mc(\Spec(A))$ pour lesquels le $G$-torseur $\ec$ est trivial sur $C_A^V$ et sur $\Spec(A[[z_v]])$ pour tout $v\in V$.
    \end{itemize}
  \end{corollaire}
\end{paragr}

\begin{paragr} \label{S:reducadel}Soit $K$ une extension de $k$. Soit $(a,t)\in \Ac_{G}(K)$. Soit $t_a\in \tgo^{\reg}[[z_\infty]](K)$ le relèvement de $t$ tel que $a$ et $\chi(t_a)$ coïncident sur $\Spec(K[[z_\infty]])$ (cf.  proposition \ref{prop:caractdesAM}). Soit $M\in \lc$ tel que $(a,t)\in \chi^M_G(\Ac_{M,\el}(K))$ (cf.  proposition \ref{prop:reunionAM}). Soit $m=(\ec,\theta,t)\in \mc(K)$ tel que  $f(m)=(a,t)$. Soit $\al_{K_v}$ et $\beta_K$ des trivialisations de $\ec$ respectivement sur $\Spec(K[[z_v]])$ et $C_K^V$. On suppose que ces trivialisations vérifient les conditions de la proposition \ref{prop:BL}. Soit $(X,(g_v)_{v\in V},t)$ le triplet qui s'en déduit par la proposition \ref{prop:BL}. On suppose de plus $X$ est un élément de $\mgo(K[C^V])$ tel que $X$ et $t_a$ sont conjugués sous $M((z_\infty))(K)$.

  \begin{proposition} \label{prop:reducadel} Soit $P\in \fc(M)$. Pour tout $v\in V$, soit $p_v\in P((z_v))(K)$ tel que 
$$g_v\in p_v G[[z_v]](K).$$
Alors il existe un triplet $(\ec_P,\theta_P,t)$ qui est une réduction de $m$ à $P$ et des trivialisations $\al_{K,v}^P$ et $\beta_K^P$ de $\ec_P$ respectivement sur $\Spec(K[[z_v]])$ et $C_K^V$ telles que 
\begin{itemize}
\item l'automorphisme  $P$-équivariant de $G\times_k K((z_v))$ défini par $(i_{K,v}^{'*}) (\beta_{K}^P)\circ (j^{*}_{K,v})(\al_{K,v}^P)^{-1}$ est donné par la translation à gauche par $p_v\in P((z_v))(K)$ ;
\item $\beta_K^P$ induit un isomorphisme entre $\Ad(j_K^*\ec_P)$ et $\pgo\times_k K[C^V]$ qui envoie $j_K^*\theta_P$ sur $X$.
\end{itemize}
    
  \end{proposition}

  \begin{preuve}
De la relation $\Ad(g_v^{-1})X\in z_v^{-d_v}\ggo[[z_v]](K)$ et du fait que $X\in \mgo(K[C^V])$, il vient
   $$\Ad(p_v^{-1})X\in z_v^{-d_v}\pgo[[z_v]](K).$$
En particulier, $\chi_P(X)\in \car_M[[z_\infty]](K)$. Comme on a supposé que $X$ et $t_a$ sont conjugués sous $M(F_\infty)$, on a $\chi_P(X)=\chi_P(t_a)$. Ainsi $(X,(p_v)_{v\in V}, t)$ est un triplet qui satisfait les conditions 1 à 3 de la proposition \ref{prop:BL} relatives à $P$. La proposition \ref{prop:BL} appliquée au groupe $P$ donne un triplet $(\ec_P,\theta_P,t)$. On vérifie que c'est bien une réduction à $P$ de $m$.
  \end{preuve}

\end{paragr}

\begin{paragr}[Fonctions $H_P$.] --- \label{S:HP}On continue avec les notations du paragraphe précédent.  

  \begin{definition} \label{def:HP}Soit $P\in \fc(T)$ et $v\in V$ et 
$$H_P\ : \ G((z_v))(K) \to \ago_{P}$$
l'application qui vérifie pour tout $\mu\in X^*(P)$ et $g\in G((z_v))(K)$, 
$$\mu(H_P(g))=-\val_v(\mu(p))$$
où $\val_v$ est la valuation usuelle sur   $K((z_v))$ et $p$ est un élément de $P((z_v))(K)$, uniquement  défini à translation à gauche près par un élément de $P[[z_v]](K)$, qui vérifie
$$g\in pG(K[[z_v]]).$$
L'existence d'un tel élément $p$ est donnée par la décomposition d'Iwasawa.\\

Plus généralement, pour une famille $(g_v)_{v\in V}$ d'éléments de $G((z_v))(K)$ on pose
$$H_P((g_v)_{v\in V})=\sum_{v\in V} H_P(g_v).$$
Cette définition garde un sens lorsque $V$ est infini pourvu que $g_v$ appartienne à $G[[z_v]]$ pour tout $v$ en dehors d'un ensemble fini. En particulier, la fonction $H_P$ est bien définie sur les points de $G$ à valeurs dans les adèles de $F$.
\end{definition}

Indiquons alors un corollaire à la proposition \ref{prop:reducadel}.

\begin{corollaire}\label{cor:reducadel}
  Avec les hypothèses et les notations du \S \ref{S:reducadel} et de  la proposition \ref{prop:reducadel}, on a l'égalité suivante pour tout $P\in \fc(M)$ et toute réduction $m_P$ de $m$ à $P$
$$\deg(m_P)=H_P( (g_v)_{v\in V}).$$
\end{corollaire}

\begin{preuve}
Soit $P\in\fc(M)$.
 La proposition   \ref{prop:reducadel} montre qu'il existe  une réduction $(\ec_P,\theta_P,t)$ de $m$ à $P$ tel que le $P$-torseur $\ec_P$ est défini par des données de recollement $p_v\in P((z_v))(K)$ où pour tout $v\in V$ on a $g_v\in p_vG[[z_v]]$. Donc pour tout caractère $\mu\in X^*(P)$, la famille $(\mu(p_v))_{v\in V}$ est une donnée de recollement pour le fibré en droites $\mu(\ec_P)$. On vérifie l'égalité 
$$\deg(\mu(\ec_P))=\sum_{v\in V} \val_v(\mu(p_v))=\mu(H_P( (g_v)_{v\in V})),$$
d'où le résultat.
\end{preuve}

\end{paragr}

\begin{paragr}[Description adélique des fibres de Hitchin.] --- \label{S:description-adelique}

On poursuit avec les notations des paragraphes précédents hormis $V$ qui désigne maintenant l'ensemble de \emph{tous} les points fermés de $C$. Soit $F$ le corps des fonctions de la courbe $C$. Soit  $\xi\in \ago_T$ et $(a,t)\in \Ac(k)$.  

\medskip

 \begin{definition}\label{def:Cata}
Soit $\mathcal{X}_{(a,t)}$ l'ensemble des  couples  $(X,(g_v)_{v\in V})$ qui vérifient 
\begin{enumerate}
\item $X$ est un  élément semi-simple $G$-régulier de $\ggo(F)$ dont la caractéristique $\chi(X)$ est égale à la restriction de $a$ au point générique de $C$ ;
   \item pour tout $v\in V$, l'élément $g_v$ appartient à  $G((z_v))(k)/G[[z_v]](k)$ et vérifie les trois conditions suivantes
     \begin{enumerate}
     \item pour tout $v\in V$ en dehors d'un ensemble fini $g_v$ est la classe triviale ;
     \item  $\Ad(g_v^{-1})X\in z_v^{-d_v}\ggo[[z_v]](k) $ ;
\item la caractéristique $\chi(t)$ est égale à la réduction modulo $z_\infty$ de $\chi(X)$ ; 
          \end{enumerate}
    \end{enumerate}
  \end{definition}
  
  \begin{remarque}
    La condition 2.(b) implique que la restriction de $\chi(X)$ à $\Spec(k((z_\infty))$ appartient à  $\car[[z_\infty]](k)$. Donc 2.(c) fait sens.
  \end{remarque}

Soit $t_a\in \tgo^{\reg}[[z_\infty]](k)$ l'unique relèvement de $t$ tel que $\chi(t_a)$ soit égal à la restriction de $a$ à $\Spec(k[[z_\infty]])$ (cf. proposition \ref{prop:caractdesAM}). Soit $M\in \lc$ tel que $(a,t)\in \chi^M_G(\Ac_{M,\el}(k))$ (cf. proposition \ref{prop:reunionAM}). 

\begin{definition}\label{def:Cata2}
Soit $\mathcal{X}^\xi_{(a,t)}$ l'ensemble des  couples  $(X,(g_v)_{v\in V})$ dans $\mathcal{X}_{(a,t)}$ qui vérifient 

\begin{enumerate}
\item (bis) $X$ est un  élément semi-simple elliptique $G$-régulier de $\mgo(F)$ qui est conjugué par un élément de $M((z_\infty))(k)$ à $t_a$ ;
\item \begin{enumerate}
\item[(d)]  la projection de $\xi$ sur $\ago_M^G$ appartient à la projection sur $\ago_M^G$ de l'enveloppe convexe des points 
$$-H_P((g_v)_{v\in V})=-\sum_{v\in V} H_P(g_v).$$
\end{enumerate}
\end{enumerate}
\end{definition}

Pour tout  $\delta\in G(F)$ et  tout couple $(X,(g_v)_{v\in V})$ dans $\mathcal{X}_{(a,t)}$  le couple $(\Ad(\delta)X,(\delta g_v)_{v\in V})$ appartient aussi à $\mathcal{X}_{(a,t)}$. On en déduit une action à gauche de $G(F)$ sur $\mathcal{X}_{(a,t)}$. L'action de $M(F)$ qui s'en déduit préserve $\mathcal{X}_{(a,t)}^\xi$. Les condition 1. (bis) et 2. (d) ci-dessus sont préservées par l'action de $M(F)$ : c'est trivial pour la première condition et la seconde résulte de l'égalité 
$$ H_P((mg_v)_{v\in V})=H_P((g_v)_{v\in V})$$
pour tout $m\in M(F)$. (Plus précisément, on a,  pour tout $\mu\in X^*(P)$ et tout $m\in M(F)$,  
$$ \mu(H_P((mg_v)_{v\in V}))=\mu(H_P((g_v)_{v\in V}))-\sum_{v\in V}   \val_v(\mu(m))$$
et la somme sur $V$ est nulle par la formule du produit.)

\medskip

Lorsqu'un groupe $G$ abstrait agit à gauche sur un ensemble $\mathcal{X}$, on note 
\begin{equation}
  \label{eq:groupoidequot}
  [G\back \mathcal{X}]
\end{equation}
le groupoïde quotient dont l'ensemble des objets est $\mathcal{X}$ et,  pour tous objets $x$ et $x'$ dans $\mathcal{X}'$, l'ensemble des morphismes de $x$ vers $x'$ est l'ensemble des $g\in G$ tels que $g\cdot x=x'$.

\begin{proposition}\label{prop:fibre}
  La construction de la proposition \ref{prop:BL} induit une équivalence de catégories entre le groupoïde quotient  $[M(F)\back \mathcal{X}^\xi_{(a,t)}]$   et la fibre $\overline{f^\xi}^{-1}(a,t)$.
\end{proposition}

\begin{preuve} Soit $m=(\ec,\theta,t)\in \overline{\mc^\xi}(k)$ tel que $f(m)=(a,t)$. Comme le corps $k$ est algébriquement clos, le $G$-torseur $\ec$ admet une trivialisation générique (cf. lemme \ref{lem:Sen}). On en fixe une et on en déduit une trivialisation générique de $\Ad_D(\ec)$, donc  un élément $X\in \ggo(F)$ image de $\theta$ par cette trivialisation. Soit $Y$ un élément semi-simple, $G$-régulier et elliptique dans $\mgo(F)$ tel que 
  \begin{itemize}
  \item   $\chi_G(Y)$ est égal à $a_\eta$ la restriction de $a$ au point générique $\eta$ de $C$ ;
  \item $Y$ est conjugué à $t_a$ par un élément de $M((z_\infty))(k)$.
 \end{itemize}
Un tel $Y$ existe (cf. corollaire \ref{cor:caractdesAM}). Il résulte du lemme \ref{lem:cjsep} que $X$ et $Y$ sont conjugués sous $G(F)$. Quitte à changer la trivialisation de $\ec$, on peut et on va supposer qu'on a  $X=Y$. 

Un choix de trivialisations de $\ec$ sur les voisinages formels $\Spec(k[[z_v]])$ donnent des éléments $g_v$ de $G((z_v))$ par la construction de la proposition \ref{prop:BL}. La classe de $g_v$ modulo $G[[z_v]]$ ne dépend pas de ce choix. On obtient ainsi un couple $(X,(g_v)_{v\in V})$ dans $\mathcal{X}_{(a,t)}$ qui vérifie en outre la condition 1.(bis) de la définition \ref{def:Cata}. En utilisant  la proposition \ref{prop:reducadel} et le fait que la trivialisation générique considérée s'étend à un ouvert de $C$, on voit qu'on a 
$$\deg(m_P)=\sum_{v\in V} H_P(g_v)$$
  pour tout $P\in \pc(M)$. Par la proposition \ref{prop:envpecvx}, la condition de $\xi$-semi-stabilité se traduit par la condition 2.(d) de la  définition \ref{def:Cata}. Ainsi $(X,(g_v)_{v\in V})$ appartient à $\mathcal{X}^\xi_{(a,t)}$. On laisse le soin au lecteur de vérifier que cette construction donne bien une équivalence de catégories comme annoncé.
  
\end{preuve}

\end{paragr}

\section{Critère valuatif : existence}\label{sec:existence}

\begin{paragr} Soit $\xi$ un élément \emph{quelconque} de $\ago_T$.   Le but de cette section est de prouver la partie existence du critère valuatif de propreté pour le morphisme $\overline{f^\xi}$ qui se résume au théorème suivant.

  \begin{theoreme}\label{thm:existence}
Soit $\kappa$ une extension  algébriquement close de $k$.  Soit $R$ un anneau de valuation discrète, complet et de corps résiduel $\kappa$. Soit $K$ le  corps des fractions de $R$ et $\overline{K}$ une clôture algébrique de $K$. 

Soit $(a,t)\in \Ac(R)$ et soit $m_K\in \overline{\mc^\xi}(K)$ tel que $f(m_K)=(a,t)$.

Il existe une extension finie $K'\subset \overline{K}$ de $K$ et $m\in \overline{\mc^\xi}(R')$, où $R'$ est la clôture intégrale de $R$ dans $K'$, tels que 
\begin{enumerate}
\item  l'image de $m$ dans $ \overline{\mc^\xi}(K')$ est isomorphe à celle de $m_K$ ;
\item  $f(m)=(a,t)$.
\end{enumerate}
   \end{theoreme}

La démonstration du théorème \ref{thm:existence} va nous occuper jusqu'à la fin de cette section.
\end{paragr}

\begin{paragr}[Quelques mots sur la démonstration.] --- Dans \cite{Langton} (\S3 théorème), Langton a montré  que le champ des fibrés vectoriels semi-stables sur une courbe projective et lisse satisfait la partie existence du critère valuatif de propreté. En adaptant les arguments de Langton, Nitsure (cf. \cite{Nitsure} section 6) a démontré le théorème \ref{thm:existence} dans le cas $G=GL(n)$ et $\xi=0$. Faltings a réussi à étendre cette démonstration à tout groupe $G$ réductif dans le cas d'un corps de base de caractéristique $0$ et pour $\xi=0$ (cf. \cite{Faltings} théorème II.4) . Sa démonstration ne semble cependant pas s'étendre au cas d'un corps de base de caractéristique positive. Dans \cite{Heinloth}, Heinloth a démontré l'analogue du théorème de Langton pour le champ des $G$-torseurs sur une courbe projective et lisse sur un corps de base caractéristique nulle ou ``pas trop petite''. Notre démonstration reprend en partie les arguments de Heinloth, eux-mêmes inspirés par ceux de Langton. \\

Esquissons les grandes lignes de notre démonstration. Les notations sont celles du  théorème \ref{thm:existence}. Soit $C_K$ et $C_\kappa$ respectivement les fibres générique et spéciale de la  courbe $C\times_k R$ sur $R$. Soit $(a,t)\in \Ac(R)$ sur cette courbe et $m_K\in \overline{\mc^\xi}(K)$ un triplet de Hitchin $\xi$-semi-stable au-dessus de $C_K$ de caractéristique $(a,t)$. Les énoncés d'existence qui suivent exigent parfois de remplacer $K$ par une extension finie, ce qu'on ne rappelera pas à chaque fois. Un tel triplet $m_K$ s'étend en un triplet $m_R$ sur toute la courbe $C_R$ (cf. \S\ref{S:existpourf}). La fibre spéciale $m_\kappa$ du prolongement n'est pas nécessairement $\xi$-semi-stable. Cependant, on peut choisir le prolongement $m_R$ de sorte  que la fibre spéciale $m_\kappa$ soit le moins ''$\xi$-instable'' possible, ce par quoi on entend que la distance $d$ de $\xi$ au convexe  $\overline{\mathcal{C}}_{m_\kappa}$ défini à la section \ref{sec:cvx}) est minimale (cf. \S\ref{S:d}). Si elle est nulle, $m_\kappa$ est $\xi$-semi-stable. Sinon, il existe $Q$ un sous-groupe parabolique maximal de $G$ tel que la face de $\overline{\mathcal{C}}_{m_\kappa}$ associée à $Q$ contienne le point où cette distance est réalisée (cf. \S\ref{S:YHN}). Le prolongement $m_R$ est donné par des ``données de recollement'' (cf. \S\ref{S:trivadm}, le résultat de \cite{DS} y joue un rôle crucial et nécessite la réduction au cas semi-simple du \S\ref{S:reducsemisimple}). En conjuguant ces données par un $K$-point convenable du centre de $Q$, on obtient un autre triplet $m'_R$ qui est en fibre générique isomorphe à $m_K$ et dont la face de $\overline{\mathcal{C}}_{m_\kappa'}$ associé au sous-groupe parabolique $Q^-$ opposé à $Q$ est égale à la face de $\overline{\mathcal{C}}_{m_\kappa}$ associé à $Q$. Par minimalité de $d$, cela force le convexe $\overline{\mathcal{C}}_{m_\kappa'}$ à être dégénéré : les faces de  $\overline{\mathcal{C}}_{m_\kappa'}$ associés à $Q$ et $Q^-$ sont dans le même hyperplan affine. Mais cela implique  que la réduction à $Q$ de $m_\kappa$ se relève en une réduction à $Q$ du triplet $m_{R/(\pi)^2}\in\mc(R/(\pi)^2)$  déduit de $m_R$ par le changement de base $R\to R/(\pi)^2$ où $(\pi)\subset R$ est l'idéal maximal. On peut réitérer ce processus en changeant de $K$-point du centre de $Q$ (cf. \S\ref{S:cstr-de-triv}). On obtient alors une réduction de $m$ à $Q$ sur le complété formel  $\varinjlim_{n} C\times_k R/\pi^n$. Par un théorème de Grothendieck, cette réduction  à $Q$ s'algébrise et cela contredit alors la $\xi$-stabilité de $m_K$ (cf. \S\ref{S:lacontradiction}).

\end{paragr}

\begin{paragr}\label{S:anneauR} On fixe désormais des objets $R$, $\kappa$, $K$ et $\overline{K}$ comme dans l'énoncé du théorème \ref{thm:existence}. Par le choix d'une uniformisante $\pi$, on identifie $R$ à l'anneau $\kappa[[\pi]]$. 
  \end{paragr}

\begin{paragr}[Le théorème \ref{thm:existence} pour le champ $\mc$.] --- \label{S:existpourf} Montrons tout d'abord que le théorème \ref{thm:existence} vaut lorsqu'on remplace le champ $\overline{\mc^\xi}$ par le champ $\mc$ tout entier. C'est bien connu. Nous ne donnons un argument que par souci d'exhaustivité. On reprend les hypothèses du théorème \ref{thm:existence} excepté qu'on ne suppose pas le point $m_K\in \mc(K)$ $\xi$-semi-stable. Le point $m_K$ s'explicite comme un triplet $(\ec_K,\theta_K,t)$. Soit $\eta$ le point générique de $C_K=C\times_k K$. Quitte à remplacer $K$ par une extension finie, on peut et on va supposer que le torseur $\ec_K$ est trivial au point $\eta$ (cf. lemme \ref{lem:Sen}). Fixons donc une telle trivialisation : on en déduit une trivialisation de $\Ad_D(\ec)$ en $\eta$ et donc un élément $Y\in \ggo(F)$, où $F$ est le corps des fonctions de $C_K$, qui est l'image de $\theta$ par cette dernière trivialisation. Notons que $\chi(Y)$ est égal à la restriction $a_\eta$ de $a$ à $\eta$.

Soit $B$ l'anneau local du point fermé de $C_R=C\times_k \Spec(R)$ défini par l'idéal engendré par $\pi$. C'est un anneau de valuation discrète de corps des fractions $F$. En prenant une section de Kostant (par exemple), on voit qu'il existe $X\in \ggo(B)$ telle que $\chi(X)=a_{|\Spec(B)}$. Puisque $\chi(X)=\chi(Y)$, le lemme \ref{lem:cjsep} entraîne que, quitte à remplacer $K$ par une extension finie, on peut supposer qu'il existe $g\in G(F)$ tel qu'on ait $\Ad(g)X=Y$. On peut alors avec cet élément $g\in G(F)$ recoller le torseur $\ec_K$ et la section $\theta$ avec le couple $(G\times_k B,X)$ sur $\Spec(B)$. On obtient ainsi un couple $(\ec_U,\theta_U)$ qui prolonge $(\ec_K,\theta_K)$ à un ouvert $U$ de $C_R$ qui contient tous les points de codimension $\leq 1$. Mais tout $G$-torseur sur un tel ouvert $U$ se prolonge de manière unique en un $G$-torseur sur $C_R$ tout entier (cf. \cite{CTS} théorème 6.13). On obtient un $G$-torseur $\ec$ sur $C_R$ qui prolonge $\ec_K$. Le fibré $\Ad_D(\ec)$ possède une section sur $U$ à savoir $\theta_U$. Mais comme $U$ contient tous les points de  codimension $\leq 1$ et que  $\Ad_D(\ec)$ est affine sur $C$, cette section se prolonge de manière unique sur $C_R$ (cf. \cite{EGAIV4} corollaire 20.4.12). On a ainsi obtenu un couple $(\ec,\theta)$. Les points $\chi(\theta)(\infty_R)$ et $\chi(t)$ sont deux points de $\car(R)$. Ils sont en fait égaux puisque leurs images dans $\car(K)$ sont égales. Ainsi le couple  $(\ec,\theta)$ se complète en un triplet $(\ec,\theta,t)\in \mc(R)$ d'où le théorème  \ref{thm:existence} pour $\mc$.
  
\end{paragr}

\begin{paragr}[Réduction au cas semi-simple.] --- \label{S:reducsemisimple} Dans ce paragraphe, on montre que le théorème  \ref{thm:existence} pour tout groupe semi-simple implique le théorème  \ref{thm:existence} pour tout groupe réductif. 

Soit $G$ un groupe réductif connexe sur $k$ d'algèbre de Lie $\ggo$ et $T$ un sous-tore maximal.  Soit $G_{\der}$ le groupe dérivé de $G$ et $\ggo_{\der}$ son algèbre de Lie. La décomposition 
$$\ggo=\ggo_{\der} \oplus \zgo $$
induit une décomposition

\begin{equation}
  \label{eq:carG1}
  \car_G=\car_{G_{\der}}\oplus \zgo
\end{equation}

Soit $A_G$ le centre connexe de $G$ et $\zgo$ son algèbre de Lie. Soit $G'=G/A_G$. C'est un groupe semi-simple d'algèbre de Lie  $\ggo'=\ggo/\zgo$ et qui admet $T'=T/A_G$ comme sous-tore maximal. Le morphisme évident $\ggo \to \ggo'$ induit des isomorphismes $\ggo_{\der}\to \ggo'$  et 
\begin{equation} \label{eq:carG2}
  \car_{G_{\der}}\simeq \car_{G'}.
\end{equation}
La  projection de $  \car_G$ sur $\car_{G_{\der}}$ définie   par (\ref{eq:carG1}) composée avec l'isomorphisme (\ref{eq:carG2}) et le morphisme évident $\tgo \to \tgo'=\tgo/\zgo$ induisent un morphisme $\Ac_G \to \Ac_{G'}$. La bijection évidente $\lc^G(T) \to \lc^{G'}(T')$ est compatible à la décomposition de la proposition \ref{prop:reunionAM}.

Soit $S$ un $k$-schéma affine et $m=(\ec,\theta,t)\in\mc_G(S)$.  Soit $\ec'$ le $G'$-torseur obtenu lorsqu'on pousse $\ec$ par le morphisme $G\to G'$. La section $\theta$ poussée par le morphisme $\ggo\to\ggo'$ fournit  une section $\theta'$ de $\Ad_D(\ec')$. Soit $t'$ l'image de $t$ par le morphisme $\tgo\to\tgo'=\tgo/\zgo$. Le triplet $m'=(\ec',\theta',t')$ appartient à $\mc_{G'}(S)$. On obtient ainsi un morphisme $\mc_G\to \mc_{G'}$ qui à $m$ associe $m'$. Il s'inscrit dans le diagramme commutatif 

$$\xymatrix{\mc_G \ar[d]^{f_G} \ar[r] & \mc_{G'} \ar[d]^{f_{G'}} \\ \Ac_G  \ar[r] & \Ac_{G'}}$$

En utilisant le morphisme $P \to P/A_G$, il est clair que toute réduction de $m$ à un sous-groupe parabolique $P$ donne naturellement une réduction  de $m'$ à $P/A_G$. Si $S$ est le spectre d'un corps, on obtient ainsi une bijection naturelle entre les réductions de $m$ et celles de $m'$ compatible à la bijection $\fc^G(T) \to \fc^{G'}(T')$. Le morphisme surjectif évident $\ago_T \to \ago_{T'}$ induit un isomorphisme $\ago_T^G\simeq \ago_{T'}$. Il est alors évident sur les définitions que le point $m$ est $\xi$-semi-stable si et seulement si $m'$ est $\xi^G$-semi-stable, où $\xi^G$ est la projection de $\xi$ sur $\ago_T^G$.

Revenons à la situation du théorème  \ref{thm:existence}. Soit $(a,t)\in \Ac_G(R)$ et $m_K= (\ec_K,\theta_K,t)\in \overline{\mc^\xi_G}(K)$ tel que $f_G(m_K)=(a,t)$. Soit $m_K'= (\ec_K',\theta_K',t')\in \overline{\mc^\xi_{G'}}(K)$ image de $m_K$ par le morphisme $\mc_G\to \mc_{G'}$. Soit $(a',t')=f_{G'}(m'_K)$. Le point $(a',t')$ appartient à $\Ac_{G'}(R)$ puisque c'est l'image de $(a,t)$.  Le théorème  \ref{thm:existence} pour le groupe semi-simple $G'$ implique que, quitte à changer $K$ par une extension finie $K'$ dans $\overline{K}$ et $R$ par sa clôture intégrale dans $K'$, on peut supposer que $m_K'$ se prolonge en un triplet $m'=(\ec',\theta',t')\in   \overline{\mc^\xi_{G'}}(R)$ qui prolonge $m_K$. Quitte de nouveau à changer $K$ et $R$ comme ci-dessus, on peut supposer que le $G'$-torseur $\ec'$ est trivial sur le spectre de $B$ l'anneau local du point de $C_R$ défini par $(\pi)$ (cf. \cite{DS} théorème 2). Fixons une telle trivialisation. On en déduit une trivialisation de $\Ad_D(\ec')$ sur $\Spec(B)$ ; celle-ci  envoie $\theta'$ sur un point $Y'\in \ggo'(B)\simeq\ggo_{der}(B)$ tel que $\chi_{G'}(Y)=a'_{|\Spec(B)}$. En écrivant $a_{|\Spec(B)}=a'_{|\Spec(B)}\oplus Z$ avec $Z\in \zgo(B)$ suivant (\ref{eq:carG1}) et (\ref{eq:carG2}), on obtient un point $Y=Y'+Z\in \ggo(B)$. Toute trivialisation de $\ec_K$ sur un ouvert $U$ de $C_K$ donne des trivialisations  de $\ec'_K$ et $\Ad_D(\ec'_K)$ sur le même ouvert ainsi que des points $X$ et $X'$  sections respectives de $\ggo \times_k U$ et $\ggo'\times_k U$ qui se déduisent de $\theta$ et $\theta'$. La condition de recollement des trivialisations de $(\ec',\theta')$ sur $\Spec(B)$ et $U$ se traduit par l'existence d'un élément $g'\in G'(F)$ tel que $\Ad(g')X'=Y'$ où $F$ est le corps des fonctions de $C_R$. Tout relèvement $g\in G(F)$ de $g'$ vérifie $\Ad(g)X=Y$ par construction de $Y$. Or l'obstruction à un tel relèvement vit dans $H^1(F,A_G)$. Quitte à remplacer $K$ par une extension, on peut supposer que cet ensemble est réduit à la classe triviale (cf. \cite{Serre} chap. III \S2) et donc qu'un tel relèvement existe. D'un tel relèvement, on déduit, comme au paragraphe \ref{S:existpourf}, un prolongement  $m=(\ec,\theta,t)\in \mc_G(R)$ de $m_K$. Par construction, l'image de $(\ec,\theta,t)$ dans $\mc_{G'}(R)$ coïncide avec $m'$ sur un ouvert de $C_R$ qui contient tous les points de codimension $\leq 1$. Par un argument déjà évoqué au  paragraphe \ref{S:existpourf}, l'image de $(\ec,\theta,t)$ est donc (isomorphe à) $m'$. Il s'ensuit qu'on  a  $m \in \overline{\mc^\xi_{G}}(R)$ comme voulu.
 \end{paragr}

 \begin{paragr}[Les points $m_K$, $(a,t)$ et $(\bar{a},\bar{t})$.]  --- Désormais et ce jusqu'à la fin de la preuve du théorème \ref{thm:existence}, on suppose que $G$ est \emph{semi-simple}.  On se place sous les hypothèses du théorème \ref{thm:existence}. Soit $(a,t)\in \Ac(R)$ et soit $m_K\in \overline{\mc^\xi}(K)$ tel que $f(m_K)=(a,t)$. Soit $(\bar{a},\bar{t})\in \Ac(\kappa)$ le point déduit de $(a,t)$ par le changement de base $\Spec(\kappa)\to\Spec(R)$. 
   
 \end{paragr}

 \begin{paragr}[Le réel $d$ et le point $m=(\ec,\theta,t)\in \mc(R)$. ] ---\label{S:d}   Pour toute extension  $K'\subset\overline{K}$ de degré fini sur $K$, soit $m_{K'}\in \overline{\mc^\xi}(K')$ déduit de $m_K$ par le changement de base $\Spec(K')\to \Spec(K)$. Soit $R'$ la clôture intégrale $R$ dans $K'$. Soit $m_{R'}\in \mc(R')$ un prolongement de $m_{K'}$ à $C\times_k R'$. L'anneau $R'$ est un anneau de valuation discrète de corps résiduel $\kappa$. Soit $m_{\kappa}$ le point de $\mc(\kappa)$ déduit de $m_{R'}$ par le changement de base $\Spec(\kappa)\to \Spec(R')$. Soit 
$$d(\overline{\mathcal{C}}_{m_{\kappa}},\xi)$$
la distance de $\xi$ au convexe $\overline{\mathcal{C}}_{m_{\kappa}}$ (cf. \ref{S:lescvxes}).

 \begin{lemme}\label{lem:d}
   L'ensemble des réels $d(\overline{\mathcal{C}}_{m_{\kappa}},\xi)$ pour tous $R'$ et $m_{R'}$ comme ci-dessus admet un plus petit élément.
 \end{lemme}

\begin{preuve}  Cet ensemble est non vide (cf. \S \ref{S:existpourf}) et  minoré par $0$. Il est même discret et fermé dans $\ago_T$ puisqu'inclus dans l'ensemble des distances de $\xi$ aux sous-espaces $X+\ago_T^Q+\ago_G$ pour $X\in X_*(Q)$ et  $Q\in \fc$ (cf. corollaire \ref{cor:HN}). Il contient donc sa borne inférieure.
 \end{preuve}

 Soit $d\geq 0$ le plus petit élément de l'ensemble considéré dans le lemme \ref{lem:d}.  Quitte à remplacer $K$ par une extension finie, on peut et on va supposer qu'il existe  
$$m=(\ec,\theta,t)\in \mc(R)$$ 
dont l'image dans  $\mc(K)$ est le point $m_K$ de départ et telle que la distance de $\xi$ au convexe $\overline{\mathcal{C}}_{m_\kappa}$ soit $d$, où $m_\kappa\in \mc(\kappa)$ est l'image de $m$. Dans toute la suite, on \emph{fixe} un tel $m$. 

Comme $m$ est $\xi$-semi-stable sur $C_K$, il sera $\xi$-semi-stable si et seulement s'il l'est aussi sur $C_\kappa$. On a donc les équivalences suivantes :
$$m\in \overline{\mc^{\xi}}(R) \Leftrightarrow m_\kappa \in \overline{\mc^{\xi}}(\kappa)  \Leftrightarrow d=0.$$
Ainsi le théorème \ref{thm:existence} est vrai pour le point $m_K$ si et seulement si $d=0$. Dans la suite, on suppose qu'on a  
$$d>0$$ 
et on va montrer qu'on aboutit à une contradiction. 
 \end{paragr}

\begin{paragr}[Le sous-groupe de Levi $L$ et les sous-groupes paraboliques $Q_0$ et $Q$.] ---  \label{S:YHN} Soit $L\in \lc$ tel que $f(m_{\kappa})\in \chi^L_G(\Ac_{L,\el}(\kappa))$ (cf. proposition \ref{prop:reunionAM}). Soit $\varrho$ le $\xi$-point de Harder-Narasimhan de $m_\kappa$. D'après la proposition \ref{prop:HN}, il existe un sous-groupe parabolique $Q_0\in \fc(L)$ tel qu'on ait 

\begin{itemize}
    \item $\xi-\varrho\in \ago_{Q_0}^+\cap \ago_T^G$ ;
    \item $\varrho$ appartient au sous-espace affine 
$$-\deg(m_{\kappa,Q_0})+\ago_T^{Q_0}.$$
    \end{itemize}
où $m_{\kappa,Q_0}$ est une réduction de $m_\kappa$ à $Q_0$ (ici $\ago_G=\{0\}$ puisque $G$ est semi-simple). Par hypothèse, la distance de $\varrho$ à $\overline{\mathcal{C}}_{m_\kappa}$ n'est pas nulle ; le sous-groupe parabolique $Q_0$ est donc propre. Dans toute la suite on fixe $Q\in \fc(Q_0)$ maximal parmi les sous-groupes paraboliques propres de $G$ qui contiennent $Q_0$. Soit $Q^-$ le sous-groupe parabolique opposé à $Q$, au sens où $Q^- \cap Q$ est l'unique facteur de Levi de $Q$ qui contient $T$.
\end{paragr}

\begin{paragr}[Quelques notations : $V$, $i_{R,v}$, $j_{R,v}$ etc.] ---  Soit $V$ est un ensemble fini de points fermés de $C_\kappa$ qui contient le support du diviseur 
$$D_\kappa=D\times_k \kappa=\sum_{v\in V} d_v \, v$$
ainsi que le point $\infty_\kappa$. Soit $C_\kappa^V=C_\kappa-V$ et $\kappa[C^V]$ l'algèbre des fonctions régulières sur $C_\kappa^V$. On reprend les notations de la section \ref{sec:description} en particulier le diagramme (\ref{eq:diagdescript}) pour la courbe $C_\kappa=C\times_k \kappa$ où sont définis les morphismes $j_A$, $i_{A,v}$ etc. pour tout point $v\in V$ et toute $\kappa$-algèbre.
\end{paragr}

\begin{paragr}[Le point $t_a$.] --- Soit $t_a$ l'unique point de $\tgo^{G\textrm{-}\reg}[[z_\infty]](R)$ qui vérifie les deux assertions suivantes :
  \begin{itemize}
  \item la réduction modulo $z_\infty$ de $t_a$ est le point $t\in \tgo^{G\textrm{-}\reg}(R)$ ;
  \item $\chi_G(t_a)=\chi_D(\theta)\circ i_{R,\infty}$.
  \end{itemize}
   L'existence et l'unicité de $t_a$ résulte du fait que $\chi$  induit un morphisme étale de $\tgo^{G\textrm{-}\reg}$ sur $\cgo^{G\textrm{-}\reg}$. 
Soit 
$$\bar{t}_a\in \tgo^{G\textrm{-}\reg}[[z_\infty]](\kappa)$$
la réduction modulo $\pi$ de $t_\infty$. Notons qu'on a $\bar{t}_a=\bar{t}_{\bar{a}}$ où $t_{\bar{a}}$ est l'unique point de   $\tgo^{G\textrm{-}\reg}[[z_\infty]](\kappa)$ de réduction $\bar{t}$ et de caractéristique $\bar{a}$.
\end{paragr}

\begin{paragr}[Trivialisation de $m$ sur un ouvert de la fibre spéciale.] --- Soit $m_\kappa=(\ec_\kappa,\theta_\kappa,\bar{t})\in \mc(\kappa)$ le point déduit de $m$ par le changement de base $\Spec(\kappa)\to \Spec(R)$. On a le lemme suivant.

  \begin{lemme} \label{lem:trivouvert} Quitte à  rajouter un nombre fini de points à $V$, on est dans la situation suivante. Il existe une trivialisation du $G$-torseur $j_\kappa^*(\ec_\kappa)$ sur $C_\kappa^V$

$$  \beta_\kappa\ :\ j_\kappa^*(\ec_\kappa)\to G\times_{k} C_\kappa^V.$$
telle que l'isomorphisme de schémas en groupes qui s'en déduit 
$$\iota_\kappa \ :\ \Aut_G(j_\kappa^*(\ec_\kappa)) \to G\times_{k}C_\kappa^V$$
satisfasse la condition : l'élément $d\iota_\kappa(j_\kappa^*\theta)$, où $d\iota_\kappa$ est le morphisme dérivé de $\iota_\kappa$, vérifie :
\begin{enumerate}
\item $d\iota_\kappa(j_\kappa^*\theta)$ appartient à $\lgo(\kappa[C_\kappa^V])$ ;
\item $d\iota_\kappa(j_\kappa^*\theta)$ est conjugué à $\bar{t}_\infty$ par un élément de $l \in L((z_\infty))(\kappa)$.
\end{enumerate}
  \end{lemme}

  \begin{preuve} Soit  $F$ est le corps des fonctions de $C_\kappa$. Comme $\kappa$ est algébriquement clos, le $G$-torseur $\ec_\kappa$ est trivial au point générique de $C_\kappa$ (cf. lemme \ref{lem:Sen}). Par le choix d'une trivialisation, on obtient aussi une trivialisation générique de $\Ad_D(\ec_\kappa)$ qui envoie $\theta$ sur un certain élément $X\in \ggo(F)$. Soit $(\bar{a},\bar{t})=f(m_\kappa)$. Cet élément  appartient à $\Ac_{L,\el}(\kappa)$ d'après la définition de $L$ (cf. \S  \ref{S:YHN}). D'après le corollaire  \ref{cor:caractdesAM}, il existe $X'$ semi-simple et $G$-régulier dans $\lgo(F)$ et $l\in L((z_\infty))(\kappa)$ tel que

  \begin{enumerate}
  \item[(A)] $\chi_G(X')=\bar{a}_{|\Spec(F)}$ ;
  \item[(B)] $\Ad(l)X'=\bar{t}_\infty$.
  \end{enumerate}

Notons qu'on a $\chi_G(X)=\chi_G(X')$ par définition de $\bar{a}$. Le lemme \ref{lem:cjsep} implique que, quitte à changer la trivialisation générique de $\ec_\kappa$, on peut supposer que $X$ vérifie les conditions (A) et (B) ci-dessus. 

On obtient alors le lemme puisque toute trivialisation générique de $\ec_\kappa$ s'étend à  un ouvert de $C_\kappa$ et que l'élément $X$ lui aussi se prolonge à un ouvert.

\end{preuve}

\end{paragr}

\begin{paragr}[Trivialisation de $m$ sur des disques formels de la fibre spéciale.] --- Dans la suite, on fixe une trivialisation 
\begin{equation}
  \label{eq:betakappa}
  \beta_\kappa\ :\ j_\kappa^*(\ec_\kappa)\to G\times_{k} C_\kappa^V.
\end{equation}
qui vérifie les conditions du lemme \ref{lem:trivouvert}.

Pour tout $v\in V$ soit 
\begin{equation}
  \label{eq:alphakappav}
  \al_{\kappa,v} \ : \ i_{\kappa,v}^*(\ec_\kappa) \to G\times_k \kappa[[z_v]].
\end{equation}
un isomorphisme $G$-équivariant. Un tel isomorphisme existe. En effet, cela revient à se donner une section du torseur $i_{\kappa,v}^*(\ec_\kappa)$ au-dessus de $\Spec(\kappa[[z_v]])$. Une telle section existe au moins au-dessus du point spécial et, par lissité, elle se prolonge à  $\Spec(\kappa[[z_v]])$.

\begin{lemme} \label{lem:triv_en_v}On peut choisir les trivialisations $ \al_{\kappa,v}$ de sorte le triplet
\begin{equation}
  \label{eq:tripletadele}
  (\bar{X},(\bar{g}_v)_{v\in V},\bar{t})
\end{equation}
associé au point $m_\kappa=(\ec_\kappa,\theta_\kappa,\bar{t})\in \mc(\kappa)$ et aux trivialisations $(\al_{\kappa,v},\beta_\kappa)$ par la bijection de la proposition \ref{prop:BL} satisfasse les conditions suivantes :
  \begin{enumerate}
  \item  $\gb_v\in Q((z_v))(\kappa)$ ;
  \item $\al_{\kappa,\infty}(\ i_{\kappa,\infty}^*\theta )=\bar{t}_a$ ;
  \item $\gb_\infty \in L((z_\infty))(\kappa)$.
  \end{enumerate}
  
\end{lemme}

\begin{preuve} Quel que soit le choix de   $ \al_{\kappa,v}$, on a 
$$\gb_v\in G((z_v))(\kappa).$$
 La décomposition d'Iwasawa
$$G((z_v))(\kappa)=Q((z_v))(\kappa) \cdot G[[z_v]](\kappa)$$
montre qu'on peut toujours modifier  $ \al_{\kappa,v}$ de sorte qu'on ait
\begin{equation}
  \label{eq:gbinQ}
  \gb_v\in Q((z_v))(\kappa)
\end{equation}
d'où l'assertion 1. 

On a 
$$\Ad(\gb_\infty)^{-1}\bar{X}\in \ggo[[z_\infty]](\kappa).$$
L'élément $\chi(\bar{X})$, après changement de base à $\Spec(k[[z_\infty]])$, appartient à $\car[[z_\infty]](\kappa)$ et sa réduction modulo $z_\infty$ coïncide avec celle de  $\chi(\bar{t}_\infty)$. Comme $\chi$ est étale au-dessus de $\car^{\reg}$, on a 
$$\chi(\Ad(\gb_\infty)^{-1}\bar{X})=\chi(\bar{X})=\chi(\bar{t}_\infty).$$ 
Le lemme \ref{lem:section} ci-dessous montre qu'on peut modifier $ \al_{\kappa,\infty}$ de sorte qu'on ait
\begin{equation}
   \label{eq:cj_t_X}
   \Ad(\gb_\infty)^{-1}\bar{X}=\bar{t}_\infty.
 \end{equation}
Mais d'après le lemme \ref{lem:trivouvert}, $\bar{X}$ et $t_\infty$ sont conjugués sous $L((z_\infty))(\kappa)$. Il s'ensuit qu'on a $\gb_\infty\in L((z_\infty))(\kappa)$.
\end{preuve}

\begin{lemme}\label{lem:section}Soit  $Y \in \ggo[[z_\infty]](R)$ tel que $\chi_G(Y)=\chi_G(t_a)$. Il existe  $g\in G[[z_\infty]](R)$ tel que $\Ad(g)Y=t_a$. \\
Le même énoncé vaut lorsqu'on remplace $R$ par $\kappa$ et $t_a$ par $\bar{t}_a$.
\end{lemme}

\begin{preuve}
  Soit $\tc$ le $R[[z_\infty]]$-schéma  défini pour toute $R[[z_\infty]]$-algèbre $B$ par
$$\tc(B)=\{h\in G(B)\  | \ \Ad(h)(Y)=t_a\}.$$

  Montrons que le schéma $\tc$ est lisse sur $\Spec(R[[z_\infty]])$. Notons tout d'abord que les  centralisateurs de $Y$ et $t_a$ dans $G\times_k R[[z_\infty]]$  sont des sous-schémas en tores :  cela résulte du fait que $\theta_\infty\in \tgo^{G\textrm{-}\reg}(R[[z_\infty]])$ et $\chi_G(t_a)=\chi_G(Y)$. Par conséquent, localement pour la topologie étale, ces centralisateurs sont conjugués. Comme $Y$ et $t_a$ ont même caractéristique, ces éléments sont localement conjugués pour la topologie étale. Ainsi le  schéma $\tc$ a localement des sections pour la topologie étale. C'est donc un torseur sous le schéma en tores $T\times_k R[[z_\infty]]$. Il est donc lisse.\\
D'après le lemme (\ref{lem:cjsep}),  le schéma $\tc\times_{R[[z_\infty]]} \kappa$ possède des sections. Par lissité de $\tc$ sur $R[[z_\infty]]$ et par complétude de $R[[z_\infty]]\simeq \kappa[[\pi,z_\infty]]$, cette  section se relève en une section $h\in G[[z_\infty]](R)$ de $\tc$.  

La seconde assertion se démontre de la même façon.

\end{preuve}

\end{paragr}

\begin{paragr}[Trivialisation admissibles de $(\ec,\theta,t)$.] --- \label{S:trivadm} 
Désormais on fixe des trivialisations $(\al_{\kappa,v},\beta_\kappa)$ de $m_\kappa$ qui satisfont les lemmes \ref{lem:trivouvert} et \ref{lem:triv_en_v}. Soit $C^V_R=C^V_\kappa \times_\kappa R$. Soit $R[C^V]=\kappa[C^V]\otimes_\kappa R$. 

  \begin{definition} \label{def:trivadm}Une \emph{trivialisation admissible} de $m=(\ec,\theta,t)$ est la donnée d'isomorphismes $G$-équivariants
\begin{equation}
  \label{eq:betaR}
  \beta_R\ :\ j_R^*\ec\to G\times_{k} R[C^V].
\end{equation}
et pour tout $v\in V$
\begin{equation}
  \label{eq:alphaRv}
  \al_{R_v} \ : \ i_{R,v}^*\ec \to G\times_k R[[z_v]]
\end{equation}
  qui satisfont les conditions suivantes
  \begin{enumerate}
  \item l'isomorphisme $\beta_R$ se spécialise en l'isomorphisme $\beta_\kappa$ de (\ref{eq:betakappa}) ;
  \item  l'isomorphisme $ \al_{R,v}$  se spécialise en l'isomorphisme $\al_{\kappa_v}$ de (\ref{eq:alphakappav}) ;
  \item $\al_{R,\infty}(i_{R,\infty}^*\theta)=t_a$.

  \end{enumerate}

\end{definition}

\begin{remarque}
  Soit $(X,(g_v)_{v\in V},t)$ le triplet associé à $m\in \mc(R)$ et aux trivialisations \ref{eq:betaR} et \ref{eq:alphaRv} par la proposition \ref{prop:BL}. Alors ces trivialisations forment une trivialisation admissible de $m$ si et seulement si
  \begin{itemize}
  \item la réduction modulo $\pi$ du triplet $(X,(g_v)_{v\in V},t)$ est le triplet $(\bar{X},(\bar{g}_v)_{v\in V},\bar{t})$ défini en (\ref{eq:tripletadele}).
  \item $\Ad(g_\infty)^{-1}X=t_a$.
  \end{itemize}
\end{remarque}

  \begin{lemme} Quitte à remplacer $K$ par une extension finie et $R$ par sa clôture intégrale dans cette extension, on peut supposer qu'il existe des trivialisations admissibles de $(\ec,\theta,t)$.    
  \end{lemme}

  \begin{preuve} Quitte à remplacer $K$ par une extension finie, on peut et on va supposer que $j^*_R\ec$ est trivial sur $C_R^V$ (comme $G$ est semi-simple c'est possible par le théorème 3 de \cite{DS}). Le torseur $j^*_R\ec$ possède donc des sections  au-dessus de $C^V_R$. On en déduit un isomorphisme $\beta_R$ comme en (\ref{eq:betaR}). Quitte à translater cette section par un élément de $G(\kappa[C^V])$, on peut supposer que $\beta_R$ satisfait la condition 1 ci-dessus.\\

Soit $v\in V$. En utilisant la lissité de $\ec$ sur $C_R$, on voit que le torseur $i_{R,v}^*\ec$ possède des sections sur $\Spec(R[[z_v]])$. On en déduit des isomorphismes $\al_{R,v}$ comme en (\ref{eq:alphaRv}). Comme précédemment, on se ramène au cas où  $\al_{R,v}$ satisfait la condition 2. \\

Par la proposition \ref{prop:BL}, on déduit des trivialisations $\beta_R$ et $\al_{R,v}$ un triplet $(X,(g_v)_{v\in V},t)$  comme ci-dessus. Soit
$$Y=\Ad(g_\infty^{-1})X.$$
C'est un élément de $\ggo[[z_\infty]](R)$ qui vérifie $\chi_G(Y)=\chi_G(X)=\chi_G(t_a)$. Le lemme \ref{lem:section} montre que $Y$ et $t_a$ sont conjugués sous $G[[z_\infty]](R)$. Cela permet de conclure. 

\end{preuve}

\end{paragr}

\begin{paragr}[Sous-groupes unipotents.] --- Soit $\Phi^+=\Phi^{N_Q}_T$ et $\Phi^-=\Phi^{N_{Q^-}}_T$ les ensembles respectifs de  racines de $T$ dans les radicaux unipotents $N_Q$ et $N_{Q^-}$. La composante neutre $A_Q$ du  centre de $Q$ est un tore de dimension $1$. Soit
  \begin{equation}
    \label{eq:lambda}
    \la : \Gm \to A_Q
  \end{equation}
  l'unique isomorphisme qui vérifie $\al(\la)>0$ pour tout $\al\in \Phi^+$. Pour tout entier $i>0$, soit
$$\Phi^+_i=\{\al\in \Phi^+ \ | \ \al(\la)=i\}$$
et 
$$\Phi^-_i=-\Phi^+_i$$
en notation additive. 

Pour tout $\al\in \Phi^+\cup\Phi^-$, soit $N_\al$ le sous-groupe unipotent associé et $\zeta_\al \ : \ \mathbb{G}_a \to N_\al$ le groupe à un paramètre additif associé à $\al$. Pour tout entier $i>0$, soit $N_i^-$ le sous-groupe unipotent engendré par les sous-groupes $N_\al$ tels que $\al(\la)\leq -i$. On obtient ainsi  une suite décroissante de sous-groupes unipotents
$$ \ldots \subset N_i^- \subset \ldots \subset N_1^-=N_{Q^-}$$
et pour $i$ assez grand on a 
$$N_i^-=\{1\}.$$
On vérifie la relation suivante sur les commutateurs
$$[N_i,N_j]\subset N_{i+j}.$$
\end{paragr}

\begin{paragr}[Condition auxiliaire sur $V$.] ---  Commençons par le lemme suivant.

  \begin{lemme}\label{lem:Adn}
    Il existe $V_0$ un ensemble fini non vide de points fermés de $C_\kappa$  tel que  pour toute $\Gamma(C_\kappa-V_0,\oc_{C_\kappa})$-algèbre $B$ et tout entier $i\geq 1$, l'application de $N_{i}^-(B)$ dans $\ngo_{i}^-(B)$ définie par
$$n\in  N_{i}^-(B)\mapsto \Ad(n^{-1})\bar{X}-\bar{X}$$
est bijective.
\end{lemme}

\begin{preuve}
  Soit $F$ le corps des fonctions de $C\times_k \kappa$. On laisse au lecteur le soin de prouver que le morphisme $n \mapsto \Ad(n^{-1})\bar{X}-\bar{X}$ induit un isomorphisme de $ N_{i}^-\times_k F$ sur $\ngo_{i}^-\times_k F$ . Il se prolonge donc en un isomorphisme à un ouvert de Zariski de $C_\kappa$ complémentaire d'un ensemble fini $V_0$. Le lemme s'en déduit.
\end{preuve}

Quitte à agrandir $V$, on suppose dans la suite que $V$ contient l'ensemble fini $V_0$ donné par le lemme \ref{lem:Adn}.
\end{paragr}

\begin{paragr}[Constructions de certaines trivialisations admissibles.]  --- \label{S:cstr-de-triv}Soit $r$ le plus petit multiple commun des entiers $i$ tels que $\Phi^+_i\not=0$. Soit $R'=R[\pi^{\frac{1}{r}}]$ et $K'$ le corps des fractions de $R'$. Pour tout entier $j$ soit $z_j\in A_Q(K')$ défini par 
$$z_j=\la(\pi^{\frac{j}{r}})$$
où $\la$ est l'isomorphisme $\Gm \to A_Q$ défini en (\ref{eq:lambda}).

Soit $(\beta_R,(\al_{R,v})_{v\in V})$ une trivialisation admissible de $m=(\ec,\theta,t)$ (cf. définition \ref{def:trivadm}). Soit $(X,(g_v)_{v\in V},t)$ le triplet qui s'en déduit par la proposition \ref{prop:BL}. 

\begin{definition}
  \label{def:j}
Soit
$$j(\beta_R,(\al_{R,v}))\in \NN\cup \{\infty\}$$
la borne supérieure de l'ensemble des entiers $j\in \NN$ qui vérifient 
\begin{enumerate}
\item $\Ad(z_j)X \in \ggo(R'[C^V])$
\item pour tout $v\in V$, 
$$z_j g_v z_j^{-1} \in G((z_v))(R').$$
\end{enumerate}
\end{definition}

\begin{proposition}\label{prop:cruciale}
  Pour tout entier $j\in \NN$, il existe   une trivialisation admissible $(\beta_R,(\al_{R,v})_{v\in V})$ de $m=(\ec,\theta,t)$ telle que
$$j(\beta_R,(\al_{R,v}))\geq j.$$
\end{proposition}

\begin{preuve} Supposons la proposition mise en défaut. Soit $j\in \NN$ le plus petit entier tel que 
$$j\geq j(\beta_R,(\al_{R,v}))$$
pour toute trivialisation admissible $(\beta_R,(\al_{R,v})_{v\in V})$ de $m$. Comme on a $\bar{X}\in \lgo(\kappa[C^V])$ et $\bar{g}_v\in Q((z_v))(\kappa)$ (cf. assertion 1 du lemme \ref{lem:trivouvert} et assertion 1 du lemme \ref{lem:triv_en_v}), on a 
$$j\geq 1.$$

Soit $(\beta_R,(\al_{R,v})_{v\in V})$   une trivialisation admissible de $m$ telle que 
$$j=j(\beta_R,(\al_{R,v})).$$

Soit $(X,(g_v)_{v\in V},t)$ le triplet qui s'en déduit par la proposition \ref{prop:BL}.

Soit $\val_{\pi}$ la valuation de $R$ normalisée par $\val_\pi(\pi)=1$. Cette valuation s'étend à $R'$ et à $R'[C^V]$ de manière évidente. Pour tout $\al\in \Phi^{\pm}$, le $k$-espace vectoriel $\ngo_\al$ est de dimension $1$. Pour tout vecteur $X_\al\in \ngo_\al(R'[C^V])$,  soit $\val_\pi(X_\al)$ la valuation du coefficient de $X_\al$ sur un vecteur non nul de $\ngo_{\al}(k)$.

Suivant la décomposition en espaces propres sous $A_Q$
$$(\oplus_{\al\in \Phi^{\pm}}\ngo_\al) \oplus \mgo_Q $$
on écrit
\begin{equation}
  \label{eq:decompoX}
  X=\sum_{\al\in \Phi^\pm}X_\al  + X_{M_Q}.
\end{equation}
Comme la réduction modulo $\pi$ de $X$ est $\bar{X}$ qui appartient à $\lgo(\kappa[C^V])$, on a 

\begin{equation}
  \label{eq:valX}
  \val(X_\al)\geq 1
\end{equation}
pour tout $\al\in \Phi^{\pm}$. De plus, la réduction modulo $\pi$ de $X_{M_Q}$ est $\bar{X}$.

On a aussi la décomposition 
$$\Ad(z_j)X=\sum_{\al\in \Phi^\pm}\pi^{\frac{j\al(\la)}{r}} X_\al  + X_{M_Q}.$$
Par la condition 1 de la définition \ref{def:j}, on a pour tout $\al\in \Phi^{\pm}$
\begin{equation}
  \label{eq:inegvaluation}
  \frac{j\al(\la)}{r}+\val_{\pi}(X_\al)\geq 0.
\end{equation}
Pour $\al\in \Phi^+$, cette inégalité est stricte. 

\begin{lemme} \label{lem:existstricte}Il existe une trivialisation admissible $(\beta_R,(\al_{R,v})_{v\in V})$ de $m$ pour laquelle
  \begin{itemize}
  \item $j=j((\beta_R,(\al_{R,v})_{v\in V}))$ ;
  \item l'inégalité  (\ref{eq:inegvaluation}) est stricte pour tout $\al\in \Phi^{-}$.
  \end{itemize}
\end{lemme}

\begin{preuve}
Partons d'une trivialisation admissible  $(\beta_R,(\al_{R,v})_{v\in V})$ de $m$ pour laquelle il existe $\al\in \Phi^-$ tel que
\begin{equation}
  \label{eq:egal=valXal}
  \frac{j\al(\la)}{r}+\val_{\pi}(X_\al)=0.
\end{equation}
Soit $i$ le plus petit entier tel qu'il existe $\al\in \Phi^-_i$ qui vérifie l'égalité (\ref{eq:egal=valXal}). Soit 
$$Y=\Ad(z_j)X\in \ggo(R'[C^V])$$
et $ \Yb$ la réduction modulo $\pi^{\frac{1}{r}}$ de $Y$.  Alors $\Yb$ appartient à $\bar{X}+\ngo^-_i(\kappa[C^V])$. D'après le lemme \ref{lem:Adn}, il existe $n_i\in N_i^-(\kappa[C^V])$ tel que 
\begin{equation}
  \label{eq:XbYb}
  \Ad(n_i^{-1}) \bar{X} = \Yb .
\end{equation}
Pour tout $\al\in \Phi_i^-$ soit $u_\al\in \kappa[C^V]$ et $n_{i+1}\in N_{i+1}^-(\kappa[C^V])$ tels qu'on ait 
\begin{equation}
  \label{eq:n}
  n_i=n_{i+1}n
\end{equation}
avec
$$n=\prod_{\al\in \Phi_i^-}\zeta_\al(u_\al).$$
Comme $n_{i+1}$ est un élément de $N_{i+1}^-(\kappa[C^V])$, on a 
$$\Ad(n_{i+1}^{-1})\bar{X}\in \bar{X}+\ngo_{i+1}^-(\kappa[C^V]).$$
En comparant avec (\ref{eq:XbYb}) et (\ref{eq:n}), on obtient
\begin{equation}
  \label{eq:AdnYb}
  \Ad(n)\Yb\in \bar{X}+\ngo_{i+1}^-(\kappa[C^V]).
\end{equation}
Remarquons que 
$$z_j^{-1}nz_j=\prod_{\al\in \Phi_i^-}\zeta_\al( \pi^{ \frac{ji}{r}} u_\al)$$
appartient à $N_i(R[C^V])$ et que sa réduction modulo $\pi$ est triviale. En effet, par hypothèse, il existe $\al\in \Phi_i^-$ tel que l'égalité (\ref{eq:egal=valXal}) soit vraie. On a donc par  (\ref{eq:valX})
$$\frac{ji}{r}=-\frac{j\al(\la)}{r}=\val_\pi(X_\al)\in \NN^*.$$

Le triplet  $(\Ad(z_j^{-1}nz_j)X,(z_j^{-1}nz_jg_v)_{v\in V},t)$ définit par la bijection de la proposition \ref{prop:BL} une trivialisation de $m$ de la forme $(\beta_R',(\al'_{R,v})_{v\in V})$ qui est encore admissible. Il est clair sur la définition \ref{def:j} qu'on a  $j(\beta_R',(\al'_{R,v})_{v\in V})\geq j$ (par maximalité de $j$ on a même égalité). Quitte à remplacer $(\beta_R,(\al_{R,v})_{v\in V})$ par $(\beta_R',(\al'_{R,v})_{v\in V})$, et vu (\ref{eq:AdnYb}) on se ramène au cas où $\Yb$ appartient à $\bar{X}+\ngo^-_{i+1}(\kappa[C^V])$. Par récurrence, on peut même supposer qu'on a $\Yb=\bar{X}$. Mais alors l'inégalité \ref{eq:inegvaluation} est stricte.
\end{preuve}

Soit $(\beta_R,(\al_{R,v})_{v\in V})$ une trivialisation admissible de $m$ qui satisfait les conclusions du lemme \ref{lem:existstricte}. Soit $(X,(g_v)_{v\in V},t)$ le triplet qui s'en déduit par la proposition \ref{prop:BL}. On a donc
 \begin{equation}
  \label{eq:firstineg}
   \min_{\al\in \Phi^-}\big(-\frac{r}{\al(\la)}\val_{\pi}(X_\al)\big)>j.
\end{equation}
On notera que, par définition de $r$, le premier membre de cet inégalité est un entier.

Soit $v\in V$ et $B_v$ l'anneau déduit de $R((z_v))$ par localisation  en l'idéal premier engendré par $\pi$.  C'est un anneau de valuation discrète de corps résiduel $\kappa((z_v))$. On note encore $\val_\pi$ la valuation normalisée par $\val_\pi(\pi)=1$. Le point $g_v\in G((z_v))(R)$ induit alors un morphisme du trait $\Spec(B_v)$ dans $G$. Comme la réduction modulo $\pi$ de $g_v$ est $\gb_v$ et que l'on a $\gb_v\in Q((z_v))(\kappa)$ (cf. l.(\ref{eq:gbinQ})), l'image du point spécial par ce morphisme appartient à l'ouvert $N_{Q} M_Q N_{Q^-}$. Il s'ensuit que l'image du trait entier est dans l'ouvert $N_{Q} M_Q N_{Q^-}$. Il existe donc  $x_v \in M_Q(B_v)$ et  $b_{\al,v}\in B_v$ tels que pour tout $\al\in \Phi^\pm$ on ait
\begin{equation}
  \label{eq:decompogv}
  g_v=\big(\prod_{\al\in \Phi^+} \zeta_\al(b_{\al,v})\big) \cdot x_v \cdot \big(\prod_{\al\in \Phi^-} \zeta_\al(b_{\al,v})\big).
\end{equation}
Soit $\xb_v\in M_Q((z_v))(\kappa)$ la réduction modulo $\pi$ de $x_v$. Puisque la réduction $\gb_v$ de $g_v$ modulo $\pi$ appartient à  $Q((z_v))(\kappa)$, on a $\val_\pi(b_{\al,v})\geq 1$ pour $\al\in \Phi^-$ et 
\begin{equation}
  \label{eq:inMN_Q}
 \gb_v\in \xb_v N_Q((z_v))(\kappa).
\end{equation}
On a
\begin{equation}
  \label{eq:z_jg}
  z_j g_v z_j^{-1} =\big(\prod_{\al\in \Phi^+} \zeta_\al(\pi^{\frac{j\al(\la)}{r}}b_{\al,v})\big) \cdot x_v \cdot \big(\prod_{\al\in \Phi^-} \zeta_\al(\pi^{\frac{j\al(\la)}{r}}b_{\al,v})\big)
\end{equation}
et cet élément appartient à $G((z_v))(R')$ par définition de $j$.
 On a donc pour tout $\al\in \Phi^{\pm}$
\begin{equation}
  \label{eq:inegvaluation2}
  \frac{j\al(\la)}{r}+\val_{\pi}(b_{\al,v})\geq 0.
\end{equation}
L'inégalité (\ref{eq:inegvaluation2}) est stricte pour $\al\in \Phi^+$.

\begin{lemme} \label{lem:existstricte2}Il existe une trivialisation admissible $(\beta_R,(\al_{R,v})_{v\in V})$ de $m$ pour laquelle
  \begin{enumerate}
  \item $j=j((\beta_R,(\al_{R,v})_{v\in V}))$ ;
  \item l'inégalité (\ref{eq:firstineg}) est satisfaite ;
  \item l'inégalité  (\ref{eq:inegvaluation2}) est stricte pour tout $\al\in \Phi^{-}$.
  \end{enumerate}
\end{lemme}

\begin{preuve} Soit $(\beta_R,(\al_{R,v})_{v\in V})$ une trivialisation admissible de $m$ qui satisfait les conclusions du lemme \ref{lem:existstricte} donc les deux premières assertions. Soit $(X,(g_v)_{v\in V},t)$ le triplet qui s'en déduit par la proposition \ref{prop:BL}. On reprend les notations utilisées ci-dessus, en particulier aux lignes (\ref{eq:decompogv}) et (\ref{eq:inMN_Q}).
Soit
$$h_v =z_j g_v z_j^{-1}\in G((z_v))(R[\pi^{\frac{1}{r}}])$$
et 
$$Y=\Ad(z_j)X \in \ggo(R[C^V][\pi^{\frac{1}{r}}]).$$

\begin{lemme} \label{lem:m'}  On a les assertions suivantes

  \begin{enumerate}
  \item le triplet $(Y, (h_v)_{v\in V},t)$ est un triplet qui satisfait les conditions de la proposition \ref{prop:BL} pour l'anneau  $R'$ ;
  \item soit $m'\in \mc(R')$  le point associé à ce triplet par la proposition \ref{prop:BL}. Alors, après changement de base à $K'$, les points $m'$ et $m$ deviennent isomorphes.
  \end{enumerate}  
\end{lemme}

\begin{preuve} Pour l'assertion 1, il s'agit de vérifier que le triplet $(Y, (h_v)_{v\in V},t)$ satisfait les conditions suivantes :
  \begin{itemize}
   \item $Y \in \ggo(R'[C^V])$ ;
   \item pour tout $v\in V$, l'élément $h_v \in G((z_v))(R')$ et vérifie
$$\Ad(h_v^{-1})Y\in z_v^{-d_v}\ggo[[z_v]](R') \ ;$$
\item $t$ est un élément de $\tgo^{\reg}(R')$ dont la caractéristique $\chi(t)$ est égale à la réduction modulo $z_\infty$ de $\chi_G(Y)$. 
   \end{itemize}
La première condition et la relation  $h_v \in G((z_v))(R')$ résultent de la définition de $j$. On a donc 
$$\Ad(h_v^{-1})Y \in  \ggo((z_v))(R').$$
En fait, on peut remplacer dans la relation précédente  $\ggo((z_v))$ par $ \ggo[[z_v]]$ car on a d'une part  
$$\Ad(h_v^{-1})Y=\Ad(z_j) (\Ad(g_v^{-1})X)$$
et d'autre part
$$\Ad(g_v^{-1})X\in z_v^{-d_v}\ggo[[z_v]](R).$$

De même, la troisième relation est vérifiée puisque $X$ et $Y$, étant conjugués, ont même caractéristique que cette relation vaut $X$.\\

Pour l'assertion 2 du lemme, on remarque qu'on a $z_j\in G(K')$ et on conclut à l'aide du corollaire \ref{cor:BL}.
\end{preuve}

\begin{lemme} \label{lem:P-} On poursuit avec les notations du lemme \ref{lem:m'}. 
  Soit $m'_\kappa\in \mc(\kappa)$ le point déduit de $m'$ par le changement de base $\Spec(\kappa)\to\Spec(R)$. On a alors les assertions suivantes :
  \begin{enumerate}
  \item $f(m')\in \Ac_{L,\el}(\kappa)$ ;
  \item pour tout $P\in \fc(L)$ et toute réduction $m'_{\kappa,P}$ de $m'_\kappa$ à $P$, on a $$\deg(m_{\kappa,P}')=H_P((\hb_v)_{v\in V})$$
où $\hb_v\in G((z_v))(\kappa)$ est la réduction modulo $\pi^{\frac{1}{r}}$ de $h_v$ ;
\item pour tout $P\in \fc(L)$ tel que $P\subset Q$ soit $ m_{\kappa,P}$ une réduction de $m_\kappa$ à $P$ et $m_{\kappa,P^-}'$ une réduction de $m_\kappa'$ à $P^-=(M_Q\cap P)N_{Q^-}$ ;  on a l'égalité 
$$\deg(m_{\kappa,P^-}')=\deg(m_{\kappa,P})$$
  \end{enumerate}

\end{lemme}

\begin{preuve}  Le triplet $(\Yb, (\hb_v)_{v\in V},\bar{t})$ est la réduction modulo $\pi^{\frac{1}{r}}$ du triplet $(Y, (h_v)_{v\in V},t)$. Par conséquent, le point $m'_\kappa$ se déduit du triplet $(\Yb, (\hb_v)_{v\in V},\bar{t})$ par la bijection de la proposition \ref{prop:BL}. L'inégalité (\ref{eq:firstineg}) implique que $\Yb$ est égal à $\bar{X}$ défini en (\ref{eq:tripletadele}). Il s'ensuit que $f(m'_\kappa)=f(m_\kappa)$ et ce point appartient à 
$\Ac_{L,\el}(\kappa)$ (cf. \S \ref{S:YHN}) d'où l'assertion 1. Rappelons que $\bar{X}\in \lgo(\kappa[C^V])$ est conjugué à $\bar{t}_a$ par $\gb_\infty\in L((z_\infty))(\kappa)$ (cf. assertion 2 et 3 du lemme \ref{lem:triv_en_v}). L'assertion 2 résulte alors du corollaire \ref{cor:reducadel}.

Prouvons l'assertion 3. Soit $P\in \fc(L)$ tel que $P\subset Q$. Soit $P^-\in \fc(L)$ défini par $P^-=(M_Q\cap P)N_{Q^-}$. On a donc $P^-\subset Q^-$. De l'égalité (\ref{eq:z_jg}) et du fait que l'inégalité (\ref{eq:inegvaluation2}) est stricte pour $\al\in \Phi^+$ on déduit
$$\hb_v\in N_{Q-}((z_v))(\kappa)\xb_v$$
et
$$\deg(m'_{\kappa,P^-})=H_{P^-}((\hb_v)_{v\in V})=H_{P^-\cap M_Q}((\xb_v)_{v\in V}).$$

Mais l'égalité (\ref{eq:inMN_Q}) implique aussi qu'on a 
$$\deg(m_{\kappa,P})=H_{P}((\gb_v)_{v\in V})=H_{P\cap M_Q}((\xb_v)_{v\in V}).$$
Comme $P\cap M_Q=P^{-}\cap M_Q$, on a bien 
$$\deg(m'_{\kappa,P^-})=\deg(m_{\kappa,P})$$
\end{preuve}

D'après les notations du paragraphe \ref{S:YHN} et les équivalences de la proposition \ref{prop:HN}, $\varrho$ est le $\xi$-point de Harder-Narasimhan de $m_\kappa$ et $Q_0$ est le  sous-groupe parabolique de $\fc(L)$ qui vérifie 
\begin{itemize}
\item[(A)] $\xi-\varrho \in \ago_{Q_0}^+\cap \ago_T^G$ ;
\item[(B)] La projection de $\varrho$ sur $\ago_L^G$ appartient à la projection sur $\ago_L^G$ de l'enveloppe convexe des points $-\deg(m_{\kappa,P})$ pour $P\in \pc(L)$ tel que  $P\subset Q_0$ et $m_{\kappa,P}$ une réduction de $m_\kappa$ à $P$. 
\end{itemize}

 L'assertion (B) ci-dessus et l'assertion 3 du lemme \ref{lem:P-} montrent que la projection de $\varrho$ sur $\ago_L^G$ appartient à la projection sur $\ago_L^G$ de l'enveloppe convexe des points $-\deg(m'_{\kappa,P^-})$ pour $P\in \fc^{Q_0}(L)$ (avec les notations du lemme \ref{lem:P-}). La proposition \ref{prop:envpecvx} implique  alors que le point $ \varrho$ appartient au convexe $\overline{\mathcal{C}}_{m'_\kappa}$ associé à $m'_\kappa$. 

Il s'ensuit qu'on a l'inégalité suivante sur les distance
$$d(\overline{\mathcal{C}}_{m'_\kappa},\xi)\leq d(\varrho_{m_\kappa},\xi)=d(\overline{\mathcal{C}}_{m_\kappa},\xi)=d.$$
Mais par définition de $d$ (cf. \S\ref{S:d}) on a nécessairement égalité dans l'inégalité ci-dessus. Donc  $\varrho$ est aussi le $\xi$-point de Harder-Narasimhan de $m'_\kappa$. Comme $\xi-\varrho \in \ago_{Q_0}^+\cap \ago_T^G$, le sous-groupe parabolique qui vérifie l'assertion 2 de la  proposition \ref{prop:HN} pour $m'_\kappa$ est nécessairement $Q_0$. La projection de $\varrho$ sur $\ago_L^G$ est donc égale à la projection sur $\ago_L^G$ de l'enveloppe convexe des points $-\deg(m'_{\kappa,P})$ associés aux  réductions $m'_{\kappa,P}$ de $m'_\kappa$ aux sous-groupes $P\in \pc(M)$ tels que $P\subset Q_0$.
En particulier, la projection de $\varrho$ sur la droite $\ago_{Q}^G$ égale au point $-\deg(m'_{\kappa,Q})$ où $m'_{\kappa,Q}$ est une réduction de $m'_\kappa$ à $Q$ (on a $\ago_G=\{0\}$ puisque $G$ est semi-simple). L'assertion (B) implique par ailleurs que cette projection est égale à  $-\deg(m_{\kappa,Q})$. On a donc 
$$\deg(m'_{\kappa,Q})=\deg(m_{\kappa,Q})$$
ce qui, combiné avec l'assertion 3 du lemme \ref{lem:P-} donne
$$\deg(m'_{\kappa,Q})=\deg(m'_{\kappa,Q^-}).$$
ou de manière équivalente par l'assertion 2 du lemme \ref{lem:P-}

\begin{equation}
  \label{eq:HQ}
  H_{Q^-}((\hb_v)_{v\in V})=H_{Q}((\hb_v)_{v\in V}).
\end{equation}
Le lemme \ref{lem:reducM} ci-dessous implique alors que pour tout $v\in V$, on a
\begin{equation}
  \label{eq:reduchbv}
  \hb_v \in M_Q((z_v))(\kappa) G[[z_v]](\kappa)
\end{equation}

Considérons d'abord le point $v=\infty_\kappa$. Par la condition 3 qui définit une trivialisation admissible, on a $\Ad(g_\infty^{-1})X=t_a$, d'où l'on tire
 $$\Ad(h_\infty^{-1})Y=t_a.$$
L'inégalité (\ref{eq:firstineg}) implique qu'on a $\Yb=\bar{X}$. En réduisant modulo $\pi^{\frac{1}{r}}$ l'inégalité ci-dessus, on obtient
 $$\Ad(\hb_\infty^{-1})\bar{X}=\bar{t}_a.$$
Mais  $\bar{X}$ et $ \bar{t}_a$ sont conjugués par $\gb_\infty\in L((z_\infty)(\kappa)$ (cf. assertion 3 du lemme \ref{lem:triv_en_v}. On a donc aussi $\hb_\infty\in   L((z_\infty))(\kappa)$. Par conséquent l'inégalité (\ref{eq:inegvaluation2}) est stricte pour $v=\infty_\kappa$ et pour tout $\al\in \Phi^\pm$.

Supposons désormais $v\not=\infty_\kappa$. Supposons aussi qu'il existe $\al\in \Phi^-$ tel que l'inégalité  (\ref{eq:inegvaluation2}) soit une égalité c'est-à-dire qu'on ait
\begin{equation}
  \label{eq:egvaluation2}
  \frac{j\al(\la)}{r}+\val_{\pi}(b_{\al,v})= 0.
\end{equation}
Soit $i$ le plus entier tel qu'il existe $\al\in \Phi^-_i$ qui vérifie l'égalité (\ref{eq:egvaluation2}).  On a donc
\begin{equation}
  \label{eq:inmN_i}
  \hb_v\in \xb_v N_i((z_v))(\kappa).
\end{equation}
En comparant avec (\ref{eq:reduchbv}), on en déduit que
\begin{equation}
  \label{eq:inmN_i2}
  \hb_v\in \xb_v N_i[[z_v]](\kappa).
\end{equation}
Pour tout $\al\in \Phi_i^-$ soit $u_\al\in \kappa[[z_v]]$ et
$$n= \prod_{\al\in \Phi_i^-}\zeta_\al(u_\al).$$
tel qu'on ait  
$$ \hb_v n \in \xb_v N_{i+1}[[z_v]](\kappa).$$
Soit $n'=z_j^{-1}nz_j$. On a donc 
$$n'=\zeta_\al( \pi^{\frac{ij}{r}} u_\al).$$
L'égalité (\ref{eq:egvaluation2}) implique que $\frac{ij}{r}$ est un entier strictement positif. Par conséquent, $n'$ appartient à $N_Q((z_v))(R)$ et sa réduction modulo $\pi$ est triviale. Quitte à remplacer $g_v$ par $g_v n'$, on voit qu'on est ramené au cas où l'on a 
$$ \hb_v\in \xb_v N_{i+1}[[z_v]](\kappa)$$
et par récurrence au cas où $\hb_v=\xb_v$. Mais, dans ce cas, l'inégalité (\ref{eq:inegvaluation2}) est stricte pour tout $\la\in \Phi^{\pm}$. Cela termine la démonstration du lemme \ref{lem:existstricte2}.
\end{preuve}

Soit $(\beta_R,(\al_{R,v})_{v\in V})$ une trivialisation admissible de $m$ qui vérifie les conditions du lemme \ref{lem:existstricte2}. Soit $(X,(g_v)_{v\in V},t)$ le triplet qui s'en déduit. Soit
$$j_1=\min\big(\min_{\al\in \Phi^-}\big(-\frac{r}{\al(\la)}\val_{\pi}(X_\al)\big),\min_{v\in V, \al\in \Phi^-}\big(-\frac{r}{\al(\la)}\val_{\pi}(b_{\al,v})\big)\big).$$
C'est un entier qui vérifie $j_1>j$,
$$\Ad(z_{j_1})X\in \ggo(R[C^V])$$
et pour tout $v\in V$
$$z_{j_1}g_v z_{j_1}^{-1}\in G((z_v))(R).$$
C'est visiblement la contradiction recherchée et cela achève la démonstration de la proposition \ref{prop:cruciale}. 
 
\end{preuve}

Le lemme suivant a été utilisé dans la preuve  précédente (preuve de la proposition \ref{prop:cruciale}).

  \begin{lemme}\label{lem:reducM}
   Soit $Q\in \fc^G$ un sous-groupe parabolique de $G$ et $Q^-$ son parabolique opposé. Pour tout $v\in V$, soit $g_v\in G((z_v))(\kappa)$. L'égalité 
   \begin{equation}
     \label{eq:egaliteH_Q}
     H_Q((g_v)_{v\in V})=H_{Q^-}((g_v)_{v\in V})
   \end{equation}
   entraîne que pour tout $v\in V$, il existe $x_v\in M_Q((z_v))(\kappa)$ et $k_v\in G[[z_v]](\kappa)$ tels que
$$g_v=x_v\,k_v$$
  \end{lemme}

  \begin{preuve} Traitons d'abord le cas où $V=\{v\}$ est un singleton. Pour alléger les notations, on omet l'indice $v$. Soit $g\in G((z))(\kappa)$. On va montrer que $-H_Q(g)+H_{Q^-}(g)$ est une somme à coefficients \emph{positifs} de coracines $\al^\vee$ dans $\Delta_Q^\vee$ et que cette somme est nulle si et seulement si  il existe $m\in M_Q((z))(\kappa)$ et $k\in G[[z]](\kappa)$ tels que $g=m\,k$. En utilisant la décomposition d'Iwasawa
$x=m \, n \, k$
avec $m\in M_Q((z))(\kappa)$, $n\in N_Q((z))(\kappa)$ et $k\in G[[z]](\kappa)$, on voit qu'on se ramène à montrer l'assertion suivante : pour tout $x\in  N_Q((z))(\kappa)$ on a
$$H_{Q^-}(x)\geq0$$
et $H_{Q^-}(x)=0$ implique que $x$ appartient à $N_Q[[z]](\kappa)$. Pour toute racine $\al$ de $T$ dans $G$ soit $\zeta_\al$ un isomorphisme de $\mathbb{G}_a$ sur le groupe radiciel correspondant à la racine $\al$. Soit $x\in N_Q((z))(\kappa)$. Il existe une famille $(\al_i)_{1\leq i\leq n}$ de racines de $T$ dans $N_Q$ deux à deux distintes et des éléments $x_i\in \kappa((z))$ tels que
$$x=\prod_{i=1}^n \zeta_{\al_i}(x_i)= \zeta_{\al_1}(x_1)\ldots \zeta_{\al_n}(x_n).$$
Si $x$ appartient à $N_Q[[z]](\kappa)$ on a   $H_{Q^-}(x)=0$. Supposons que $x$ n'appartient pas à $N_Q[[z]](\kappa)$. Il existe donc $i$ tel que $\val(x_i)<0$. Quitte à translater $x$ à droite par un élément de $N_Q[[z]](\kappa)$, ce qui ne modifie pas la valeur de $H_{Q^-}(x)$, on peut supposer  $\val(x_n)<0$ . Pour toute racine $\al$ et tout $a\in \kappa((z))$ de valuation strictement négative, on vérifie l'assertion suivante
$$\zeta_{\al}(a)\in \al^\vee(a) \zeta_{-\al}(a) G[[z]](\kappa).$$
En utilisant cette relation  pour $\al_n$ et $x_n$, on voit que pour tout $1\leq i\leq n-1$, il existe $x'_i \in \kappa((z))$ tel que 
$$H_{Q^-}(x)=H_{Q^-}(\al^\vee(x_n) \prod_{i=1}^{n-1}\zeta_{\al_i}(x'_i)).$$
En raisonnant par récurrence, on voit qu'il existe un entier $1\leq n'\leq n$ tel que pour tout $1\leq i\leq n'$  il existe 
\begin{itemize}
\item un élément $y_i\in \kappa((z))$ tel que  $\val(y_i)<0$ ;
\item une racine $\beta_i$ de $T$ dans $N_Q$
\end{itemize}
tels que 
$$H_{Q^-}(x)=H_{Q^-}\big(\prod_{i=1}^{n'} \beta^\vee_i(y_i)\big).$$
Soit $\gamma_i$ la projection de $\beta_i$ sur $\ago_Q$ : c'est un vecteur non nul. On obtient donc  
$$H_{Q^-}(x)=-\sum_{i=1}^{n'}\val(y_i) \gamma_i^\vee$$
ce qui montre que $H_{Q^-}(x)$ est une combinaison à coefficients positifs d'éléments de $\Delta_Q$ et même strictement positifs. Donc  $H_{Q^-}(x)$ est non nul. Cela conclut la preuve lorsque $V$ est un singleton.\\

Revenons au cas général d'un ensemble $V$ fini non vide.  On vient de voir  que pour tout $v\in V$, la différence $-H_Q(g_v)+H_{Q^-}(g_v)$ est  une combinaison à coefficients positifs d'éléments de $\Delta_Q$. Il s'ensuit que l'égalité (\ref{eq:egaliteH_Q}) est vraie si et seulement si pour tout $v\in V$
$$ H_Q(g_v)=H_{Q^-}(g_v).$$
On est donc ramené au cas du singleton.
  \end{preuve}

\end{paragr}

  \begin{paragr}[Où l'on obtient la contradiction cherchée.] --- \label{S:lacontradiction}Pour tout entier $n$, soit
$$R_n=R/ \pi^n R.$$
Soit $n$ un entier $\geq 1$. Soit $(\beta_R, \al_{R,v})$ une trivialisation admissible de $m$ qui vérifie
\begin{equation}
  \label{eq:j>n}
  j(\beta_R, \al_{R,v})\geq n,
\end{equation}
(cf. proposition \ref{prop:cruciale}. Soit $(X,(g_v)_{v\in V},t)$ le triplet associé à $m$ et à cette trivialisation par la proposition \ref{prop:BL}. Soit $(X_n,(g_{n,v})_{v\in V},t_n)$ le triplet qui s'en déduit par réduction modulo $\pi^n$. Ce triplet vérifie les conditions suivantes :
\begin{enumerate}
\item $X_n\in \qgo(R_n[C^V])$ ;
\item pour tout $v\in V$, $g_{n,v}\in Q((z_v))(R_n)$ ;
\item $\Ad(g_{n,\infty})^{-1}X_n= t_{a,n}$ où $t_{a,n}$ est la réduction modulo $\pi^n$ de $t_a$.
\end{enumerate}

De ce triplet, on déduit, par  la proposition \ref{prop:BL} appliquée au groupe $Q$, un triplet $m_{Q,n}=(\ec_{Q,n},\theta_{Q,n},t_n)$ formé d'un $Q$-torseur $\ec_{Q,n}$ sur
$$C_n=C\times_k R_n,$$
d'une section $\theta_{Q,n}$ de $\Ad_D(\ec_{Q,n})$ tel que $\chi_Q(\theta_{Q,n}(\infty_{R_n}))=\chi_Q(t_n)$. En outre, ce triplet est une réduction à $Q$ du point $m_n\in \mc(R_n)$ obtenu par changement de base à $R_n$. 

L'algorithme utilisée dans la preuve de la proposition \ref{prop:cruciale} consiste à conjuguer $X$ (et translater à droite $g_v$ en conséquence), resp. translater à gauche $g_v$, par des éléments de la forme $\zeta_\al(u_\al)$ où $u_\al$ est un élément de $R[C^V]$, resp. de $R((z_v))$, de valuation $\pi$-adique $\geq n$ de façon à obtenir un triplet $(X',(g'_v)_{v\in V},t)$ associée à une trivialiation admissible dont le $j$ associé est plus grand que $n+1$. En particulier la réduction modulo $\pi^n$ de  $(X',(g'_v)_{v\in V},t)$ est égale à $(X_n,(g_{n,v})_{v\in V},t_n)$. Soit $(X_{n+1},(g_{n+1,v})_{v\in V},t_{n+1})$ la réduction modulo $\pi^{n+1}$ de  $(X',(g'_v)_{v\in V},t)$. Soit $m_{Q,n+1}=(\ec_{Q,n+1},\theta_{Q,n+1},t_{n+1})$ la réduction à $Q$ de $m_{n+1}\in \mc(R_{n+1})$ qui s'en déduit comme plus haut. Le changement de base de  $m_{Q,n+1}$ à $C_n$ redonne $m_{Q,n}$.

Par récurrence, on obtient pour tout entier $n$ une réduction  $m_{Q,n}=(\ec_{Q,n},\theta_{Q,n},t_n)$ à $Q$ du point $m_n=(\ec_n,\theta_n,t_n) \in \mc(R_n)$ déduit de $m\in \mc(R)$ par changement de base telle que l'image de $m_{n+1}$ dans $\mc(R_n)$ soit isomorphe à $m_n$. Le morphisme $Q$-équivariant
$$\ec_{Q,n} \to \ec_{Q,n}\times^G_k Q \simeq \ec_n$$
fournit par passage au quotient une section
$$\sigma_n \ : \ C_n \to \ec_n/Q.$$
Lorsque $n$ varie, les sections obtenues s'inscrivent dans un diagramme commutatif
$$\xymatrix{\ec_{n}/Q \ar[r]  & \ec_{n+1}/Q \\ C_{n} \ar[r]  \ar[u]^{\sigma_{n}} & C_{n+1} \ar[u]_{\sigma_{n+1}} } $$
où les morphismes horizontaux sont les flèches canoniques de transition. Elles forment donc un système inductif $\hat{\sigma}=(\sigma_n)_{n\geq 1}$. Soit 
 $$\hat{C}=\varinjlim_{n} C_n$$ et 
$$\widehat{\ec/Q}=\varinjlim_{n}\ec_n/Q.$$
Le système inductif $\hat{\sigma}$ donne un élément, encore noté  $\hat{\sigma}$, de 
$$\Hom_{\Spf(R)}(\hat{C},\widehat{\ec/Q}).$$
Comme $C_R$ est propre sur $\Spec(R)$ et que $\ec/Q$ est séparé et de type fini sur $\Spec(R)$, on sait, d'après le théorème d'algébrisation des morphismes de Grothendieck (théorème 5.4.1 de \cite{EGAIII1}), que l'application de prolongement aux complétés
$$\Hom_{\Spec(R)}(C_R,\ec/Q)\to \Hom_{\Spf(R)}(\hat{C},\widehat{\ec/Q})$$
est une bijection. Il existe donc $\sigma\in \Hom_{\Spec(R)}(C_R,\ec/Q)$ qui s'envoie sur $\hat{\sigma}$ par la bijection ci-dessus. En utilisant de nouveau le théorème d'algébrisation, on voit que $\sigma$ est en fait une section de $\ec/Q$. Soit $\ec_Q$ le $Q$-torseur au-dessus de $C_R$ défini par 
$$\ec_Q=C_R\times_{\sigma,\ec/Q}\ec.$$
Par construction, $\ec_Q$ est une réduction de $\ec$ à $Q$.

Soit 
$$\widehat{\Ad_D(\ec_Q)}=\varinjlim_{n}\Ad_D(\ec_{Q,n}).$$
Pour tout entier $n\geq 1$, la section $\theta_n$ de $\Ad_D(\ec_n)$ se factorise par la section $\theta_{Q,n}$ de $\Ad_D(\ec_{Q,n})$. Par construction, les sections $\theta_n$ forment un système inductif et donc un élément noté de 
$$\Hom_{\Spf(R)}(\hat{C},\widehat{\Ad_D(\ec_Q)}).$$
Par le théorème d'algébrisation déjà cité, on en déduit que $\theta$ se factorise par une section $\theta_Q$ de $\Ad_D(\ec_Q)$.

Par construction, le réduction modulo $\pi^n$ de $\chi_Q(\theta_{Q})(\infty_R)$ est égale à $\chi_Q(\theta_{Q,n})(\infty_{R_n})=\chi_Q(t_{n})$. Il s'ensuit qu'on a
 $$\chi_Q(\theta_{Q})(\infty_R)=\chi_Q(t).$$
Par conséquent le triplet $m_Q=(\ec_Q,\theta_Q,t)$ est une réduction de $m$ à $Q$. Soit $m_{Q,K}$ et $m_{Q,\kappa}$ les réductions à $Q$ de $m_K$ et $m_\kappa$ qui s'en déduisent par changement de base. Par platitude de $\ec_Q$ sur $C_R$, on a
$$\deg(m_{Q,K})=\deg(m_{Q,\kappa}).$$
Or $m_K$ est $\xi$-semi-stable : on a donc 
$$\xi\in -\deg(m_{Q,K})-\, \overline{^+ \ago_Q}+ \ago_T^Q.$$
Au paragraphe \ref{S:YHN}, on a introduit le  sous-groupe parabolique $Q_0\subset Q$ qui vérifie $\xi-\rho\in \ago_{Q_0}^+$ et $\varrho\in -\deg(m_{Q_0,\kappa})+\ago_T^{Q_0}$. On a donc 
$$\xi \in -\deg(m_{Q_0,\kappa})+\ago_{Q_0}^+ + \ago_T^{Q_0}.$$ 
Il s'ensuit qu'on a 
$$\xi\in \big(-\deg(m_{Q,K})-\, \overline{^+ \ago_Q}+ \ago_T^Q \big) \cap \big(-\deg(m_{Q_0,\kappa})+\ago_{Q_0}^+ + \ago_T^{Q_0}\big).$$
De l'égalité $\deg(m_{Q,K})=\deg(m_{Q,\kappa}),$ on tire $\big(-\, \overline{^+ \ago_Q}+ \ago_T^Q \big) \cap \big(\ago_{Q_0}^+ + \ago_T^{Q_0}\big)\not=\emptyset$. En projetant sur $\ago_Q$, on trouve $-\, \overline{^+ \ago_Q} \cap \ago_Q^+\not=\emptyset$ ce qui n'est pas : c'est la contradiction cherchée (cf. fin du \S\ref{S:d}).

\end{paragr}

\section{Séparation du morphisme $f^\xi$}\label{sec:separation}

\begin{paragr} Le but de cette section est de prouver le théorème suivant. 

  \begin{theoreme} \label{thm:separation} Supposons $G$ semi-simple. Pour tout $\xi\in \ago_T$, le morphisme 
$$f^\xi\ : \ \mc^\xi \to \Ac$$
est séparé.  
  \end{theoreme}

  \begin{remarque}
    L'énoncé et la démonstration du théorème \ref{thm:separation} s'inspirent d'un théorème de Langton (cf. \cite{Langton}, \S3 théorème). Faltings a donné une preuve du théorème \ref{thm:separation} pour $\xi=0$ et lorsque le corps de base est de caractéristique nulle (cf. \cite{Faltings} théorème II.4). 
  \end{remarque}

\medskip

Dans toute la suite, on suppose $G$ \emph{semi-simple}. D'après le critère valuatif de séparation, il est équivalent de prouver l'énoncé suivant.

\begin{proposition} Soit $\kappa$ une extension de $k$ algébriquement close. Soit $R$ un anneau de valuation discrète complet, de corps des fractions $K$ et de corps résiduel $\kappa$. Soit $m$ et $m'$ deux éléments de $\mc^\xi(R)$ et $m_{K}$ et $m_{K}'$ les éléments de $\mc^\xi(K)$ qui s'en déduisent par changement de base. Supposons que ces objets satisfont les deux conditions suivantes : 
  \begin{itemize}
  \item $f^\xi(m)=f^\xi(m')$ ;
  \item il existe un isomorphisme $\phi_K : m_{K}' \to m_{K}$.
  \end{itemize}
Alors il existe un unique isomorphisme $\phi : m_{}' \to m_{}$ qui prolonge $\phi_K$.
\end{proposition}
\end{paragr}

\begin{paragr}[Où l'on se ramène à un problème sur un  trait.] --- Dans toute la suite, soit deux triplets $m=(\ec,\theta,t)$ et $m'=(\ec',\theta',t)$ dans $\mc^\xi(R)$ tels que $f^\xi(m)=f^\xi(m')$. Soit  $m_K=(\ec_K,\theta_K,t)$ et $m'=(\ec'_K,\theta'_K,t)$ les triplets de $\mc^\xi(K)$ qui s'en déduisent par changement de base. Soit $\phi_K : m_{K}' \to m_{K}$ un isomorphisme. On note encore  $\phi_K$ l'isomorphisme $G$-équivariant sous-jacent
$$\phi_K \ : \ \ec_K'\to \ec_K.$$
L'isomorphisme qui s'en déduit
$$\Ad_D(\phi_K) \ :\ \Ad_D( \ec_K') \to \Ad(\ec_K)$$
envoie $\theta_K'$ sur $\theta_K$.  On cherche à prolonger $\phi_K$ en un isomorphisme $G$-équivariant
$$\phi \ : \ \ec'\to \ec$$
de sorte que l'isomorphisme qui s'en déduit
$$\Ad_D(\phi) \ :\ \Ad_D( \ec) \to \Ad(\ec)$$
envoie $\theta'$ sur $\theta$. Soit $\pi$ une uniformisante de $R$ et $B$ l'anneau local de $C_R$ en le point de codimension $1$ défini par l'idéal $(\pi)$. C'est donc un anneau de valuation discrète de corps résiduel le corps des fonctions $\kappa(C)$ de la courbe $C_\kappa$ et de corps des fractions le corps des fonctions $F$ de la courbe $C_K$.

\begin{lemme}\label{lem:reductrait}
  Pour que $\phi_K$ se prolonge à $C_R=C\times_k R$ il faut et il suffit qu'il se prolonge à $\Spec(B)$. De plus, de tels prolongements, s'ils existent, sont uniques.
\end{lemme}

\begin{preuve} Un isomorphisme $G$-équivariant de $\ec'$ sur $\ec$ n'est autre  qu'une section globale  du $C_R$-schéma
$$\mathrm{Isom}_G(\ec',\ec)=(\ec'\times_{C_R} \ec)/G$$
où $G$ agit diagonalement à droite sur $\ec'\times_{C_R} \ec $. L'isomorphisme $\phi_K$ définit  une section $\sigma_K$ de $\mathrm{Isom}_G(\ec',\ec)$ au-dessus de $C_K$. Tout prolongement de cette section à $C_R$ est nécessairement unique, d'où l'unicité de l'énoncé. Supposons que $\phi_K$ se prolonge à $\Spec(B)$. La section $\sigma_K$ se prolonge donc à $\Spec(B)$ et par suite à un ouvert de $C_R$ qui contient tous les points de codimension $\leq 1$. Comme $\mathrm{Isom}_G(\ec',\ec)$ est affine sur $C_R$  --- c'est même un torseur sous $\Aut_G(\ec)$ ---, une telle section se prolonge automatiquement à $C_R$ (cf. \cite{EGAIV4} corollaire 20.4.12).
\end{preuve}

Notons que si $\phi_K$ se prolonge en un isomorphisme $\phi$, l'isomorphisme $\Ad_D(\phi)$ envoie nécessaire\-ment $\theta'$ sur $\theta$ puisque les sections $\Ad_D(\phi)(\theta')$ et $\theta$, qui coïncident sur un ouvert de $C_R$, sont en fait égales.
\end{paragr}

\begin{paragr}[Étude sur le trait $\Spec(B)$] ---\label{S:trait}
Pour toute extension $K'$ finie de $K$ soit $R'$ la clôture intégrale de $R$ dans $K'$. Soit $F'$ le corps des fonctions de la courbe relative $C_{R'}=C\times_k R'$ et  $B'$  la clôture intégrale de $B$ dans $F'$. Dans la démonstration du lemme \ref{lem:reductrait}, on a vu que $\phi_K$ s'interprète comme une section du $C_R$-schéma affine $\mathrm{Isom}_G(\ec',\ec)$ et que prolonger $\phi_K$ revient à prolonger cette section à $\Spec(B)$. Comme $B$ est normal, $\phi_K$ se prolonge à $\Spec(B)$ si et seulement si le changement de base de $\phi_K$ à $K'$ se prolonge à $B'$. Dans la suite, on pourra toujours, s'il le faut, remplacer $K$ par $K'$.

Soit $(a,t)\in \Ac_G(R)$ défini par $(a,t)=f^\xi(m)$. Soit $t_a\in \tgo^{G\textrm{-}\reg}[[z_\infty]](R)$ l'unique élément dont la  caractéristique est la restriction de $a$ à $\Spec(R[[z_\infty]]$ et dont la réduction modulo $z_\infty$ est $t$. 

Soit $M\subset L$ les deux sous-groupes de Levi dans $\lc$  tels que  $(a,t)$ appartienne en fibre générique à $\chi^L_G(\Ac_{L,\el}(K))$ et à $\chi^M_G(\Ac_{M,\el}(\kappa))$ en fibre spéciale (cf. proposition \ref{prop:reunionAM}). Soit $(a_L,t)\in \Ac_L(R)$ l'unique point tel que $(a,t)=\chi^L_G((a_L,t))$ (cf. proposition \ref{prop:immersionfermee}).

\begin{lemme} \label{lem:specB} Quitte à remplacer $K$ par un extension finie $K'$ assez grande et $R$ par $R'$, on est dans la situation suivante 
\begin{enumerate}
\item Il existe un point  $X\in \lgo(B)$ qui vérifie
  \begin{enumerate}
  \item $\chi_L(X)$ coïncide avec $a_L$ sur $\Spec(B)$ ;
  \item la réduction $\bar{X}$ de $X$ modulo $\pi$ appartient à $\mgo(\kappa(C))$ et est conjugué à la réduction modulo $\pi$ de $t_a$ par un élément de $M((z_\infty)(\kappa)$ ;
  \end{enumerate}

\item Il existe des trivialisations des $G$-torseurs  $\ec$ et $\ec'$  sur $\Spec(B)$ tels que les trivialisations de $\Ad_D(\ec)$ et $\Ad_D(\ec')$ qui s'en déduisent envoient $\theta$ et $\theta'$ sur $X$.
\end{enumerate}
  \end{lemme}

  \begin{preuve}
    La section de Kostant montre qu'un point $X\in \lgo(B)$ qui vérifie 1.(a) existe. Par la proposition \ref{prop:caractdesAM} et le lemme \ref{lem:cjsep}, on peut supposer que, quitte à conjuguer $X$ par un élément de $L(\kappa(C))$, l'assertion 1.(b) est aussi vérifiée. 

 Si l'extension $K'$ est assez grande les torseurs $\ec$ et $\ec'$ sont triviaux sur $\Spec(B')$ (théorème 2 de \cite{DS}). Quitte à remplacer $R$ par $R'$, on peut  supposer que   les torseurs $\ec$ et $\ec'$ sont triviaux sur $\Spec(B)$ et même, d'après le lemme \ref{lem:cjgenerique}ci-dessous qu'il existe des trivialisations qui vérifient l'assertion2.

  \end{preuve}

 \begin{lemme}\label{lem:cjgenerique}
   Soit $a\in \car^{\reg}(B)$. Soit $X$ et $Y$ deux éléments de $\ggo(B)$ tels que 
$$\chi_G(X)=\chi_G(Y).$$
Il existe une extension finie $K'$ de $K$ telle que  $X$ et $Y$ sont conjugués par un élément de $G(B')$, où $B'$ est la clôture intégrale de $B$ dans le corps des fonctions de $C_{K'}$.
 \end{lemme}

 \begin{preuve}
   Pour tout $B$-anneau $A$ soit 
$$\tc(A)=\{g\in G(A) \ | \ \Ad(g)X=Y\}.$$
Rappelons que $B$ a pour corps des fractions le corps $F$ des fonctions de $C_K$ et pour corps résiduel  le corps de fonctions $\kappa(C)$ de $C_\kappa$. D'après le lemme \ref{lem:cjsep},  l'ensemble $\tc(\kappa(C))$ est non vide et, quitte à remplacer $K$ par une extension finie, on peut supposer qu'il en est de même pour $\tc(F)$. On en déduit que le foncteur $\tc$ est représenté par un torseur encore noté $\tc$ sous le schéma en groupes $T_X$ qui centralise $X$. Or vu l'hypothèse de   $a\in \car^{\reg}(B)$, ce schéma en groupes est un schéma en tores. Comme le torseur $\tc$ est trivial en fibre générique $\Spec(F)$, il est trivial globalement (cf. par exemple \cite{CTS2} proposition 2.2). Il possède donc une section sur $\Spec(B)$.
 \end{preuve}

Désormais on fixe $X\in \ggo(B)$ et des trivialisations de $\ec$ et $\ec'$ sur $\Spec(B)$ qui vérifient les assertions du lemme \ref{lem:specB}. On obtient alors une trivialisation de $\phi_K$ qui s'identifie à la translation à gauche par un élément $\delta\in G(F)$ qui vérifie
$$\Ad(\delta^{-1})X=X.$$

Il résulte du lemme \ref{lem:reductrait} et des considérations ci-dessus  qu'on a le lemme suivant.

\begin{lemme}\label{lem:deltadansGB}
  Le morphisme $\phi_K$ se prolonge si et seulement si $\delta\in G(B)$.
\end{lemme}

Comme la restriction de $a$ définit un élément de $\car^{\reg}(B)$, la réduction modulo $\pi$ de $X$ est encore semi-simple régulière. Le centralisateur de $X$ dans $G\times_k B$ est donc un schéma en tores sur $B$. Comme $B$ est normal, on sait qu'un tel schéma en tores est isotrivial. Autrement dit, le tore $T_{X,F}$ se déploie sur une extension séparable et non ramifiée. On en déduit qu'il existe $\la\in X_*(T_X)$ et un élément $h\in T_X(F)\cap G(B)$ tel que 
$$\delta=\la(\pi) h.$$ 
Comme $h\in G(B)$ centralise $X$, on peut  composer la trivialisation de $\ec$ avec $h$ sans que la condition 2 du lemme \ref{lem:specB} soit affectée.  Quitte à changer de trivialisation, on peut et on va supposer que $h=1$. Notons que $\la$ est nécessairement fixe sous $\Gal(F_s/F)$. La caractéristique $\chi_L(X)$ est égale à la restriction de $a_L$ à $\Spec(B)$. Comme $(a_L,t)$ définit en fibre générique un point de  $\Ac_{L,\el}(K)$, le tore  $T_{X,F}$ est un sous-$F$-tore elliptique de $L$ (cf. corollaire \ref{prop:caractdesAM} et lemme \ref{lem:cjsep}) ce qui se traduit par l'égalité
  $$X_*(T_X)^{\Gal(F_s/F)}=X_*(A_L)$$
où $A_L$ est le centre connexe de $L$. Ainsi $\la\in X_*(A_L)$. Le lemme suivant, combiné au lemme \ref{lem:deltadansGB}, montre l'existence du prolongement de $\phi_K$.

\begin{lemme} \label{lem:la=0}Sous les hypothèses ci-dessus, on a $\la=0$.
  \end{lemme}

\end{paragr}

\begin{paragr}[Preuve du lemme \ref{lem:la=0}.] --- On continue avec les notations du paragraphe précédent. Notons tout d'abord que le lemme est évident si $L=G$ puisque dans ce cas $X_*(A_G)=0$. Dans la suite, on suppose donc qu'on a $L\not=G$. On raisonne par l'absurde en supposant $\la\not=0$ et on cherche une contradiction. Soit $P$, resp. $\bar{P}$, le sous-groupe parabolique dans $\fc(L)$  défini par la condition suivante : les racines $\al$ de $T$ dans $P$ satisfont l'inégalité $\al(\la)\geq 0$, resp. $\al(\la)\leq 0$.  Ces sous-groupes sont opposés au sens où l'intersection $P\cap \bar{P}$ un sous-groupe de Levi commun. . On notera que cette intersection contient le sous-groupe de Levi $M$ défini au \S\ref{S:trait} et que ces sous-groupes paraboliques sont propres puisque $\la$ n'est pas nul. Soit $r$ le plus grand entier $\al(\la)$ lorsque $\al$ parcourt l'ensemble $\Phi_T^G$ des racines de $T$ dans $G$. Soit
$$\ggo=\bigoplus_{\al\in X^*(T)}\ggo_\al$$
le décomposition de $\ggo$ en espaces propres pour l'action de $T$. Pour tout entier $i$ qui vérifie $- r\leq i \leq r$, soit
$$\ggo_i=\bigoplus_{\{\al\in X^*(T) \ |\ \al(\la)=i \}} \ggo_\al$$
et
$$\ggo_+=\bigoplus_{\{-r<i\leq r \}} \ggo_i.$$
Les sous-espaces $\ggo_+$ et $\ggo_r$ sont stables par $P$ et $\ggo_-r$ est stable par $\bar{P}$.

L'inclusion
\begin{equation}
  \label{eq:uneinclusion}
    \pi^r \Ad(\la(\pi))\ggo(B)\subset \ggo(B)
  \end{equation}
  est une inclusion entre deux sous-$B$-modules de $\ggo(F)$, libres et de rang maximal, qui fournissent chacun un  prolongement à $C_R$ du fibré vectoriel $\Ad_D(\ec_K)$ sur  $C_K$ (cf. proposition 6 de \cite{Langton}). L'un, celui associé à $\ggo(B)$,  n'est autre que $\Ad_D(\ec)$ et on note l'autre $\vc$. Dans la suite, on note $\ec_\kappa$ et $\vc_\kappa$ le $G$-torseur et le fibré vectoriel sur $C_\kappa=C\times_k \kappa$ déduits de $\ec$ et $\vc$ par changement de base. On déduit de l'inclusion (\ref{eq:uneinclusion})  un morphisme de fibrés vectoriels
  \begin{equation}
    \label{eq:vctoAd}
    \vc \to \Ad_D(\ec).
  \end{equation}
  Soit $\mathcal{W}$ et $ \mathcal{Q}$ le noyau et l'image du  morphisme $\vc_\kappa \to  \Ad_D(\ec_\kappa)$.  Notons que $\mathcal{W}$ est localement facteur direct de $\vc_\kappa$. Ce n'est pas le cas pour $ \mathcal{Q}$ qui est simplement un sous-faisceau localement libre de  $\Ad_D(\ec_\kappa)$. On a 
$$\deg(\vc_\kappa)=0.$$ 
En effet, par platitude on a $\deg(\vc_\kappa)=\deg(\vc_K)$ et d'autre part $\deg(\vc_K)=\deg(\Ad_D(\ec_K))=0$. Comme $\deg(\vc_\kappa)=0$, on a 
\begin{equation}
  \label{eq:degWQ}
  \deg(\mathcal{W})=-\deg(\mathcal{Q} )
\end{equation}

Soit  
$$\Ad_D(\ec')\to \vc$$
l'unique isomorphisme qui, en fibre générique,  est le composé  de l'isomorphisme 
$$\Ad_D(\phi_K) \ : \ \Ad_D(\ec'_K)\to \Ad_D(\ec_K)$$
par l'homothétie de rapport $\pi^r$ et qui sur $\Spec(B)$ coïncide avec  l'isomorphisme 
$$\ggo(B) \to  \pi^r \Ad(\la(\pi))\ggo(B)$$ 
donné par $\pi^r\Ad(\la(\pi))$. En composant cet isomorphisme avec le morphisme $\vc\to \Ad_D(\ec)$ défini en (\ref{eq:vctoAd}), on obtient un morphisme de fibrés vectoriels sur $C_R$
\begin{equation}
  \label{eq:unprolongement}
  \Ad_D(\ec')\to \Ad_D(\ec),
\end{equation}
qui donne, en fibre spéciale, un morphisme
\begin{equation}
  \label{eq:fibspeciale}
  \Ad_D(\ec'_\kappa)\to \Ad_D(\ec_\kappa)
\end{equation}
dont on note $\mathcal{W}'$ le noyau. Par construction, on a $\mathcal{W}'\simeq \mathcal{W}$ et l'image de ce morphisme est le fibré vectoriel $\Qc$. Vu (\ref{eq:degWQ}), on a
\begin{equation}
  \label{eq:degWQ'}
  \deg(\mathcal{W}')=-\deg(\mathcal{Q} )
\end{equation}

On a fixée au  \S \ref{S:trait} des trivialisations de $\ec$ et $\ec'$ sur $\Spec(B)$ qui vérifient les assertions du lemme \ref{lem:specB}. Via ces choix, le sous-groupe parabolique $P$ de $G$ définit une réduction générique de $\ec'_\kappa$ à $P$ qui  s'étend automatiquement en une réduction  notée $\ec'_{P,\kappa}$ de  $\ec'_{\kappa}$ à $P$ sur $C_\kappa$. De même, on note $\ec_{\bar{P},\kappa}$ la réduction de $\ec_\kappa$ à $\bar{P}$ qui prolonge la réduction générique donnée par $\bar{P}\subset G$.

\begin{lemme}\label{lem:egalite-inj}
  Avec les notations ci-dessus, on a
  \begin{enumerate}
  \item une égalité $\ec'_{P,\kappa}\times^{P,\Ad}_k \ggo_+=\mathcal{W}' \ ;$
  \item une injection $\mathcal{Q} \hookrightarrow \ec_{\bar{P},\kappa} \times_k^{\bar{P},\Ad} \ggo_{-r}$ entre fibrés vectoriels de même rang. 
  \end{enumerate}
\end{lemme}

\begin{preuve}
  Comme $\ec'_{P,\kappa}\times^{P,\Ad}_k \ggo_+$ et $\mathcal{W}'$ sont des sous-fibrés vectoriels de 
$\Ad_D(\ec_\kappa')$, il suffit de vérifier 1 au point générique de $C_\kappa$. De même, comme $\mathcal{Q}$ et  $\ec_{\bar{P},\kappa} \times_k^{\bar{P},\Ad} \ggo_{-r}$ sont des sous-faisceaux localement libres de $\Ad_D(\ec_\kappa')$ et que  $\ec_{\bar{P},\kappa} \times_k^{\bar{P},\Ad} \ggo_{-r}$ est même localement facteur direct, il suffit de vérifier 2  au point générique de $C_\kappa$. Or, dans  les trivialisations $\ec$ et $\ec'$ sur $\Spec(B)$ qu'on a fixée au \S \ref{S:trait}, le morphisme (\ref{eq:unprolongement}) est le morphisme $\ggo(B) \to \ggo(B)$ donnée par $\pi^r\Ad(\la(\pi))$. En particulier, en réduction modulo $\pi$, ce morphisme est la projection de $\ggo(\kappa(C))$ sur $\ggo_{-r}(\kappa(C))$ de noyau  $\ggo_+(\kappa(C))$. Le lemme est alors clair.
\end{preuve}

\begin{lemme} \label{lem:inegdegres} Soit $\mu\in X^*(P)$ le déterminant de l'action adjointe de $P$ sur $\ggo_+$. On a l'inégalité 
\begin{equation}
  \label{eq:inegWQ}
  \deg(\mu(\ec'_{P,\kappa}))\geq \deg(\mu(\ec_{\bar{P},\kappa})).
\end{equation}
  \end{lemme}

  \begin{remarque} Le caractère $\mu$ est un plus haut poids dans la représentation de $G$ sur $\wedge^{\dim(\ggo_+)}\ggo$. Il s'écrit donc dans la base $\hat{\Delta}_P$ avec des coefficients positifs dont l'un au moins est non nul.
      \end{remarque}

\medskip

  \begin{preuve} D'après l'assertion 1 du lemme \ref{lem:egalite-inj}, on a 
$$\deg(\mathcal{W}')=\deg(\ec'_{P,\kappa}\times^{P,\Ad}_k \ggo_+)=\deg(\mu(\ec'_{P,\kappa})).$$
D'après  l'assertion 2 du lemme \ref{lem:egalite-inj}, on a
$$\deg(\mathcal{Q})\leq   \deg(\ec_{\bar{P},\kappa} \times_k^{\bar{P},\Ad} \ggo_{-r})=-\deg(\mu(\ec_{\bar{P},\kappa})).$$
On peut alors conclure puisque $\deg(\mathcal{W}')=-\deg(\mathcal{Q})$ (cf. (\ref{eq:degWQ'}).
      \end{preuve}

On a défini au \S\ref{S:trait} un  sous-groupe de levi $M$. Comme $P$ et $\bar{P}$ contiennent tous deux $M$, il résulte de la proposition \ref{prop:Preduc} et de sa démonstration que  les torseurs $\ec'_{P,\kappa}$ et $\ec_{\bar{P},\kappa}$ sont les premiers facteurs de triplets $m'_{P,\kappa}$ et $m_{\bar{P},\kappa}$ qui sont des réductions à $P$ et $\bar{P}$ de $m_\kappa'$ et $m_\kappa$. Or, par hypothèse, ces triplets de Hitchin sont $\xi$-stables. On a donc, d'une part,
$$\xi\in -\deg(m'_{P,\kappa}) - \,^+\ago_P+ \ago_T^M$$
(rappelons que $\ago_G=\{0\}$ puisque $G$ est semi-simple) d'où (cf. la remarque qui suit le lemme \ref{lem:inegdegres}) 
$$\mu(\xi) < - \mu(\deg(m'_{P,\kappa}))=-\deg(\mu(\ec'_{P,\kappa}))$$
et d'autre part
$$\xi\in  -\deg(m_{\bar{P},\kappa}) - \,^+\ago_{\bar{P}}+ \ago_T^M$$
ce qui implique l'inégalité 
$$ -\mu(\xi) < \mu(\deg(m_{\bar{P},\kappa}))=\deg(\mu(\ec_{\bar{P},\kappa})).$$
On a donc
$$\deg(\mu(\ec'_{P,\kappa}))< -\mu(\xi) < \deg(\mu(\ec_{\bar{P},\kappa}))$$
ce qui contredit l'inégalité (\ref{eq:inegWQ})du lemme \ref{lem:inegdegres}. C'était la contradiction recherchée au début de ce paragraphe.
\end{paragr}

\section{$\mc^\xi$ est un champ de Deligne-Mumford}\label{sec:DM}

\begin{paragr} L'objet de cette section est de prouver le théorème suivant.

  \begin{theoreme}\label{thm:DM}
Pour tout groupe  semi-simple $G$ et tout paramètre $\xi$ en position générale (cf. définition \ref{def:posgenerale}), le champ $\mc^\xi_G$ est un champ de Deligne-Mumford.
  \end{theoreme}

  La démonstration se trouve au paragraphe suivant.
\end{paragr}

\begin{paragr} On fixe un  groupe  semi-simple $G$ et un paramètre $\xi$ en position générale. Soit $(\ec,\theta,t)\in \mc^\xi(k)$ et $\Aut_G(\ec,\theta)$ le sous-schéma en groupes de $\Aut_G(\ec)$ qui centralise $\theta$. Le foncteur qui, à un $k$-schéma $S$, associe le groupe des sections sur $C\times_k S$ du schéma en groupes $\Aut_G(\ec,\theta)\times_k S$ est représenté par un schéma en groupes sur $k$. Il s'agit de voir que ce groupe est non ramifié, ou encore, par un critère différentiel, que son espace tangent sur $k$ est trivial. Or cet espace tangent admet la description suivante : c'est l'espace des sections $\varphi\in H^0(C,\Ad(\ec))$ qui commutent à $\theta$ c'est-à-dire telles que le crochet $[\theta,\varphi]$ est nul. La nullité de l'espace tangent résulte donc de la proposition suivante.

  \begin{proposition}
    Pour toute section $\varphi\in H^0(C,\Ad(\ec))$ telle que $[\theta,\varphi]=0$ on a
$$\varphi=0.$$
  \end{proposition}

  \begin{preuve} Soit  $\varphi$ comme ci-dessus. Notons que,  par le lemme évident ci-dessous,  la caractéristique de $\varphi$ est constante.

    \begin{lemme}
      L'espace des sections $H^0(C,\car_G)$ s'identifie à $\car_G(k)$.
    \end{lemme}

Soit $M \in \lc$ tel que $f(\ec,\theta,t)\in \chi^M_G(\Ac_{M,\el}(k))$ et $F$ le corps de fonctions de la courbe $C$.

\begin{lemme}\label{lem:trivphi}
Il existe une trivialisation générique de $\ec$ qui vérifie les deux assertions suivantes :
\begin{enumerate}
\item l'identification de $\Ad_D(\ec)$ avec $\ggo(F)$ qui en résulte envoie $\theta$ sur un élément $X\in \mgo(F)$ qui est semi-simple, $G$-régulier et elliptique dans $M$ et qui vérifie l'assertion 1 de la définition \ref{def:Cata} ;
\item l'identification de $\Ad(\ec)$ avec $\ggo(F)$ qui en résulte envoie  $\varphi$ sur un élément $Y\in \tgo(k)$.
\end{enumerate}
\end{lemme}

\begin{preuve}
  On a vu dans la preuve de la proposition \ref{prop:fibre} qu'il existe une  trivialisation générique de $\ec$ qui vérifie l'assertion 1. Dans l'identification de $\Ad(\ec)$ avec $\ggo(F)$, la section  $\varphi$ s'envoie sur un élément de $\ggo(F)$ noté $Z$. Comme $Z$ commute à $X$ et que le centralisateur de $X$ dans $\ggo$ est inclus dans $\mgo$, on en déduit que $Z$ appartient à $\mgo(F)$. La caractéristique $\chi_M(Z)\in \car_M(F)$ s'envoie par $\chi_G^M$ sur la caractéristique de $\varphi$ qui par le lemme précédent est un élément de $\car_G(k)$. Il s'ensuit que $\chi_M(Z)$ appartient à $\car_M(k)$ et qu'on peut trouver un élément $Y\in \tgo(k)$ tel que $\chi_M(Y)=\chi_M(Z)$. Par une variante du lemme \ref{lem:cjsep}, on en déduit que $Y$ et $Z$ sont conjugués par un élément de $M(F)$. Le lemme est alors évident.
  \end{preuve}

On fixe désormais une  trivialisation générique de $\ec$ qui vérifie le lemme \ref{lem:trivphi}. Soit $X\in \mgo(F)$ et $Y\in \tgo(k)$ les éléments qui s'en déduisent. Soit $L$ le centralisateur de $Y$ dans $G$. On vérifie que $L$ est un élément de $\lc^G(T)$ qui contient $M$. Soit $Q$ un sous-groupe parabolique de $G$ de Levi $L$. On notera le lemme suivant, dont la démonstration est laissée au lecteur.

\begin{lemme}\label{lem:NQ-nQ}
  Le morphisme $N_Q \to \ngo_Q$ donné par $n\mapsto \Ad(n^{-1})Y-Y$ est un isomorphisme.
\end{lemme}

Soit $V$ l'ensemble des points fermés de $C$. On reprend les notations de la section \ref{sec:description}. On déduit de la trivialisation générique du torseur $\ec$ et de l'existence d'une trivialisation sur un voisinage formel de tout point $v\in V$ une famille $(g_v)_{v\in V}$ telle que 
\begin{enumerate}
\item $g_v\in G((z_v))(k)/ G[[z_v]](k)$ est trivial sauf pour un nombre fini de $v$ ;
\item $\Ad(g_v^{-1})Y \in \ggo[[z_v]](k)$ ;
\item  la projection de $\xi$ sur $\ago_M^G$ appartient à la projection sur $\ago_M^G$ de l'enveloppe convexe des points 
$$-H_P((g_v)_{v\in V})=-\sum_{v\in V} H_P(g_v).$$
\end{enumerate}
La condition 2 traduit le fait que $\varphi$ est une section globale de $\Ad(\ec)$ et la condition 3 la $\xi$-stabilité du triplet $(\ec,\theta,t)$ (cf. condition (d) de la définition \ref{def:Cata}).

\begin{lemme}
  Pour tout $v\in V$, la classe $g_v$ se relève en un élément de $L((z_v))(k)$.
\end{lemme}

\begin{preuve}
  Soit $Q$ un sous-groupe parabolique de $G$ de Levi $L$.  Par la décomposition d'Iwasawa, l'élément $g_v$ se relève en un élément $l_v n_v $ avec $l_v\in L((z_v))(k)$ et $n_v \in N_Q((z_v))(k)$. Le fait que $L$ centralise $Y$ et la relation 1 ci-dessus entraînent alors qu'on a 
$$\Ad(n_v^{-1})Y\in \ggo[[z_v]](k)$$
d'où 
$$\Ad(n_v^{-1})Y-Y\in \ngo_Q[[z_v]](k).$$ 
Par le lemme \ref{lem:NQ-nQ}, on a alors $n_v\in  N_Q[[z_v]](k)$. D'où le lemme. 
\end{preuve}

Le lemme précédent implique que le triplet $(\ec,\theta,t)$ possède  une ``réduction'' au sous-groupe de Levi $L$. Il implique aussi que la projection sur $\ago_L^G$ de l'enveloppe convexe de 3 ci-dessus est réduite à un point qui est nécessairement un point du réseau $X_*(L)$ défini au \S\ref{S:espacesa}. La $\xi$-stabilité entraîne que la projection de $\xi$ sur $\ago_L^G$ appartient à ce réseau. Si $L\neq G$, cela contredit le fait que $\xi$ soit en position générale. On a donc $L=G$ et $Y=0$ d'où $\varphi=0$. 
\end{preuve}
\end{paragr}

\section{Comptage des points rationnels dans une fibre}\label{sec:comptage}

\begin{paragr}[Le théorème principal.] ---   Voici la situation qui va nous occuper jusqu'à la fin de cette section. Rappelons qu'on a fixé en \S \ref{S:lacourbe} une courbe $C$ sur $k$ ainsi qu'un diviseur $D$ et un point $\infty\in C(k)$. On suppose que $C$ provient par changement de base d'une courbe  $C_0$ sur $\Fq$. Soit $\tau$ l'automorphisme de Frobenius de $k$ donné par $x\mapsto x^q$. On note encore $\tau$ l'automorphisme $1\times \tau$ de $C=C_0\times_{\Fq} k$. On suppose que $D$ et $\infty$ sont fixes sous $\tau$.

Soit $F$ le corps des fonctions de $C$. L'automorphisme $\tau$ de $C$  induit un automorphisme de $F$ encore noté $\tau$ dont l'ensemble des points fixes est précisément le corps $F_0$ de la courbe $C_0$. Soit $V$, resp. $V_0$, l'ensemble des points fermés de $C$, resp. $C_0$. Pour tout $v_0\in V_0$, soit $F_{0,v_0}$ le complété de $F_0$ en  $v_0$ ; on note simplement $\tau$ l'automorphisme $\tau\otimes 1$ de $$F\otimes_{F_0} F_{0,v_0}=\prod_{v\in V, \ v|v_0}F_v$$
où $F_v$ est le complété de $F$ en $v$. On identifie $F_v$ à $k((z_v))$ et $F_{0,v_0}$ à $\Fq((z_{v_0}))$ par le choix d'uniformisantes. Soit $\AAA$ l'anneau des adèles de $F$. Cet anneau est naturellement muni d'une action de $\tau$ pour laquelle le sous-anneau des points fixes sous $\tau$ s'identifie à l'anneau $\AAA_0$ des adèles de $F_0$.

On suppose que que le groupe $G$ est \emph{semi-simple} et que le couple $(G,T)$ (cf. \S \ref{sec:morcar}) provient d'un couple analogue $(G_0,T_0)$ défini sur $\Fq$ tel que le tore $T_0$ soit déployé sur $\Fq$. Dans ce cas, tous les sous-groupes de Levi ou paraboliques de $G$ qui contiennent $T$ proviennent d'objets analogues définis sur $\Fq$ relatifs à $(G_0,T_0)$. On note encore $\tau$ l'automorphisme de $G=G_0\times_{\Fq}k$ donnée par $1\times \tau$. On en  déduit un automorphisme toujours noté $\tau$ des  groupes de points de $G$ à valeurs dans $F$, $\AAA$ ou  $F\otimes_{F_0} F_{0,v_0}$ pour tout  $v_0\in V_0$. 

Comme le diviseur $D$ sur $C$ est fixe par $\tau$, il se descend en un diviseur $D_0$ sur $C_0$ de la forme 
$$\sum_{v\in V_0} d_v v.$$
Soit $\mathbf{1}_{D_0}$ la fonction sur $\ggo(\AAA_0)$ qui est la fonction caractéristique de l'ensemble
$$\prod_{v\in V_0} z_v^{-d_v} \ggo[[z_v]](\Fq).$$
Pour alléger un peu les notations, on pose
$$K= \prod_{v\in V_0} G[[z_v]](\Fq),$$
c'est un sous-groupe ouvert et compact de $G(\AAA_0)$. De plus, la fonction $\mathbf{1}_{D_0}$ est invariante sous l'action adjointe de $K$.

Pour tout $\delta\in M(F)$, soit $\tau_\delta$ l'automorphisme de $G(F)$, resp. de $\ggo(F)$, donné par $\Int(\delta)\circ\tau$, resp. $\Ad(\delta)\circ \tau$. Le groupe $G$ muni de l'automorphisme $\tau_\delta$ définit par descente un groupe noté $G^\delta$ sur $\Fq$.

L'automorphisme de Frobenius $\tau$ agit également sur les espaces $\Ac$ et $\Ac_M$ pour tout $M\in \lc$. Soit $(a,t)\in \Ac(k)$ un élément fixe par $\tau$. Il existe $M\in \lc$ et $(a_M,t)\in \Ac_{M,\el}(k)$ tels que 
$$(a,t)=\chi^M_G((a_M,t)).$$
D'après les propositions \ref{prop:immersionfermee} et \ref{prop:reunionAM}, de tels éléments existent et sont uniques : en particulier ils sont fixes sous $\tau$.

Voici le principal théorème de cette section. Il exprime le nombre de points rationnels d'une fibre  de $\overline{f^\xi}$ en termes \emph{d'intégrales orbitales pondérées} d'Arthur.

\begin{theoreme}
    \label{thm:comptage}
Soit $\xi$ \emph{en position générale}. Soit $M\in \lc$ un sous-groupe de Levi, $(a_M,t)\in \Ac_{M,\el}(k)$  fixe par $\tau$ et $(a,t)=\chi^M_G((a_M,t))$. Le cardinal du groupoïde des points rationnels de la fibre de $\overline{f^\xi}$ en $(a,t)$
 est égal à
$$ \vol(\ago_M/X_*(M))^{-1}\cdot\vol(a,t)\cdot  \sum_{ h} \, \sum_{X} \, J^G_M(\Ad(h^{-1})X,\mathbf{1}_{D_0})$$ 
où
\begin{itemize}
\item  $h$  parcourt l'ensemble des doubles classes 
$$h\in M(F)\back M(\AAA) / M(\AAA_0)$$
qui vérifient
$$\delta=h\tau(h)^{-1} \in M(F) \ ;$$
\item  $X$ parcourt un système de représentants de l'ensemble des classes de $M^\delta(F_0)$-conjugaison dans l'ensemble 
$$\{X\in\mgo^\delta(F_0)\ | \ \chi_M(X)={a_{M}}_{|\Spec(F_0)}  \}$$

des $X\in\mgo^\delta(F_0)$ dont la caractéristique $\chi_M(X)$ est égale à la restriction de $a_M$ à $Spec(F_0)$ ;
\item $J^G_M(\Ad(h^{-1})X,\mathbf{1}_{D_0})$ est l'intégrale orbitale pondérée d'Arthur définie au \S\ref{S:poidsArthur} l.(\ref{eq:IOPv}) ;
\item le volume $\vol(a,t)$ est défini au \S\ref{S:Haar} l.(\ref{eq:vol(a,t)}) ;
\end{itemize}

 \end{theoreme}

\begin{remarque}
  Les intégrales orbitales pondérées considérées ci-dessus dépendent de choix de mesures de Haar sur $\ago_M$ (cf. \S\ref{S:poidsArthur}) et sur certains tores (cf. \S\ref{S:Haar}). Comme ces choix interviennent aussi dans les facteurs $ \vol(\ago_M/X_*(M))$ et $\vol(a,t)$, la somme  qui apparaît dans le théorème \ref{thm:comptage} ne dépend d'aucun choix.
\end{remarque}

\medskip
Avant d'entamer la démonstration du théorème \ref{thm:comptage} qui va nous occuper jusqu'à la fin de cette section, nous allons donner quelques rappels sur le cardinal d'un groupoïde quotient et sur ses points fixes sous un automorphisme.
\end{paragr}

\begin{paragr}[Cardinal d'un groupoïde.] --- \label{S:cardinal}Soit $\xc$ un ensemble et $G$ un groupe abstrait qui agit à gauche sur $\xc$. Le cardinal du groupoïde quotient $[G\back \xc']$ (cf. l(\ref{eq:groupoidequot}) du \S\ref{S:description-adelique}) est l'élément $\card([G\back \xc])\in \NN\cup \{\infty\}$ défini par
$$\card([G\back \xc])=\sum_{x\in G\back \xc} \frac{1}{|\stab_G(x)|}$$
où 
\begin{itemize}
\item $G\back \xc$ est (un système de représentants de) l'ensemble des orbites de $\xc$ sous $G$ ;
\item $\displaystyle \frac{1}{|\stab_G(x)|}$ vaut $0$ si le stabilisateur de $x$ dans $G$ est infini et l'inverse de son cardinal sinon.
\end{itemize}
  
\end{paragr}

\begin{paragr}[Points fixes sous un automorphisme.] --- \label{S:ptfixe}On continue avec les notations du paragraphe précédent. Soit $\tau$ un automorphisme de $G$. Soit une bijection de $\xc$, encore notée $\tau$, compatible à l'automorphisme de $G$ au sens où l'on a 
$$\tau(g.x)=\tau(g).\tau(x)$$
  Par définition, le groupoïde des points fixes sous $\tau$ de $[G\back \xc]$ est le groupoïde $[G\back \xc^\tau]$ où
  \begin{itemize}
  \item $\xc^\tau$ est l'ensemble des couples $(x,g)\in \xc\times G$ tels que 
$$g.\tau(x)=x \ ;$$
\item $G$ agit à gauche sur $\xc^\tau$ par 
$$h.(x,g)=(hx,hg\tau(h)^{-1}).$$
  \end{itemize}
\end{paragr}

\begin{paragr}[Première étape.] --- Dans la suite, on fixe $M\in \lc$ et des couples $(a_M,t)$ et $(a,t)$ comme dans l'énoncé  du   théorème \ref{thm:comptage}. Soit $\xi\in \ago_T$ un élément \emph{quelconque}. La fibre du morphisme de Hitchin $\overline{f^\xi}$ au-dessus de $(a,t)$ s'identifie par la proposition  \ref{prop:fibre} au groupoïde quotient $[M(F)\back\xc_{(a,t)}^\xi ]$ où l'ensemble $\xc_{(a,t)}^\xi$ est celui de la définition \ref{def:Cata} du paragraphe \ref{S:description-adelique}. Nos hypothèses entraînent que $\tau$ induit une bijection de cet ensemble compatible à l'action de $\tau$ sur $M(F)$ (au sens du \S\ref{S:ptfixe}). On cherche une expression pour le nombre de points sur $\Fq$ de cette fibre, c'est-à-dire le cardinal du groupoïde des points fixes sous $\tau$ de $[M(F)\back\xc_{(a,t)}^\xi ]$, en termes d'intégrales orbitales pondérées.

Suivant le paragraphe \ref{S:ptfixe}, on introduit l'ensemble $\xc_{(a,t)}^{\xi,\tau}$ formé des triplets $(X,(g_v)_{v\in V},\delta)$ tels que $(X,(g_v)_{v\in V})\in \xc_{(a,t)}^{\xi}$ et $\delta\in M(F)$ vérifient les relations 
  \begin{equation}
    \label{eq:relation1}
    \Ad(\delta)\tau(X)=X
  \end{equation}
  et
\begin{equation}
    \label{eq:relation2}
    \delta\tau((g_v)_{v\in V})=(g_v)_{v\in V}.
  \end{equation}

Le groupe $M(F)$ agit sur les deux premiers facteurs par l'action décrite au \S \ref{S:description-adelique} et par $\tau$-conjugaison sur le deuxième facteur (c'est-à-dire un élément $\gamma\in M(F)$ agit sur $\delta\in M(F)$ par $\gamma \delta \tau(\gamma)^{-1}$). La première étape du comptage consiste à décrire un système de représen\-tants des orbites de $\xc_{(a,t)}^{\xi,\tau}$ sous l'action de $M(F)$. Pour cela, on introduit quelques notations supplémen\-taires.

Pour tout $g \in G(\AAA)$, soit $\mathbf{1}_{M,g}$ la fonction sur $\ago_M$ caractéristique de l'enveloppe convexe des points $-H_P(g)$ pour $P\in \pc(M)$. La définition \ref{def:HP} du \S\ref{S:HP}, appliquée au groupe $M$ et à son sous-groupe parabolique $P=M$, donne un morphisme
$$H_M \ :\ M(\AAA) \to \ago_M.$$
Soit $\xi_M$ la projection de $\xi$ sur $\ago_T$ suivant la décomposition $\ago_T=\ago_M\oplus \ago_T^M$.

\begin{lemme}\label{lem:systrep}
  Il existe une bijection entre l'ensemble $M(F)\back \xc_{(a,t)}^{\xi,\tau}$ et l'ensemble des triplets 
$$(X, g, \delta=h\tau(h)^{-1})$$
qui vérifient les trois conditions suivantes 
\begin{enumerate}
\item  $h$  parcourt l'ensemble des doubles classes 
$$h\in M(F)\back M(\AAA) / M(\AAA_0)$$
qui vérifient
$$\delta=h\tau(h)^{-1} \in M(F) \ ;$$

\item $X$ parcourt un système de représentants de l'ensemble des classes de $M^\delta(F_0)$-conjugaison dans l'ensemble 
$$\{X\in\mgo^\delta(F_0)\ | \ \chi_M(X)={a_{M}}_{|\Spec(F_0)}  \}$$

des $X\in\mgo^\delta(F_0)$ dont la caractéristique $\chi_M(X)$ est égale à la restriction de $a_M$ à $\Spec(F_0)$ ;

\item $g$ parcourt l'ensemble  des doubles classes 
$$ \Int(h^{-1}) T_X(F_0) \back G(\AAA_0) /  K,$$
(où $T_X\subset M^\delta$ est le centralisateur de $X$ dans $G^\delta$) qui vérifient les deux conditions 
\begin{enumerate}
\item  $\mathbf{1}_{D_0}(\Ad(g^{-1}) \Ad(h^{-1})X)=1$ ;
\item $\mathbf{1}_{M,g}(\xi_M+H_M(h))=1$.
\end{enumerate}
\end{enumerate}
De plus, dans cette bijection, le cardinal du stabilisateur dans $M(F)$ d'un élément de $M(F)\back \xc_{(a,t)}^{\xi,\tau}$ qui correspond au triplet  $(X, g, h\tau(h)^{-1})$ est donné par le cardinal de l'ensemble fini 
$$ \big(\Int(h^{-1}) T_X(F_0)\big)\cap gK g^{-1}.$$
  
\end{lemme}

\begin{preuve}
  Partons d'un triplet $(X,(g_v)_{v\in V},\delta)$ dans $\xc_{(a,t)}^{\xi,\tau}$. Rappelons que $g_v$ désigne une classe dans $G((z_v))(k)/ G[[z_v]](k)$.  Soit un élément $g\in G(\AAA)$ tel que pour tout $v$ la composante en $v$ relève $g_v$. Par la relation (\ref{eq:relation2}) ci-dessus, on a
 $$g^{-1} \delta \tau(g)\in \prod_{v\in V} G[[z_v]](k).$$
Par le théorème de Lang (appliqué aux quotients (connexes) du $k$-groupe pro-algébrique  $G[[z_v]]$), on sait qu'un tel élément s'écrit  $ x \tau(x^{-1})$  pour un certain $x\in  \prod_{v\in V} G[[z_v]](k)$. Quitte à remplacer  $g$ par $gx$, on peut et on va supposer qu'on a
\begin{equation}
  \label{eq:underline}
  g^{-1} \delta \tau(g)=1
\end{equation}
autrement dit la classe de $\tau$-conjugaison de $\delta$ dans $G(\AAA)$ est triviale. Le lecteur  vérifiera que cela implique que   la classe de $\tau$-conjugaison de $\delta$ dans $M(\AAA)$ est triviale. Il existe donc  $h \in M(\AAA)$ tel que 
$$\delta=h\tau(h)^{-1}.$$
Posons $g_0=h^{-1}g$. Alors la relation (\ref{eq:underline}) se traduit par 
$$g_0\in G(\AAA_0).$$

D'après la relation (\ref{eq:relation1}), l'élément $X$ appartient à $\mgo^\delta(F_0)$. Soit $t_a\in \tgo^{\reg}[[z_\infty]](k)$ l'unique relèvement de $t$ dont la caractéristique est $a$. De même, on peut définir un point $t_{a_M}$ associé au couple $(a_M,t)$. Mais par unicité de $t_a$, on a $t_a=t_{a_M}$. Ainsi la caractéristique $\chi_M(t_a)$ est $a_M$. La condition 1. (bis) de la définition \ref{def:Cata2} entraîne que $X$ est conjugué à $t_a$ par un élément de  $M((z_\infty)(k)$ d'où 
$$\chi_M(X)=\chi_M(t_a)=a_M.$$
Les conditions  2.(b) et  2.(d) des définitions  \ref{def:Cata} et \ref{def:Cata2} se traduisent par les condition 3 (a) et (b) ci-dessus. 

On laisse le soin au lecteur de vérifier que la construction qui à  $(X,(g_v)_{v\in V},\delta)$ associe  $(X, g_0, h\tau(h)^{-1})$ définit par passage au quotient une bijection avec les propriétés annoncées. On notera que le groupe $\Int(h^{-1}) T_X(F_0)$ est bien inclus dans $G(\AAA_0)$.

Le centralisateur de  $(X,(g_v)_{v\in V},\delta)$ dans $M(F)$ est le groupe 
\begin{eqnarray*}
  \{t \in T_X(F_0) \ | \  (tg_v)_{v\in V}=(g_v)_{v\in V}\} &=& T_X(F_0) \cap \prod_{v\in V} g_v G[[z_v]](k) g_v^{-1}\\
&\simeq & \big(\Int(h^{-1}) T_X(F_0)\big)\cap g_0K g^{-1}_0
\end{eqnarray*}
ce qui donne la dernière assertion.
\end{preuve}
\end{paragr}

\begin{paragr}[Mesures de Haar.] --- \label{S:Haar} On munit $G(\AAA_0)$ de la mesure de Haar normalisée par 
$$\vol(K)=1.$$
Plus généralement pour tout sous-groupe parabolique semi-standard $P$ de $G$, les groupes $M_P(\AAA_0)$ et $N_P(\AAA_0)$ sont munis des mesures de Haar normalisées par
$$\vol(M_P(\AAA_0)\cap K)=1 \ \ \ \text{  et   }\ \ \  \vol(N_P(\AAA_0)\cap K)=1.$$

Soit $X_0\in \mgo(F_0)$ dont la caractéristique est égale à la restriction de $a_M$ à $\Spec(F_0)$. Un tel $X_0$ existe par la section de Kostant relative à $M$. On munit $T_{X_0}(\AAA_0)$ d'une mesure de Haar. Soit 
$$T_{X_0}(\AAA_0)^1=T_{X_0}(\AAA_0)\cap \Ker(H_M)\ ;$$
c'est un sous-groupe ouvert de $T_{X_0}(\AAA_0)$ qu'on munit de la mesure induite. On munit $T_{X_0}(F_0)$ de le mesure de comptage et $T_{X_0}(F_0)\back T_{X_0}(\AAA_0)^1$ de la mesure quotient.

\begin{lemme} Le volume du quotient ci-dessous est fini
  \begin{equation}
  \label{eq:volume}
  \vol(T_{X_0}(F_0)\back T_{X_0}(\AAA_0)^1)<\infty.
\end{equation}
\end{lemme}

\begin{preuve}  D'après le  corollaire \ref{cor:caractdesAM}, il existe un élément $X\in \mgo(F)$, semi-simple, $G$-régulier et elliptique dans $\mgo(F)$. D'après le  lemme \ref{lem:cjsep}, $X$ et $X_0$ sont conjugués sous $M(F)$. On en déduit que le plus grand sous-tore  $F_0$-déployé de $T_{X_0}$ est inclus dans le centre connexe de $M\times F_0$  (et en fait égal). Ainsi, le $T_{X_0}$ est $F_0$-elliptique dans $M$ et l'on sait bien que le quotient $T_{X_0}(F_0)\back T_{X_0}(\AAA_0)^1$ est alors compact d'où la finitude du volume.
\end{preuve}

Soit $\delta\in G(F)$ et $X\in \mgo^\delta(F_0)$ de caractéristique $a_M$. La condition sur la caractéristique entraîne que $X_0$ et $X$ sont conjugués par un élément $m\in M(F)$ (cf. lemme \ref{lem:cjsep}). L'automorphisme $\Int(m)$ induit un $F_0$-isomorphisme entre les tores $T_{X_0}$ et $T_X$. Par transport par $\Int(m)$, on déduit de la mesure sur $ T_{X_0}(\AAA_0)$ une mesure de Haar sur  $ T_{X}(\AAA_0)$ pour laquelle on a 

\begin{equation}
  \label{eq:volumesegaux}
\vol(T_{X}(F_0)\back T_{X}(\AAA_0)^1)  =\vol(T_{X_0}(F_0)\back T_{X_0}(\AAA_0)^1).
\end{equation}
On note 
\begin{equation}
  \label{eq:vol(a,t)}
\vol(a,t)
\end{equation}
le volume ci-dessus. 
\end{paragr}

\begin{paragr}[Réseaux dans $\ago_M$.] --- \label{S:reseaux}Soit $G_{\scnx}$ le revêtement simplement connexe de $G$. Soit $L\in \lc^G(T)$ un sous-groupe de Levi semi-standard. Soit $L_{\scnx}$ l'image réciproque de $L$ dans $G_{\scnx}$. Soit $L_{\der}$ le groupe dérivé de $L$ et $L_{\SCNX}$ le revêtement simplement connexe de $L_{\der}$. On ne confondra pas  $L_{\scnx}$ et $L_{\SCNX}$. Le groupe $L_{\SCNX}$ est en fait le groupe dérivé de  $L_{\scnx}$. On note encore $T_?$ le tore obtenu par image réciproque de $T$ dans $?$ où $?$ peut être l'un des trois groupes $L_{\der}$, $L_{\SCNX}$ ou $L_{\scnx}$.
Soit 
$$\pi_1(L_{\der})=\mathrm{coker}(X_*(T_{L_{\SCNX}})\to  X_*(T_{L_{\der}})).$$

Le morphisme de restriction $X^*(L)\to X^*(T)$ et donne dualement un morphisme $X_*(T) \to X_*(L)$ où l'on a posé $X_*(L)=\Hom_{\ZZ}(X^*(L),\ZZ).$

\begin{lemme}\label{lem:DC}
  On a un diagramme commutatif à lignes et colonnes exactes.

$$\xymatrix{  &  0   \ar[d]     &    0   \ar[d]     &   0   \ar[d]     & \\
0 \ar[r] &  X_*(T_{L_{\SCNX}})     \ar[d]    \ar[r]  &  X_*(T_{G_{\scnx}})     \ar[d]    \ar[r] &   X_*(L_{\scnx})  \ar[d]    \ar[r] & 0 \\
0 \ar[r] &  X_*(T_{L_{\der}})     \ar[d]    \ar[r]  &  X_*(T)     \ar[d]    \ar[r] &   X_*(L)  \ar[d]    \ar[r] & 0 \\
0 \ar[r] &  \pi_1(L_{\der})     \ar[d]    \ar[r]  & \pi_1(G)     \ar[d]    \ar[r] &   X_*(L)/X_*(L_{\scnx})  \ar[d]    \ar[r] & 0\\
 &  0       &    0        &   0        &
}$$ 

\end{lemme}

\begin{preuve}
  Traitons l'exactitude de la deuxième ligne. Il s'agit d'une suite de $\ZZ$-module libre dont la suite duale est la suite exacte
$$\xymatrix{0 \ar[r] &  X^*(L)     \ar[r]  &  X^*(T)        \ar[r] &   X^*(T_{L_{\der}})  \ar[r] & 0}.$$
Seule la surjectivité est moins évidente. Elle résulte des isomorphismes 
$$T/T_{L_{\der}}=Z_L^0/ (Z_L^0\cap L_{\der})=L/L_{\der},$$
où $Z_L^0$ est le centre connexe de $L$, et de l'égalité  $ X^*(L)= X^*(L_{\der})$.

L'exactitude de la première ligne se déduit de celle de la deuxième lorsqu'on remplace $L$ par $L_{\scnx}$. L'injectivité des deux premières flèches verticales est bien connue. On rappelle que $X_*(T_{G_{\scnx}})$ est le sous-$\ZZ$-module de $X_*(T)$ engendré par les coracines de $T$ dans $G$. Le reste du diagramme se déduit alors du lemme du serpent et de l'injectivité de 
$$ \pi_1(L_{\der})  \longrightarrow    \pi_1(G)$$
c'est-à-dire de l'inclusion
$$X_*(T_{L_{\der}}) \cap   X_*(T_{G_{\scnx}}) \subset X_*(T_{L_{\SCNX}}).$$
Cette dernière est évidente car $X_*(T_{L_{\SCNX}})$ est d'indice fini dans $X_*(T_{L_{\der}})$ et d'autre part  $X_*(T_{L_{\SCNX}})$ est facteur direct dans $X_*(T_{G_{\scnx}})$. 
\end{preuve}

Le lemme \ref{lem:DC} montre que la projection $\ago_T\to \ago_M$ envoie $X_*(T)$ surjectivement sur $X_*(M)$. D'autre part, elle envoie le sous-$\ZZ$-module $X_*(T_{G_{\scnx}})$ surjectivement sur le sous-$\ZZ$-module $X_*(M_{\scnx})$ de $X_*(M)$. 

\end{paragr}

\begin{paragr}[Le poids $\mathrm{w}^\xi_M$.] --- \label{S:wMxi} Pour tout $g\in G(\AAA)$, on introduit le poids 
$$\mathrm{w}^\xi_M( g)=|\{\mu \in X_*(M) \ | \ \mathbf{1}_{M,g}(\xi_M+\mu)=1\}|$$
qui est le nombre (fini) de points de l'ensemble $\xi_M+X_*(M)$ qui appartiennent à l'enveloppe convexe des points $-H_P(g)$ pour $P\in \pc(M)$.

\begin{lemme} \label{lem:R_M}  La fonction poids $g\in G(\AAA)\mapsto \mathrm{w}^\xi_M( g)$ est invariante à gauche par $M(\AAA)$.

Soit $h$ et $X$ deux éléments respectivement des ensembles décrits en 1 et 2  du lemme \ref{lem:systrep}. On a une suite exacte 
$$ 1 \longrightarrow T_X(\AAA_0)^1 \longrightarrow T_X(\AAA_0)\overset{H_M}{\longrightarrow} X_*(M) \longrightarrow 0.$$
\end{lemme}

\begin{preuve} L'action à gauche de $h\in M(\AAA)$ translate l'enveloppe conexe des points  $-H_P(g)$ par le vecteur $H_M(h)$. Comme ce dernier appartient à $X_*(M)$, cela ne change pas $\mathrm{w}^\xi_M( g)$.

 Dans la suite exacte, seule la surjectivité à droite n'est pas évidente. Il suffit de prouver cette surjectivité pour la restriction de $H_M$ à $T_X(\Fq((z_\infty)))$. En combinant le lemme \ref{lem:cjsep} et la proposition \ref{prop:caractdesAM}, on montre que $X$ et $t_a$ sont conjugués sous $M((z_\infty))(k)$ ($t_a\in \tgo^{\reg}[[z_\infty]](k)$ est l'élément défini dans  la proposition \ref{prop:caractdesAM}). On remarquera qu'on a même $t_a\in \tgo^{\reg}[[z_\infty]](\Fq)$. L'automorphisme $\Int(m)$ induit alors un $\Fq((z_\infty))$-isomorphisme entre les $\Fq((z_\infty))$-tores déduits de $T_X$ et $T$ par changement de base. La surjectivité est alors évidente. 
\end{preuve}
\end{paragr}

\begin{paragr}[Intégrales orbitales pondérées $J_M^\xi$.] ---  \label{S:IOPxi}  Soit $h$ et $X$ deux éléments respectivement des ensembles décrits en 1 et 2  du lemme \ref{lem:systrep}. Posons $Y=\Ad(h^{-1})X$ et $T_Y(\AAA_0)=\Int(h^{-1}) T_X(\AAA_0)$. Le groupe $T_X(\AAA_0)$ a été muni d'une  mesure de Haar en \S\ref{S:Haar}. Le groupe  $T_Y(\AAA_0)$ est muni de la mesure déduite par transport. On introduit alors l'intégrale orbitale pondérée suivante 
 \begin{equation}
   \label{eq:IOPw}
   J^\xi_M(Y, \mathbf{1}_{D_0} )=\int_{T_Y(\AAA_0)\back G(\AAA_0)}   \mathbf{1}_{D_0}(\Ad(g^{-1}) Y) \mathrm{w}^\xi_M( g) dg.
 \end{equation}
 On a muni $T_Y(\AAA_0)\back G(\AAA_0)$ de la mesure quotient. L'intégrale converge car l'intégrande est à support compact. Pour alléger les notations, on omet la fonction  $\mathbf{1}_{D_0}$ et on pose
$$  J^\xi_M(Y)=J^\xi_M(Y, \mathbf{1}_{D_0} ).$$

\begin{remarque}
  Cette intégrale orbitale pondérée n'est pas celle qu'Arthur utilise habituellement. Le lien avec celle d'Arthur se fera au \S\ref{S:comparaisonIOP}.
\end{remarque}
\end{paragr}

\begin{paragr}[Un premier comptage.] --- Voici une première expression pour le cardinal d'une fibre.  
\begin{proposition}\label{prop:1ercomptage}
  On a l'égalité
$$\card([M(F)\back \xc_{(a,t)}^{\xi,\tau}])=\vol(a,t) \sum_{ h} \, \sum_{X} \, J^\xi_M(\Ad(h^{-1})X)$$
où $h$ et $X$ parcourent respectivement les ensembles décrits en 1 et 2 du lemme \ref{lem:systrep} .
\end{proposition}

\begin{preuve}
  D'après la formule donnée au \S\ref{S:cardinal} et le lemme \ref{lem:systrep}, le cardinal du groupoïde $[M(F)\back \xc_{(a,t)}^{\xi,\tau}]$ s'écrit comme la somme sur les éléments $h$ et $X$ des ensembles respectivement décrits en 1 et 2 du lemme \ref{lem:systrep} de 
$$ \sum_{g\in T_Y(F_0) \back G(\AAA_0) / K} \frac{1}{|T_Y(F_0)\cap  g Kg^{-1} |}\mathbf{1}_{D_0}(\Ad(g^{-1}) Y) \mathbf{1}_{M,g}(\xi_M+H_M(h)), $$
où l'on a posé  $Y=\Ad(h^{-1})X$ et $T_Y(F_0)=\Int(h^{-1}) T_X(F_0)$ comme ci-dessus.  D'après notre choix de mesure sur $G(\AAA_0)$, cette expression s'écrit encore
$$\int_{  T_Y(F_0) \back G(\AAA_0) }  \mathbf{1}_{D_0}(\Ad(g^{-1}) Y) \mathbf{1}_{M,g}(\xi_M+H_M(h)) \,  dg.$$
Introduisons le groupe  $T_Y(\AAA_0)=\Int(h^{-1}) T_X(\AAA_0)$ qui est muni de la mesure de Haar déduite de  $T_X(\AAA_0)$. Comme $g\mapsto   \mathbf{1}_{D_0}(\Ad(g^{-1}) Y)$ est invariante à droite par $T_Y(\AAA_0)$, l'expression précédente s'écrit 
$$\int_{  T_Y(\AAA_0) \back G(\AAA_0) }  \mathbf{1}_{D_0}(\Ad(g^{-1}) Y) \int_{T_Y(F_0)\back  T_Y(\AAA_0)} \mathbf{1}_{t g }(\xi_M+H_M(h)) \,dt \, \frac{dg}{dt}.$$
 En utilisant l'égalité 
$$ \mathbf{1}_{t g }(\xi_M+H_M(h))= \mathbf{1}_{g }(\xi_M+H_M(h)+H_{M}(t))$$
nos choix de mesures et le lemme \ref{lem:R_M}, on voit qu'on a
$$\int_{T_Y(F_0)\back  T_Y(\AAA_0)} \mathbf{1}_{t g }(\xi_M+H_M(h)) \,dt=\vol(T_{X_0}(F_0)\back T_{X_0}(\AAA_0)^1) \cdot \mathrm{w}^\xi_M( g).$$
Le résultat s'en déduit.
\end{preuve}
  
\end{paragr}

\begin{paragr}[Une formule pour le poids $\mathrm{w}_M^\xi$.] --- Le but de cette section est de donner une formule analytique pour le poids $\mathrm{w}^\xi_M$. Cela nous permettra ensuite de comparer ce poids à celui qu'Arthur considère dans ses travaux. Pour cela, on suit, à peu de choses près, Arthur (cf. \cite{localtrace} \S6). On a introduit au \S\ref{S:reseaux} le sous-$\ZZ$-module $X_*(M_{\scnx})$ de $X_*(M)$. Soit $P\in \pc(M)$. L'ensemble $\Delta_P^\vee$ défini au \S\ref{S:deltachapeau} fournit une base de $X_*(M_{\scnx})$. Tout $\la\in \ago_M$  s'écrit de manière unique
  \begin{equation}
    \label{eq:la=laP}
    \la=[\la]_P + \{\la\}_P
  \end{equation}
  avec $[\la]_P\in X_*(M_{\scnx})$ et 
$$\{\la\}_P=\sum_{\al\in \Delta_P} r_\al \al^\vee$$
avec $0\leq r_\la<1$.
  
Pour tout $\Lambda\in \ago_M^*$, soit
$$\mathrm{c}_P(\Lambda)=\prod_{\al\in \Delta_P} (\exp(\Lambda(\al^\vee))-1).$$
Lorsque $M=G$, ce produit vaut $1$ (par convention tout produit sur l'ensemble vide vaut $1$). 
 Dans la suite, on dit que $\Lambda\in \ago_M^*$ est générique si $\mathrm{c}_P(\Lambda) \not=0$ pour tout $P\in \pc(M)$ c'est-à-dire $\Lambda(\al^\vee)\not=0$ pour tous $P\in \pc(M)$ et $\al^\vee\in \Delta_P^\vee$.

\begin{proposition} \label{prop:wM}Pour tout $g\in G(\AAA)$, on a l'égalité suivante
$$\mathrm{w}^\xi_M(g) =\lim_{\Lambda \to 0} \sum_{\mu_0\in X_*(M)/X_*(M_{\scnx}) } \sum_{P\in \pc(M)} \mathrm{c}_P(\Lambda)^{-1}\exp(-\Lambda(H_P(g) + [\mu_0+\xi_M]_P))$$
où la limite est prise sur les $\Lambda\in \ago_M^*$ génériques.
\end{proposition}

\begin{remarque}
 Dans la limite, la somme intérieure dépend du choix d'un système de représentants de $X_*(M)/X_*(M_{\scnx}) $. Changer $\mu_0$ en $\mu_0+\mu$ avec $\mu\in X_*(M_{\scnx})$ multiplie la somme intérieure par $\exp(-\Lambda(\mu))$. Cela n'affecte donc pas la limite en $\Lambda=0$.
\end{remarque}

\medskip

\begin{preuve} Elle repose entièrement sur les méthodes d'Arthur. Donnons quelques indications pour la commodité du lecteur. Soit $\Lambda\in \ago_M^*$ générique et  $P\in \pc(M)$. Soit $(\varpi_{\al})_{\al\in \Delta_P}$ la base de $\ago_M$ duale de $\Delta_P^\vee$. Soit 
$$\Delta_P^\Lambda=\{\al\in \Delta_P \ | \ \Lambda(\al^\vee)<0\}$$
et $\varphi_P^\Lambda$ la fonction caractéristique des $\la\in \ago_M$ tels que $\varpi_\al(\la)>0$ pour tout $\al\in \Delta_P^\Lambda$ et $\varpi_\al(\la)\leq 0$ pour $\al\in \Delta_P - \Delta_P^\Lambda$. 
Soit $g\in G(\AAA)$. D'après un lemme dû à Langlands (cf. \cite{localtrace} formule (3.8) p.22), la fonction caractéristique de l'enveloppe convexe des points $-H_P(g)$ est égale à la fonction
$$\sum_{P\in \pc(M)} (-1)^{|\Delta_P^\Lambda |} \varphi_P^\Lambda(\mu +H_P(g))$$
de la variable $\mu\in \ago_M$. Cette formule vaut dans notre contexte car la famille $(-H_P(g))_{P\in \pc(M)}$ est orthogonale positive au sens d'Arthur (cf. \cite{localtrace} pp.19-20). Il s'ensuit qu'on a
$$\mathrm{w}^\xi_M(g) = \sum_{\mu\in X_*(M)} \sum_{P\in \pc(M)} (-1)^{|\Delta_P^\Lambda |} \varphi_P^\Lambda(\mu +\xi_M+H_P(g)).$$

À la suite d'Arthur (cf. \cite{localtrace} \S6), on introduit la série de Fourier  pour $p\in \pc(M)$
 $$S_P(\Lambda)=\sum_{\mu\in X_*(M)} \varphi_P^\Lambda(\mu +\xi_M+H_P(g)) \exp(\Lambda(\mu)).$$
Cette série est absolument convergente et définit une fonction continue sur le complémentaire dans $\ago_M^*$ des hyperplans d'équation $\al^\vee$ pour $\al\in \Delta_P$. De plus, la limite en $\Lambda=0$ de cette série et la somme
$$\sum_{P\in \pc(M)} (-1)^{|\Delta_P^\Lambda |}S_P(\Lambda)$$
a pour limite en $\Lambda=0$ (où la limite est prise sur les $\Lambda$ génériques) le poids $\mathrm{w}^\xi_M(g)$. 
Pour calculer la somme de la série $S_P(\Lambda)$, il est plus commode de sommer d'abord sur $X_*(M_{\scnx})$ et ensuite sur un système de représentants de  $X_*(M)/X_*(M_{\scnx})$. Soit  $\mu_0\in X_*(M)/X_*(M_{\scnx})$. La contribution de $\mu_0$ à $S_P(\Lambda)$ est 
$$\sum_{\mu\in X_*(M_{\scnx})} \varphi_P^\Lambda(\mu +\mu_0+\xi_M+H_P(g)) \exp(\Lambda(\mu+\mu_0)).$$
Par un changement de variable et l'utilisation de la décomposition
$$\mu_0+\xi_M=[\mu_0+\xi_M]_P+\{\mu_0+\xi_M\}_P,$$
on voit que la somme ci-dessus est égale à 
$$\exp(\Lambda(\mu_0-[\mu_0+\xi]_P-H_P(g)))\sum_{\mu\in X_*(M_{\scnx})}  \varphi_P^\Lambda(\mu +\{\mu_0+\xi_M\}_P) \exp(\Lambda(\mu)).$$

Comme on a $\{\mu_0+\xi_M\}_P=\displaystyle \sum_{\al\in \Delta_P}r_\al \al^\vee$ avec $0\leq r_\al <1$ et que $\Delta_P^\vee$ est une base de $X_*(M_{\scnx})$, la somme ci-dessus est égale à 
$$\sum_{\ \ (m_\al)_{\al\in \Delta_P} }  \exp(\Lambda(\sum_{\al\in \Delta_P} m_\al \al^\vee))$$
 où l'on somme sur les familles d'entiers  $(m_\al)_{\al\in \Delta_P}$ qui vérifient $m_\al\geq 0$ pour $\al\in \Delta_P^\Lambda$ et $m_\al\leq -1$ pour $\al\in \Delta_P-\Delta_P^\Lambda$. Ces séries géométriques ont pour somme 
$$\sum_{\ \ (m_\al)_{\al\in \Delta_P} }  \exp(\Lambda(\sum_{\al\in \Delta_P} m_\al \al^\vee))= (-1)^{|\Delta_P^\Lambda|} \mathrm{c}_P(\Lambda)^{-1}.$$
Le résultat annoncé s'en déduit.
\end{preuve}
\end{paragr}

\begin{paragr}[Poids et intégrales orbitales pondérées d'Arthur.] --- \label{S:poidsArthur} Comme au \S\ref{S:HN}, on munit $\ago_T$ d'un produit scalaire invariant par $W$. Pour tout $L\in \lc$, on munit $\ago_L$ de la mesure de Haar qui donne le covolume $1$ aux réseaux engendrés par des bases orthonormales dans $\ago_L$. Par définition, pour tout $g\in G(\AAA)$, le poids d'Arthur 
$$\mathrm{v}_M(g)$$
est le volume dans $\ago_M$ de l'enveloppe convexe des points $-H_P(g)$. On va rappeler l'expression analytique due à Arthur de ce poids. Auparavant, pour tous $P\in \pc(M)$ et $\Lambda\in \ago_M^*$, on introduit le polynôme
$$\mathrm{d}_P(\Lambda)=\vol(\ago_M/ X_*(M_{\scnx}))^{-1} \prod_{\al\in \Delta_P} \Lambda(\al^\vee).$$ 
Pour $M=G$, on a, par convention, $\vol(\ago_G/ X_*(G_{\scnx}))=1$ et le polynôme ci-dessus est le polynôme constant égal à $1$.

\begin{proposition}(Arthur \cite{dis_series} pp. 219-220) \label{prop:vM}Pour tout $g\in G(\AAA)$, on a l'égalité suivante
$$\mathrm{v}_M(g) =\lim_{\Lambda \to 0} \sum_{P\in \pc(M)} \mathrm{d}_P(\Lambda)^{-1}\exp(-\Lambda(H_P(g))$$
où la limite est prise sur les $\Lambda\in \ago_M^*$ génériques.
  
\end{proposition}

\begin{preuve}Le poids $\mathrm{v}_M(g)$ est égal la limite en $\Lambda=0$ de la transformée de Fourier 
$$\int_{\ago_M}  \sum_{P\in \pc(M)} (-1)^{|\Delta_P^\Lambda |} \varphi_P^\Lambda(\mu +H_P(g)) \exp(\Lambda(\mu)) \, d\mu.$$
Or celle-ci se calcule et vaut $\displaystyle \sum_{P\in \pc(M)} \mathrm{d}_P(\Lambda)^{-1}\exp(-\Lambda(H_P(g))$.
\end{preuve}

Reprenons les notations du \S\ref{S:IOPxi}. Lorsqu'on fait agit $T_Y(\AAA_0)$ par translation à gauche sur $G(\AAA_0)$,  l'enveloppe convexe des points $(-H_P(g))$ pour $P\in \pc(M)$ subit une translation par un vecteur de $H_M(T_Y(\AAA_0))$. Le volume $\mathrm{v}_M(g)$ est donc invariant à gauche par $T_Y(\AAA_0)$. L'intégrale orbitale pondérée d'Arthur (cf. \cite{ar1} \S8) est définie  par la formule

 \begin{equation}
   \label{eq:IOPv}
   J_M(Y)=J_M(Y,\mathbf{1}_{D_0})=\int_{T_Y(\AAA_0)\back G(\AAA_0)}   \mathbf{1}_{D_0}(\Ad(g^{-1}) Y) \mathrm{v}_M( g) dg.
 \end{equation}

Ce sont en fait des analogues pour les algèbres de Lie des intégrales orbitales pondérées d'Arthur. Ces intégrales apparaissent naturellement dans un analogue pour les algèbres de Lie de la formule des traces d'Arthur-Selberg (cf. \cite{PHC}). 
\end{paragr}

\begin{paragr}[Les $(G,M)$-familles.] --- Dans \cite{trace_inv}, Arthur a introduit la notion de $(G,M)$-famille. Une $(G,M)$-famille est une famille de fonctions $(\mathrm{b}_P)_{P\in \pc(M)}$ lisses sur $\ago_M^*$ et qui vérifient pour un couple $(P,P')$ de sous-groupes paraboliques adjacents la condition de ``recollement''
$$\mathrm{b}_P(\Lambda)=\mathrm{b}_{P'}(\Lambda)$$
sur l'hyperplan $\Lambda(\al^\vee)=0$ où $\al^\vee$ est l'unique élément de $\Delta_P^\vee\cap (-\Delta_{P'}^\vee)$.  On peut alors définir \emph{a priori} uniquement pour $\Lambda$ générique la fonction
\begin{equation}
  \label{eq:bM}
  \mathrm{b}_M(\Lambda)=\sum_{P\in \pc(M)} \mathrm{d}_P(\Lambda)^{-1} \mathrm{b}_P(\Lambda)
\end{equation}
et Arthur montre que cette fonction se prolonge en une fonction lisse sur $\ago_M^*$ (\cite{trace_inv}, lemme 6.2). On pose
$$\mathrm{b}_M=\mathrm{b}_M(0).$$
Les exemples de $(G,M)$-familles que l'on considérera sont les deux suivants. Le premier exemple, qui dépend de $g\in G(\AAA)$, est la famille $(\mathrm{v}_P(g))_{P\in \pc(M)}$ où, pour tout $P\in \pc(M)$, on introduit la fonction de la variable $\Lambda$
\begin{equation}
  \label{eq:vP}
\mathrm{v}_P(g,\Lambda)=\exp(-\Lambda(H_P(g))).
\end{equation}
La condition de recollement résulte du lemme \ref{lem:pos}. Le second exemple est la famille  $(\mathrm{w}_P(\mu))_{P\in \pc(M)}$ qui dépend de $\mu\in \ago_M$ et qui est définie par
\begin{equation}
  \label{eq:wP}
  \mathrm{w}_P(\mu,\Lambda)=\frac{\mathrm{d}_P(\Lambda)}{\mathrm{c}_P(\Lambda)} \exp(-\Lambda([\mu]_P).
\end{equation}
Le lecteur vérifiera immédiatement que le lemme ci-dessous implique que la famille précédente est une $(G,M)$-famille.

\begin{lemme}
  Soit $\mu\in \ago_M$. Soit $P$ et $P'$ deux sous-groupes paraboliques dans $\pc(M)$ adjacents et soit $\al^\vee$ l'unique élément de $\Delta_P^\vee\cap (-\Delta_{P'}^\vee)$. Pour tout $\Lambda\in \ago_M^*$ qui vérifie $\Lambda(\al^\vee)=0$, on a
$$\Lambda([\mu]_P)=\Lambda([\mu]_{P'}).$$
\end{lemme}

\begin{preuve}
  L'hyperplan de $\ago_M^*$ d'équation $\Lambda(\al^\vee)=0$ n'est autre que le sous-espace $\ago_Q^*$ où $Q$ est le plus petit sous-groupe parabolique de $G$ qui contient à la fois $P$ et $P'$. La projection $p_Q$ de $\ago_M$ sur $\ago_Q$ selon la décomposition $\ago_M=\ago_Q\oplus \ago_M^Q$ induit une bijection de $\Delta_P^\vee-\{\al^\vee\}$, resp.   $\Delta_{P'}^\vee-\{-\al^\vee\}$, sur $\Delta_Q^\vee$. Soit $\mu\in \ago_M$. Écrivons cet élément dans les bases $\Delta_P^\vee$ et  $\Delta_{P'}^\vee$
$$\mu=\sum_{\beta\in \Delta_P} x_\beta \beta^\vee= \sum_{\gamma\in \Delta_{P'}} y_\gamma \gamma^\vee.$$
Pour tous $\beta\in \Delta_P$ et $\gamma\in  \Delta_{P'}$ distincts de $\pm\al$  tels que $p_Q(\beta^\vee)=p_Q(\gamma^\vee)$ on donc $x_\beta=y_\gamma$. En utilisant la notation $[\cdot]$ pour la partie entière, on a donc 
$$p_Q([\mu]_P)=\sum_{\beta \in \Delta_P-\{\al\}}   [x_\beta] p_Q(\beta^\vee)= \sum_{\gamma\in \Delta_{P'}-\{-\al\}} [y_\gamma]p_Q(\gamma^\vee)=p_Q([\mu]_{P'})$$
d'où
$$\Lambda([\mu]_P)=\Lambda([\mu]_{P'})$$
pour tout $\Lambda\in \ago_Q^*$.
\end{preuve}
\end{paragr}

\begin{paragr}[Produit de deux $(G,M)$-familles.] --- \label{S:produit}Le produit de deux $(G,M)$-familles est évidemment une $(G,M)$-famille. On note $(\mathrm{v\cdot w})(g,\mu)$ le produit des familles $\mathrm{v}(g)$ et $\mathrm{w}(\mu)$ définies en (\ref{eq:vP}) et (\ref{eq:wP}) ci-dessus. Dans le lemme suivant, on reformule les propositions \ref{prop:wM} et \ref{prop:vM}.

  \begin{lemme}
    \label{lem:reformulation}
Pour tout $g\in g(\AAA)$, les poids $\mathrm{w}_M^\xi(g)$ et $\mathrm{v}(g)$ sont respectivement les valeurs en $\Lambda=0$ des fonctions lisses
$$\sum_{\mu\in X_*(M)/X_*(M_{\scnx})} (\mathrm{v\cdot w})_M(g,\mu+\xi_M,\Lambda)\cdot $$
et 
$$\mathrm{v}_M(g,\Lambda).$$
\end{lemme}

La clef pour notre problème de comparaison d'intégrales orbitales pondérées est une formule pour  le produit $(\mathrm{v\cdot w})_M$ due à Arthur. Avant de pouvoir l'énoncer dans le lemme \ref{lem:produit}, nous aurons besoin de quelques notations supplémentaires.

Soit  $\mathrm{b}$ une  $(G,M)$-famille et $L\in \lc(M)$. Pour tout $Q\in \pc(L)$, la fonction sur $\ago_L^*$ obtenue par restriction de  $\mathrm{b}_P$ ne dépend pas du choix du sous-groupe parabolique  $P\in \pc^Q(M)$. La famille $(\mathrm{b}_Q)_{Q\in \pc(L)}$ est une $(G,L)$-famille. Par la formule (\ref{eq:bM}) appliquée au sous-groupe de Levi $L$, on obtient une fonction lisse notée $ \mathrm{b}_L(\Lambda)$  sur $\ago_L^*$ et on pose
$$ \mathrm{b}_L =\mathrm{b}_L(0).$$

\begin{lemme}
  \label{lem:wL} Soit $L\in \lc(M)$. Pour tout $\mu\in \ago_M$  soit  $(\mathrm{w}_{Q}(\mu,\Lambda))_{Q\in \pc(L)}$ la  $(G,L)$-famille déduite comme ci-dessus de la $(G,M)$-famille $(\mathrm{w}_{P}(\mu,\Lambda))_{P\in \pc(M)}$ définie en (\ref{eq:wP}). Alors on a 
$$\sum_{\mu\in X_*(M)/X_*(M_{\scnx})}\mathrm{w}_{L}(\mu+\xi_M)=  \frac{\vol(\ago_L/X_*(L))}{\vol(\ago_M/X_*(M))}\cdot \mathrm{w}_{L}^\xi(1)$$
où 
\begin{itemize}
\item par convention, on pose $\vol(\ago_G/X_*(G))=1$ ;
\item $1$ est l'élément neutre de $G(\AAA_0)$ ;
\item le second membre est celui défini au \S \ref{S:wMxi} lorsqu'on remplace $M$ par $L$.
\end{itemize}
\end{lemme}

\begin{preuve}
  Le membre de gauche vaut par définition
$$\lim_{\Lambda \to 0} \sum_{\mu_0\in X_*(M)/X_*(M_{\scnx}) } \, \ \sum_{Q\in \pc(L)} \mathrm{d}_Q(\Lambda)^{-1}\cdot \frac{ \mathrm{d}_P(\Lambda)}{ \mathrm{c}_P(\Lambda)}\cdot \exp(-\Lambda([\mu_0+\xi_M]_P))$$
où
\begin{itemize}
\item la limite est prise sur les $\Lambda\in \ago_L^*$ génériques ;
\item pour chaque $Q\in \pc(L)$, on choisit $P\in \pc^Q(M)$ ;
\item la fraction  $\frac{ \mathrm{d}_P(\Lambda)}{ \mathrm{c}_P(\Lambda)}$ est bien définie sur les points génériques de $\ago_M^*$ et se prolonge en une fonction lisse sur $\ago_M^*$.
\end{itemize}
Pour de tels $\Lambda$, $Q$ et $P$, on vérifie les formules
$$\frac{ \mathrm{d}_P(\Lambda)}{ \mathrm{c}_P(\Lambda)}=  \frac{\vol(\ago_L/X_*(L_{\scnx}))}{\vol(\ago_M/X_*(M_{\scnx}))}\cdot \frac{ \mathrm{d}_Q(\Lambda)}{ \mathrm{c}_Q(\Lambda)}$$
et
$$[\mu_0+\xi_M]_P=[\mu_0'+\xi_L]_Q$$
où $\mu_0'$ est la projection de $\mu_0$ sur $\ago_L$. Le lemme \ref{lem:DC} implique que la projection de $\ago_M$ sur $\ago_L$ induit une surjection de $ X_*(M)/X_*(M_{\scnx})$ sur  $X_*(L)/X_*(L_{\scnx})$. L'ordre de son noyau se combine au rapport des covolumes des réseaux $X_*(L_{\scnx})$ et $X_*(M_{\scnx})$ pour donner le facteur
$$\frac{\vol(\ago_L/X_*(L))}{\vol(\ago_M/X_*(M))}$$
On conclut la démonstration en appliquant la proposition \ref{prop:wM} au sous-groupe de Levi $L$ et à $g=1$.

\end{preuve}

Soit $Q$ un sous-groupe parabolique de Levi $L$.  Arthur définit également une $(L,M)$-famille $\mathrm{b}^Q$ de la façon suivante : pour tout $P\in \pc^L(M)$, $\mathrm{b}^Q_P(\Lambda)$ est la fonction lisse sur $\ago_M^*$ définie par
$$\mathrm{b}^Q_P(\Lambda)=\mathrm{b}_{PN_Q}(\Lambda),$$
où le groupe $PN_Q$ est l'élément de $\pc(M)$ engendré par $P$ et $N_Q$.
On peut alors former une fonction lisse 
$$\mathrm{b}_M^Q(\Lambda)= \sum_{P\in \pc^L(M)}      \mathrm{d}_P^Q(\Lambda)^{-1} \mathrm{b}^Q_P(\Lambda)$$ 
en utilisant le coefficient
 $$\mathrm{d}_P^Q(\Lambda)=\vol(\ago_M^Q/ X_*(M'))^{-1} \prod_{\al\in \Delta_P^Q} \Lambda(\al^\vee)$$
où $\Delta_P^Q\subset \Delta_P$ est le sous-ensemble des racines dans $L$ et $M'$ est l'image réciproque de $M$ dans le revêtement simplement connexe du groupe dérivé de $L$. Le $\ZZ$-module  $ X_*(M')$ est alors le réseau de $\ago_M^Q$ engendré par $\Delta_P^Q$.

Notons le lemme suivant qui est dû à Arthur et dont on laisse la démonstration au lecteur.

\begin{lemme} \label{lem:LMfamille} Soit $L\in \lc(M)$ et $Q\in \pc(L)$. Soit $g\in G(\AAA_0)$.   

Soit $\mathrm{v}^Q(g)$ la $(L,M)$-famille déduite de  $\mathrm{v}(g)$. 

Pour tout $l\in L(\AAA_0)$, soit $\mathrm{v}^L(l)$ la $(L,M)$-famille définie par $\mathrm{v}^L_P(l,\Lambda)=\exp(-\Lambda(H_P(l))$.

On a l'égalité de $(L,M)$-familles
$$\mathrm{v}^Q(g)=\mathrm{v}^L(l_Q(g))$$
où $l_Q(g) \in L(\AAA_0)$ est donné par la décomposition d'Iwasawa
$$g\in l_Q(g)N_Q(\AAA_0) K.$$
En particulier, on a l'égalité
\begin{equation}
  \label{eq:vMQ=vML}
  \mathrm{v}^Q_M(g)=\mathrm{v}^L_M(l_Q(g))
\end{equation}
\end{lemme}

Arthur déduit  également d'une  $(G,M)$-famille $\mathrm{b}$ une fonction lisse sur $\ago_Q^*$ notée   $\mathrm{b}_Q'(\Lambda)$ (cf. \cite{trace_inv} \S6(6.3) et  lemme 6.1). Nous ne rappelerons pas sa définition ; seules les deux propriétés suivantes nous seront utiles.

\begin{lemme}(cf. \cite{trace_inv} lemme 6.3 et corollaire 6.4) \label{lem:produit}Soit  $\mathrm{b}$ une  $(G,M)$-famille. 
  \begin{enumerate}
  \item Pour tout $L\in \lc(M)$ et $\Lambda\in \ago_L^*$, on a
$$\mathrm{b}_L(\Lambda)=\sum_{Q\in \pc(L)} \mathrm{b}_Q'(\Lambda).$$
\item Soit   $\mathrm{e}$ une  $(G,M)$-famille et $\mathrm{e}\cdot\mathrm{b}$ le produit de  $\mathrm{e}$ et  $\mathrm{b}$. Pour tout $\Lambda\in \ago_M^*$, on a l'égalité 
$$( \mathrm{e}\cdot\mathrm{b})_M(\Lambda)=\sum_{Q\in \fc(M)} \mathrm{e}_M^Q(\Lambda) \mathrm{b}_Q'(\Lambda).$$
\end{enumerate}
\end{lemme}

\end{paragr}

\begin{paragr}[Comparaison d'intégrales orbitales pondérées.] --- \label{S:comparaisonIOP}Soit $L\in \lc(M)$ et $Q\in \pc(L)$. On reprend les notations du \S\ref{S:IOPxi}. Les ``poids'' $\mathrm{v}^Q_M$ et $\mathrm{v}^L_M$ déduits des $(L,M)$-familles éponymes permettent de définir les intégrales orbitales pondérées ci-dessus
$$ J^Q_M(Y, \mathbf{1}_{D_0} )=\int_{T_Y(\AAA_0)\back G(\AAA_0)}   \mathbf{1}_{D_0}(\Ad(g^{-1}) Y) \mathrm{v}^Q_M( g)\, dg.$$
et 
   $$J_M^L(Y ,\mathbf{1}_{D_0})=\int_{T_Y(\AAA_0)\back L(\AAA_0)}   \mathbf{1}_{D_0}(\Ad(l^{-1}) Y) \, \mathrm{v}_M^L( l) \,dl.$$

Pour $Q=G$ ou $L=G$, on retrouve l'intégrale $J_M$ définie en (\ref{eq:IOPv}). Les deux intégrales sont reliées par la formule de descente énoncée ci-dessous.

\begin{lemme}
  \label{lem:descente}
Avec les notations ci-dessous, on a 
$$ J^Q_M(Y, \mathbf{1}_{D_0} )= q^{\dim(\ngo_Q)\deg(D_0)} J_M^L(Y ,\mathbf{1}_{D_0})$$
\end{lemme}

\begin{preuve} Elle repose sur la formule (\ref{eq:vMQ=vML}) du lemme \ref{lem:LMfamille}, la décomposition d'Iwasawa $G(\AAA_0)=L(\AAA_0) N_Q(\AAA_0) K$, qui est compatible à nos choix de mesures ainsi que, pour $l\in L(\AAA_0)$,  sur le changement de variables $n \mapsto \Ad(ln)^{-1}Y-Y$ qui induit un isomorphisme entre les espaces mesurés $N_Q(\AAA_0)$ et $\ngo_Q(\AAA_0)$, ce dernier étant muni de la mesure de Haar qui donne le volume $1$ au produit $\prod_{v\in V_0}\ngo_Q[[z_v]](\Fq)$. Le facteur en $q$ provient de l'égalité
$$\int_{\ngo_Q(\AAA_0)}  \mathbf{1}_{D_0}(Z) \, dZ= q^{\dim(\ngo_Q)\deg(D_0)}.$$
  \end{preuve}

On peut alors énoncer le théorème suivant.

  \begin{theoreme} Avec les notations ci-dessus, on a 
$$J_M^\xi(Y)=\sum_{L\in \lc(M)} q^{1/2(\dim(G)-\dim(L))} \cdot \frac{\vol(\ago_L/X_*(L))}{\vol(\ago_M/X_*(M))}\cdot J_M^L(Y) \cdot \mathrm{w}_L^\xi(1)$$
où $1\in G(\AAA_0)$ est l'élément neutre.
\end{theoreme}

      \begin{preuve}
Soit $g\in G(\AAA_0)$ et $\mu\in \ago_M$. D'après l'assertion 2 du lemme  \ref{lem:produit}, on a 
$$ (\mathrm{v\cdot w})_M(g,\mu,\Lambda)=\sum_{Q\in \fc(M)} \mathrm{v}_M^Q(g,\Lambda) \mathrm{w}_Q'(\mu,\Lambda).$$

En utilisant cette formule et le lemme \ref{lem:reformulation}, on obtient
$$
  \mathrm{w}^\xi_M(g)=\sum_{Q\in \fc(M)} \mathrm{v}_M^Q(g) \cdot\big(\sum_{\mu\in  X_*(M)/X_*(M_{\scnx})} w_Q'(\mu+\xi_M)\big)$$
d'où l'on déduit la relation
\begin{equation}
  \label{eq:avecQ}
  J_M^\xi(Y)=\sum_{Q\in \lc(M)} J_M^Q(Y)  \cdot\big(\sum_{\mu\in  X_*(M)/X_*(M_{\scnx})} w_Q'(\mu+\xi_M)\big)
\end{equation}
En utilisant la formule de descente du lemme \ref{lem:descente}, la ligne (\ref{eq:avecQ}) ci-dessus se réécrit
\begin{equation}
  \label{eq:avecL}
  J_M^\xi(Y)=\sum_{L\in \lc(M)} q^{1/2(\dim(G)-\dim(L))} J_M^L(Y)  \cdot\big(\sum_{\mu\in  X_*(M)/X_*(M_{\scnx})} \sum_{Q\in \pc(L) }w_Q'(\mu+\xi_M)\big).
\end{equation}
D'après l'assertion 1 du lemme \ref{lem:produit}, on a 
$$\sum_{Q\in \pc(L) }w_Q'(\mu+\xi_M)=w_L(\mu+\xi_M)$$
et d'après le lemme \ref{lem:wL}, on a 
$$\sum_{\mu\in  X_*(M)/X_*(M_{\scnx})}w_L(\mu+\xi_M)=\frac{\vol(\ago_L/X_*(L))}{\vol(\ago_M/X_*(M))}\cdot\mathrm{w}_L^\xi(1).$$
Cela conclut la démonstration.
      \end{preuve}

      \begin{corollaire}
        \label{cor:comparaison}
Avec les notations ci-dessus, on a pour tout $\xi\in \ago_T$ \emph{en position générale} l'égalité
$$J_M^\xi(Y)= \vol(\ago_M/X_*(M))^{-1}\cdot  J_M^G(Y).$$
      \end{corollaire}
      
\begin{preuve} Soit $L\in \lc(M)$. Par définition (cf. \S\ref{S:wMxi}), le poids $\mathrm{w}_L^\xi(1)$ vaut $1$ si $0$ appartient à $\xi_L+X_*(L)$ et $0$ sinon. Or comme $\xi$ est en position générale, le premier cas n'est possible que si $L=G$ (cf. définition \ref{def:posgenerale}). D'où le corollaire.        
      \end{preuve}

\end{paragr}

\begin{paragr}[Démonstration du théorème \ref{thm:comptage}.] Le théorème  \ref{thm:comptage} est simplement une reformulation de la proposition \ref{prop:1ercomptage},  pour des $\xi$ en position générale, à l'aide du corollaire \ref{cor:comparaison}.

\end{paragr}

\bibliographystyle{plain}
\bibliography{biblio}

\bigskip

\begin{flushleft}
Pierre-Henri Chaudouard \\
CNRS et Universit\'{e} Paris-Sud \\
 UMR 8628 \\
 Math\'{e}matique, B\^{a}timent 425 \\
F-91405 Orsay Cedex \\
France \\
\bigskip

Gérard.Laumon \\
CNRS et Universit\'{e} Paris-Sud \\
 UMR 8628 \\
 Math\'{e}matique, B\^{a}timent 425 \\
F-91405 Orsay Cedex \\
France \\
\bigskip

Adresses électroniques :\\
Pierre-Henri.Chaudouard@math.u-psud.fr \\
Gerard.Laumon@math.u-psud.fr\\
\end{flushleft}
\end{document}